\documentclass[12pt]{article}
\usepackage{eurosym}
\usepackage{amsfonts}
\usepackage{amsmath,amssymb}
\usepackage[a4paper, total={7.0in, 9.6in}]{geometry}
\usepackage[utf8]{inputenc}
\usepackage[english]{babel}
\usepackage{amsthm}

\newtheorem{Lemma}{Lemma}[section]
\newtheorem{theorem}[Lemma]{Theorem}
\newtheorem{lemma}[Lemma]{Lemma}

\newtheorem{remark}[Lemma]{Remark}
\newtheorem{definition}[Lemma]{Definition}

\begin{document}

\begin{center}
\bigskip {\Large \textbf{Weak solution for Stochastic Degasperis-Procesi
Equation}} \vspace{3mm}\\[0pt]

\textsc{\ Nikolai V. Chemetov}{\footnote{%
Department of Computing and Mathematics, University of S{\~a}o Paulo,
14040-901 Ribeir{\~a}o Preto - SP, Brazil, E-mail: nvchemetov{\char'100}
gmail.com.}},

\textsc{Fernanda Cipriano}{\footnote{%
Center for Mathematics and Applications (NOVA Math) and Department of
Mathematics, NOVA SST, Lisbon, Portugal E-mail: cipriano{\char'100}%
fct.unl.pt.}}
\end{center}

\date{\today }

\begin{abstract}
This article is concerned with the existence of solution to the stochastic
Degasperis-Procesi equation on $\mathbb{R}$ with an infinite dimensional
multiplicative noise and integrable initial data.

Writing the equation as a system composed of a stochastic nonlinear
conservation law and an elliptic equation, we are able to develop a method
based on the conjugation of kinetic theory with stochastic compactness
arguments. More precisely, we 
apply the stochastic Jakubowski-Skorokhod representation theorem to show the
existence of a weak kinetic martingale solution. 
In this framework, the solution is a stochastic process with sample paths in
Lebesgue spaces, which are compatible with peakons and wave breaking
physical phenomenon.
\end{abstract}

\vspace{3mm} \textbf{Key words.} Stochastic Degasperis-Procesi equation,
kinetic method, solvability.\newline


\vspace{3mm}

\thispagestyle{empty}

\noindent 

\section{Introduction}

\setcounter{equation}{0} 
This work is devoted to the study of a multiplicative white noise
perturbation of the deterministic Degasperis-Procesi (DP) equation 
\begin{equation*}
\partial _{t}(\left( 1-\partial _{xx}^{2}\right) u)=\left( -4u\partial
_{x}u+3\partial _{x}u\partial _{xx}^{2}u+u\partial _{xxx}^{3}u\right)
\end{equation*}%
introduced by Degasperis, Processi in \cite{DP1998}. The DP equation has a
similar structure and shares analogous properties to the Camassa-Holm
shallow water wave equation. Namely, Degasperis, Holm, Hone in \cite{DHH}
proved the exact integrability of DP equation by constructing a Lax pair,
they obtained a bi-Hamiltonian structure and an infinite hierarchy of
symmetries and conservation laws. Moreover, they verified that the equation
admits exact peakon solutions analogous to the Camassa-Holm peakons. Later
on Lundmark, Szmigielski \cite{LS1,LS2} showed that the DP equation admits
global weak solutions representing a wave train of peakons.

The Cauchy problem was discussed for the first time by Yin \cite{Yin1}.
Considering initial conditions $u_{0}$ in $H^{s}(\mathbb{R}),$ $s>3/2$, he
proved the existence and uniqueness of local-in-time strong solutions, and
derived a blow-up result and the precise blow-up scenario for the equation.
The global-in-time existence, $L^{1}$-stability and uniqueness results for
weak solutions in $L^{1}(\mathbb{R})\cap BV(\mathbb{R})$ and in $L^{2}(%
\mathbb{R})\cap L^{4}(\mathbb{R})$ with Kruzkov's entropy condition were
obtained by Coclite, Karlsen in \cite{CK06} (see also Escher, Liu, Yin \cite%
{Li}), and the periodic case was analysed in \cite{CK15}. In the article 
\cite{CK07}, the same authors show that the uniqueness result can be
achieved by replacing the entropy condition by an Ole\u{\i}nik-type
estimate. Let us emphasize that an important feature of weak solutions is
that they can develop discontinuities in finite time, which is interely
compatible with the physical phenomenon of wave breaking, arising from the
collision of a peakon with an anti-peakon, where the solution of the DP
equations develop shock waves with jump discontinuities (see Escher, Liu, Z.
Yin \cite{E1, E2}, Himonas, Holliman, Grayshan \cite{HHH}, Khorbatlyab,
Molinet \cite{K}, Liu, Yin \cite{LiuYin}, Lundmark \cite{L} for a detailed
analysis). Regarding the global existence of strong solutions and global
weak solutions of DP equation, let us also mention the studies by Lin, Liu 
\cite{LinLiu} and Yin \cite{Yin2, Yin3}.

\bigskip

In this work we are concerned with an existence result for the stochastic
Degasperis-Processi (SDP) equation with multiplicative noise, compatible
with the the physical phenomenon of wave breaking. More precisely, we
consider

\begin{equation}
d\left( 1-\partial _{xx}^{2}\right) u=-\left( 4u\partial _{x}u-3\partial
_{x}u\partial _{xx}^{2}u-u\partial _{xxx}^{3}u\right) dt+\left( 1-\partial
_{xx}^{2}\right) \sigma (u)d\mathcal{W}(t).  \label{Hyp1}
\end{equation}%
Setting $p=\left( 1-\partial _{xx}^{2}\right) ^{-1}\frac{3u^{2}}{2}$ and
taking into account that 
\begin{equation*}
\left( 1-\partial _{xx}^{2}\right) \left[ u\partial _{x}u+\partial _{x}p%
\right] =4u\partial _{x}u-3\partial _{x}u\partial _{xx}^{2}u-u\partial
_{xxx}^{3}u,
\end{equation*}
we can apply the operator $\left( 1-\partial _{xx}^{2}\right) ^{-1}$ to %
\eqref{Hyp1}, therefore the equation \eqref{Hyp1} can be formally written as 
\begin{equation}
\left\{ 
\begin{array}{l}
du=-\left( u\partial _{x}u+\partial _{x}\mathrm{p}\right) dt+\sigma (u)d{%
\mathcal{W}}_{t}, \\ 
\left( 1-\partial _{xx}^{2}\right) \mathrm{p}=\frac{3}{2}u^{2}.\vspace{2mm}%
\end{array}%
\right.  \label{Hyp2S}
\end{equation}%
This system being constituted by a stochastic nonlinear conservation law,
and an elliptic equation for the so-called pressure $p$ (by analogy with
Euler's equations), requires a subtil analysis and just few works can be
found in the literature. The SDP system \eqref{Hyp2S} with a linear
multiplicative noise $\sigma (u)=u$ has been considered by {\ Chen, Gao in
the article \cite{CG15},} where the authors studied the global
well-posedness in the space of stochastic processes with sample paths in $%
C([0,\infty ),H^{s}(\mathbb{R}))$, $s>3/2,$ however this regularity of the
sample paths are not suitable to model peakons and wave breaking. Recently,
Arruda, Chemetov, Cipriano \cite{ACC21} addressed the additive case ($\sigma
(u)=1$) in more irregular framework. Namely, they proved the existence of a
pathwise solution with values in a Lebesgue space $L^{p}$, by applying a
strategy based on the kinetic theory conjugated with the pathwise
deterministic method developed by Bessaih, Flandoli in \cite{Fla99} to study
the stochastic Euler equation.

{\ The purpose of this work is to show the existence 
of the solution to the SDP system \eqref{Hyp2S} with multiplicative noise,
for an initial condition in a Lebesgue space $L^{2}\cap L^{2q}$ for some $%
q>1 $. The solution will be a stochastic process, with integrable
trajectories suitable for modeling the peakons and wave breaking. Comparing
with previous approaches to the SDP equation, here we propose a different
method. Essentially, we 
introduce the notion of weak martingale kinetic solution, which is defined
at a suitable probability space by applying the Jakubowski-Skorokhod
representation theorem. 
One of the main dificulties we face reside on the fact that instead of just
one stochastic differential equation, we have a hyperbolic-elliptical
system, and the additional elliptic equation substantially complicates the
analysis. In fact, we are able to prove the existence of solution,
nevertheless due to the presence of the pressure term the uniqueness of the
solution for this system remains an open problem. Another dificulty is
related to the presence of the stochastic noise, which does not allow to
derive an $L^\infty-$ estimate for $u$ on the stochastic parameter. As a
consequence we have to use the more sophisticated kinetic method developed
by Lions, Perthame, Tadmor \cite{LPT94} to study deterministic conservation
law equations, instead of the standard Kruzkov approach. We refer Debussche,
Vovelle \cite{DV10}, and Breit, Feireisl, Hofmanov \cite{BEF18} for the
extension of the methods in \cite{LPT94} to the case of a stochastic
differential equation, and a clear exposition of the stochastic setting. }

Let us mention that the kinetic theory for various deterministic ($\sigma =0$%
) hyperbolic-elliptic systems in bounded and unbounded domains was developed
in \cite{chem3}-\cite{chem1}, \cite{per}.

The structure of the article is as follows. Section \ref{Sec2} formulates
the problem, presents the notion of solution and enunciates the main
results. In Section \ref{Sec3}, we introduce the regularized viscosity
system and derive uniform estimates with respect to the viscosity parameter.
Section \ref{Sec4} is devoted to the reformulation of the viscous system in
terms of viscous kinetic transport equations, and the deduction of uniform
estimates for the corresponding solutions. Taking into account the uniform
estimates obtained in the previous sections, Section \ref{Sec5} establishes
the tightness criterion, which will be used in Section \ref{Sec6} to apply
the Jakubowski-Skorokhod representation theorem to perform the limit
transition by the kinetic method. Finally, in the Appendix, we recall some
results known in the literature, necessary for our analysis.

.

\bigskip

\section{Formulation of the problem and main result}

\label{Sec2} \setcounter{equation}{0}

We introduce the appropriate functional spaces. Let $\mathbb{X}$ be some
topological space. We denote the space of bounded smooth functions from $%
\mathbb{X}$ to $\mathbb{R}$ by $C_{b}^{\infty }(\mathbb{X})$ and the
subspace of $C_{b}^{\infty }(\mathbb{X})$ {\ of functions }with a compact
support by $C_{c}^{\infty }(\mathbb{X}).$

We consider the space $L^{s}(\mathbb{R})$ for $s\geqslant 1$ and denote by $%
\Vert \cdot \Vert _{L^{s}}$ the corresponding norm. The usual inner product
in $L^{2}(\mathbb{R})$ and also in $L^{2}(\mathbb{R}^{2})$ are denoted by $%
(\cdot ,\cdot )$. {\ We hope that this abuse of notation does not bring
difficulties in reading, being clear within the context.}

Let $X$ be a real Banach space with norm $\left\Vert \cdot \right\Vert _{X}.$
We denote by $L^{r}(0,T;X)$ the space of $X$-valued measurable $r-$%
integrable functions defined on $[0,T]$ for $r\geqslant 1$.

Let $(\Omega ,\mathcal{F},P,\{\mathcal{F}_{t}\},\mathcal{W}_{t}^{k})$, $t\in
\lbrack 0,T]$, $k\in \mathbb{N}$, be a stochastic basis. In this article, we
consider a generic Hilbert space $H$, and a cylindrical Wiener process in $H$%
, which is written by 
\begin{equation*}
\mathcal{W}_{t}=\sum_{k}\mathcal{W}_{t}^{k}\text{e}_{k},
\end{equation*}%
where $\{\text{e}_{k}\}$ is an orthonormal basis of $H$, and $\mathcal{W}%
_{t}^{k},\;k\in \mathbb{N},$ are independent $\mathbb{R}-$valued Brownian
motions. We recall that $\mathcal{W}_{t}$ is a Wiener process
with values in a bigger space $\mathbb{H}$ such that the inclusion $%
H\hookrightarrow \mathbb{H}$ is a Hilbert-Schmidt operator.

This work is devoted to show the existence of the solution for the
stochastic Degasperis-Procesi equation with an infinite dimensional
multiplicative noise, in the one-dimensional domain $\mathbb{R} $. The
evolution system reads 
\begin{equation}
\left\{ 
\begin{array}{l}
du=-\left( u\partial _{x}u+\partial _{x}\mathrm{p}\right) \ dt+\sigma (u)d{%
\mathcal{W}}_{t},\vspace{2mm} \\ 
\left( 1-\partial _{xx}^{2}\right) \mathrm{p}=\frac{3}{2}u^{2}\qquad \qquad 
\text{in }\mathbb{R}_{T}=(0,T)\times \mathbb{R},\vspace{2mm} \\ 
u(0)=u_{0}\qquad \qquad \qquad \text{ in }\mathbb{R},%
\end{array}%
\right.  \label{Hyp3}
\end{equation}%
where $u=u(t,x)$ denotes the velocity of the fluid, $u_{0}=u_{0}(x)$ is the
initial velocity, and $\sigma $ is the diffusion coefficient.

We assume 
\begin{equation}
u_{0}\in L^{2}(\mathbb{R})\cap L^{2q}(\mathbb{R})\qquad \text{ for some \ }%
\,q>1.  \label{02}
\end{equation}

{Let us specify the diffusion coefficient $\sigma $. We consider a sequence
of continuous functions $\sigma _{k}:\mathbb{R}^{2}\rightarrow \mathbb{R}$, $%
k\in \mathbb{N}$, such that there exist some positive constants $C,$ $C_{0},$
satisfying 
\begin{eqnarray}
&&\sigma ^{2}(x,v):=\sum_{k}\sigma _{k}^{2}(x,v)\leqslant C_{0}|v|^{2},\quad
\forall x,v\in \mathbb{R};  \label{ass1} \\
&&D_{2s}z:=\sum_{k}\Vert \sigma _{k}(\cdot ,z(\cdot ))\Vert^2
_{L^{2s}}\leqslant C\Vert z\Vert^2 _{L^{2s}},\quad \forall z\in L^{2s}(%
\mathbb{R})\quad \text{for }s=1\text{ and }s=q;  \label{ass2} \\
&&\sum_{k}|\sigma _{k}(x,v)-\sigma _{k}(y,u)|^{2}\leqslant C\left(
|x-y|^{2}+|v-u|h(|x-y|)\right) .  \label{20:01:21:59}
\end{eqnarray}%
The last inequaility is valid for all $x,y,v,u\in \mathbb{R}$ and some
non-decreasing function $h:\mathbb{R}^{+}\rightarrow \mathbb{R}$, such that $%
h(0)=0.$ }

For any $u\in L^2(\mathbb{R}),$ we consider the linear mapping $\sigma
(u):H\rightarrow L^2(\mathbb{R})$ defined on the elements of the basis $%
\{e_k\}$ by $\sigma (u)($e$_{k})=\sigma _{k}(\cdot ,u)$, $k\in \mathbb{N}.$
Due to \eqref{ass1}, $\forall u\in L^2(\mathbb{R})$, the linear operator $%
\sigma (u)$ belongs to the space of Hilbert-Schmidt operators $HS(H;L^{2}(%
\mathbb{R}))$ defined in $H$ with values in $L^2(\mathbb{R})$. Namely, we
have 
\begin{eqnarray*}
\Vert \sigma (u)\Vert _{HS(H;L^{2}(\mathbb{R}))}^{2} :=\sum_{k}\Vert \sigma
_{k}(\cdot ,u(\cdot))\Vert _{L^{2}}^{2}\leqslant C_{0}\|u\|_{L^2}^{2}<\infty
.
\end{eqnarray*}

Furthermore, we introduce the differential operator 
\begin{equation}
A=4-\partial _{xx}^{2}.  \label{20:01:15:28}
\end{equation}%
Setting $z=A^{-1}v$ and 
\begin{equation}
S(z)=4\Vert z\Vert _{L^{2}}^{2}+5\Vert \partial _{x}z\Vert
_{L^{2}}^{2}+\Vert \partial _{xx}^{2}z\Vert _{L^{2}}^{2},
\label{17:01:17:22}
\end{equation}%
we recall from \cite{CK06} the validality of the following inequality 
\begin{equation}
C_{1}\Vert v\Vert _{L^{2}}^{2}\leqslant S(z)\leqslant C_{2}\Vert v\Vert
_{L^{2}}^{2},\qquad v\in L^{2}(\mathbb{R}).  \label{20:01:14:59}
\end{equation}%
We assume the existence of a positive constant $C$, such that 
\begin{eqnarray}
\sum_{k}\bigl(4\Vert A^{-1}\sigma _{k}(x,v(x))\Vert _{L^{2}}^{2} &+&5\Vert
\partial _{x}A^{-1}\sigma _{k}(x,v(x))\Vert _{L^{2}}^{2}  \notag \\
&+&\Vert \partial _{xx}^{2}A^{-1}\sigma _{k}(x,v(x))\Vert _{L^{2}}^{2}\bigr)%
\leqslant CS(z)  \label{17:01:17:33}
\end{eqnarray}%
for any $v\in L^{2}(\mathbb{R})$.

\begin{remark}
An example of our general setting \eqref{ass1}-\eqref{20:01:21:59} and %
\eqref{17:01:17:33}, corresponds to the case $H=\mathbb{R}$, and the
diffusion coefficient defined by%
\begin{equation*}
\sigma (v)=v.\qquad
\end{equation*}
\end{remark}

Let us introduce the notion of the stochastic weak kinetic solution.

\begin{definition}
\label{weak_kinetic} Assume the hypothesis \eqref{02}-\eqref{20:01:21:59}
and \eqref{17:01:17:33}. A system 
\begin{equation*}
( \widetilde{\Omega },\widetilde{\mathcal{F}},\widetilde{P},\{\widetilde{%
\mathcal{F}}_{t}\}_{t\in \lbrack 0,T]},\widetilde{\mathcal{W}},(\widetilde{u}%
,\widetilde{p}))
\end{equation*}%
is a stochastic weak kinetic solution of the system \eqref{Hyp3} if the
following properties hold:

\begin{enumerate}
\item $(\widetilde{\Omega},\widetilde{\mathcal{F}}, \widetilde{P})$ is a
complete probability space;

\item $\widetilde{\mathcal{W}}$ is a Brownian motion defined on the
probability space $(\widetilde{\Omega},\widetilde{\mathcal{F}}, \widetilde{P}%
)$ for the filtration $\{\widetilde{\mathcal{F}}_t\}_{t\in[0,T]}$;

\item $(\widetilde{u},\widetilde{p})=(\widetilde{u}(t),\widetilde{p}(t))$, $%
t\in \lbrack 0,T],$ is a predictable process defined on the filtered
probability space $(\widetilde{\Omega },\widetilde{\mathcal{F}},\widetilde{P}%
,\{\widetilde{\mathcal{F}}_{t}\}_{t\in \lbrack 0,T]})$, and 
there exist some positive constants $C_{s},$ such that for a.e. $t\in
\lbrack 0,T]$, 
\begin{equation*}
{\widetilde{\mathbb{E}}\Vert \widetilde{u}}(t)\Vert _{L^{2s}}^{2s}\leqslant
C_{s}\quad \text{for }s=1\ \text{and }s=q,
\end{equation*}
{where $\widetilde{\mathbb{E}}$ denotes the mathematical expectation with
respect to the measure $\widetilde{P}$},

\item there exists a positive measure $\widetilde{m}=\widetilde{m}(t,x,c)$
on $\mathbb{R}_{T}^{2},$ such that 
\begin{equation*}
{\widetilde{\mathbb{E}}}\int_{\mathbb{R}_{T}^{2}}|c|^{2(s-1)}\,d\widetilde{m}%
(t,x,c)<\infty ;
\end{equation*}

\item a.s. in $\widetilde{\Omega }$ \ the pair $(\widetilde{u},\widetilde{p}%
) $ verifies\ the system, consisting from\ the transport equality for the
kinetic function $\widetilde{f}(t,x,c)=sign^{+}(\widetilde{u}-c):$ 
\begin{eqnarray}
(\widetilde{f}(t),\varphi ) &=&\left( f_{0},\varphi (0)\right) +\int_{0}^{t}%
\left[ (\widetilde{f},c\partial _{x}\varphi )-(\partial _{x}\widetilde{p}\;%
\widetilde{f},\partial _{c}\varphi )\right] ds  \notag \\
&+&\sum_{k}\int_{0}^{t}\int_{\mathbb{R}^{2}}\varphi (s,x,\widetilde{u}%
(s,x))\sigma _{k}(x,\widetilde{u}(s,x))\,dxd\widetilde{{\mathcal{W}}}{%
_{s}^{k}}  \label{10:36:25:03:210} \\
&+&\frac{1}{2}\int_{0}^{t}\int_{\mathbb{R}^{2}}\partial _{c}\varphi (s,x,%
\widetilde{u})\sigma ^{2}(x,\widetilde{u})\,dxdt-\int_{0}^{t}\int_{\mathbb{R}%
^{2}}\partial _{c}\varphi \ d\widetilde{m}\quad \text{with }%
f_{0}=sign^{+}(u_{0}-c)   \notag
\end{eqnarray}%
for any  $\varphi =\varphi (t,x,c)\in C([0,T];C_{c}^{1}(\mathbb{%
R}^{2})),$ such that $\varphi |_{t=T}=0,$ and the  following equality
holds
\begin{equation}
\left( 1-\partial _{xx}^{2}\right) \widetilde{p}=\frac{3}{2}\widetilde{u}%
^{2}\qquad \text{ a.e. in }\mathbb{R}_{T}.  \label{10:36:25:03:21}
\end{equation}
\end{enumerate}
\end{definition}

Next, we formulate the main result of the article.

\begin{theorem}
\label{the_1} Under the hypothesis \eqref{02}-\eqref{20:01:21:59} and %
\eqref{17:01:17:33} there exists a stochastic weak kinetic solution 
\begin{equation*}
\left( \widetilde{\Omega },\widetilde{\mathcal{F}},\widetilde{P},\{%
\widetilde{\mathcal{F}}_{t}\}_{t\in \lbrack 0,T]},\widetilde{\mathcal{W}},(%
\widetilde{u},\widetilde{p})\right)
\end{equation*}%
of the system \eqref{Hyp3} in the sense of the Definition \ref{weak_kinetic}%
, satisfying also the entropy inequality a.s. in $\widetilde{\Omega }:$ 
\begin{eqnarray}
&&\int_{\mathbb{R}_{T}^{2}}\left\{ |\widetilde{u}-c|^{+}\partial _{t}\psi
+sign^{+}(\widetilde{u}-c)\left( \frac{\widetilde{u}^{2}-c^{2}}{2}\right)
\partial _{x}\psi dt-|\widetilde{u}-c|^{+}\partial _{x}\widetilde{p}\partial
_{c}\psi \right\} dtdxdc  \notag \\
&&+\int_{\mathbb{R}_{T}}|u_{0}-c|^{+}\psi (0)dxdc+\sum_{k}\int_{\mathbb{R}%
_{T}^{2}}sign^{+}(\widetilde{u}-c)\sigma _{k}(x,\widetilde{u})\psi \,dxdcd{%
\mathcal{\widetilde{W}}_{t}}  \notag \\
&&+\frac{1}{2}\int_{\mathbb{R}_{T}^{2}}sign^{+}(\widetilde{u}-c)\sigma
^{2}(x,\widetilde{u})\partial _{c}\psi \,dtdxdc\geqslant 0,\quad
\label{15:01:20:55}
\end{eqnarray}%
for any positive $\psi =\psi (t,x,c)\in C^{1}([0,T];C_{c}^{1}(%
\mathbb{R}^{2})),$ such that $\psi |_{t=T}=0.$
\end{theorem}

\section{Existence of the regularized problem and a priori estimates}

\label{Sec3} \setcounter{equation}{0}

In order to show the existence of the solution, we follow the viscosity
approach by introducing for any $\varepsilon \in (0,1)$ \ the following
stochastic {viscous} system 
\begin{equation}
\left\{ 
\begin{array}{l}
du_{\varepsilon }=\left[-\left( u_{\varepsilon }\partial _{x}u_{\varepsilon
}+\partial _{x}\mathrm{p}_{\varepsilon }\right) +\varepsilon \partial
_{xx}^{2}u_{\varepsilon }\right]\ dt+\sigma (u_{\varepsilon })d{\mathcal{W}}%
_{t}, \\ 
\\ 
\left( 1-\partial _{xx}^{2}\right) \mathrm{p}_{\varepsilon }=\frac{3}{2}%
u_{\varepsilon }^{2}\qquad \qquad \text{in }\mathbb{R}_{T}, \\ 
\\ 
u_{\varepsilon }(0)=u_{\varepsilon ,0}\qquad \qquad \qquad \quad \text{ in }%
\mathbb{R}.%
\end{array}%
\right.  \label{Sys}
\end{equation}%
Having (\ref{02}) and using the regularization procedure we can assume $%
u_{\varepsilon ,0}\in C_{c}^{\infty }(\mathbb{R}),$ such that 
\begin{eqnarray}
\Vert u_{\varepsilon ,0}\Vert _{L^{2s}(\mathbb{R})} &\mathbf{\leqslant }%
&C\Vert u_{0}\Vert _{L^{2s}(\mathbb{R})}\quad \text{for }s=1\text{ and }s=q,
\notag \\
u_{\varepsilon ,0} &\rightarrow &u_{0}\quad \text{strongly in }\quad L^{2}(%
\mathbb{R})\cap L^{2q}(\mathbb{R}).  \label{03}
\end{eqnarray}

{Applying, for instance, the arguments in \cite{CK}, \cite{DaPrato}, \cite%
{GR}, \cite{LS}, we can show the existence and the uniqueness of solution
for the viscous system.}

\begin{theorem}
\label{Visc_Sys} Under the hypothesis \eqref{02}-\eqref{20:01:21:59}, %
\eqref{17:01:17:33} and \eqref{03}    there exists a unique
distributional solution 
\begin{equation*}
u_{\varepsilon }\in L^{2}(\Omega ;\,C([0,T];L^{2}(\mathbb{R})\cap L^{2q}(%
\mathbb{R}))\cap L^{2}(0,T;H^{2}(\mathbb{R})))
\end{equation*}%
for the system (\ref{Sys}).
\end{theorem}

\medskip

\begin{remark}
\label{to zero} Due to the embedding $H^{2}(\mathbb{R})\hookrightarrow C(%
\mathbb{R}),$ the constructed solution $u_{\varepsilon }(t,\cdot )$ for the
system (\ref{Sys}) is a continuous function on the variable $x\in \mathbb{R}$%
; then the following property holds 
\begin{equation*}
u_{\varepsilon }(t,x)\rightarrow 0\text{\qquad as }|x|\rightarrow +\infty
\end{equation*}%
by the integrability $u_{\varepsilon }(t,\cdot )\in L^{2}(\mathbb{R})$\ \
for a.s. $(\omega ,t)\in \Omega \times (0,T).$
\end{remark}

\begin{remark}
\label{U0U} Taking into account that $u_{\varepsilon }\in C([0,T];L^{2}(%
\mathbb{R})\cap L^{2q}(\mathbb{R}))$ a.s. in $\Omega $ the first equation in %
\eqref{Sys} is understood in the following distributional sense 
\begin{eqnarray}
\int_{\mathbb{R}}u_{\varepsilon }(x,t)\varphi (x)dx &=&\int_{\mathbb{R}%
}u_{\varepsilon ,0}(x,0)\varphi (x)\,dx  \notag \\
&-&\int_{0}^{t}\int_{\mathbb{R}}\left( u_{\varepsilon }(x,s)\partial
_{x}u_{\varepsilon }(x,s)+\partial _{x}\mathrm{p}_{\varepsilon
}(x,s)-\varepsilon \partial _{xx}^{2}u_{\varepsilon }(x,s)\right) \varphi
(x)\,dxds  \notag \\
&+&\sum_{k}\int_{0}^{t}\int_{\mathbb{R}}\sigma _{k}(x,u_{\varepsilon
}(s,x))\varphi (x)\,dxd{\mathcal{W}_{s}^{k}},\quad \forall \varphi \in
C_{c}^{\infty }(\mathbb{R}).  \label{uuu}
\end{eqnarray}
\end{remark}

\subsection{$\protect\varepsilon -$Uniform estimates}

The following lemma collects the $\varepsilon -$uniform estimates for $%
(u_{\varepsilon },\mathrm{p}_{\varepsilon }).$

\begin{lemma}
\label{lem1}Under the hypothesis \eqref{02}-\eqref{20:01:21:59}, %
\eqref{17:01:17:33} and \eqref{03} the following estimates are valid for $%
s=1 $ and $s=q:$%
\begin{eqnarray}
\mathbb{E}\Vert u_{\varepsilon }\Vert _{L^{\infty }(0,T;L^{2}(\mathbb{R}%
))}^{2s}+\mathbb{E}\int_{\mathbb{R}_{T}}\varepsilon \left\vert
u_{\varepsilon }(t,x)\right\vert ^{2(s-1)}\left\vert \partial
_{x}u_{\varepsilon }(t,x)\right\vert ^{2}dxdt &\leqslant &C,  \notag \\
\mathbb{E}\Vert u_{\varepsilon }(t,x)\Vert _{L^{\infty }(0,T;L^{2q}(\mathbb{R%
}))}^{2q} &\leqslant &C  \label{Z}
\end{eqnarray}%
$\quad $ and%
\begin{eqnarray}
\mathbb{E}\Vert \mathrm{p}_{\varepsilon }\Vert _{L^{\infty }(0,T;W^{1,r}(%
\mathbb{R}))} &\leqslant &C,\qquad \forall r\in \lbrack 1,+\infty ],  \notag
\\
\mathbb{E}\Vert \partial _{xx}^{2}\mathrm{p}_{\varepsilon }\Vert _{L^{\infty
}(0,T;L^{q}(\mathbb{R}))}^{q} &\leqslant &C  \label{Z1}
\end{eqnarray}%
where $C$ are positive constants independent of $\varepsilon $.
\end{lemma}

\textbf{Proof. } We split the proof into three steps.

\textit{Step 1.} In this step we will prove the first estimate in \eqref{Z}.
Let $u_{\varepsilon }$ be the solution of the viscous system (\ref{Sys}). We
introduce the auxiliar stochastic process 
\begin{equation}
\ z_{\varepsilon }=A^{-1}u_{\varepsilon },  \label{y}
\end{equation}%
where $A$ is defined by \eqref{20:01:15:28}.  Let us show that $%
z_{\varepsilon }$ verifies the property 
\begin{equation}
z_{\varepsilon }(t,x)\rightarrow 0\text{\qquad as }|x|\rightarrow +\infty 
\text{\quad for a.s. }(\omega ,t)\in \Omega \times (0,T).  \label{zero}
\end{equation}%
Indeed, $z_{\varepsilon }$ has the explicit representation%
\begin{equation*}
z_{\varepsilon }(t,x)=\int_{\mathbb{R}}u_{\varepsilon }(t,y)\exp (-2|x-y|)\
dy,\qquad (t,x)\in \mathbb{R}_{T},
\end{equation*}%
and by Remark \ref{to zero} we have that for a.s. $(\omega ,t)\in \Omega
\times (0,T)$\ and $\forall \delta >0$ $\exists M=M(\delta )>0,$ such that $%
|u_{\varepsilon }(t,x)|<\frac{\delta }{2}$ for $|x|>M$. Denoting $%
N=max_{|x|\leqslant M}|u_{\varepsilon }(t,x)|,$ we obtain%
\begin{eqnarray*}
|z_{\varepsilon }(t,x)| &\leqslant &N\int_{|y|\leqslant M}\exp (-2|x-y|)\ dy+%
\frac{\delta }{2}\int_{\mathbb{R}}\exp (-2|x-y|)\ dy \\
&=&N~I(x)+\frac{\delta }{2}.
\end{eqnarray*}%
Since the integral $I(x)$ vanishes as $|x|\rightarrow +\infty ,$ there
exists $K>M,$ such that $|I(x)|<\frac{\delta }{2N}$ for $|x|>K,$ hence 
\begin{equation*}
|z_{\varepsilon }(t,x)|\leqslant \delta \text{\qquad as }|x|>K,
\end{equation*}%
implying \eqref{zero}. 

Applying the operator $\partial _{x}^{i}A^{-1}$ {\ (with $\partial
_{x}^{2}=\partial _{xx}$, $\partial _{x}^{1}=\partial _{x}$, $\partial
_{x}^{0}=$identity)}, for each $i=0,1,2,$ \ to the first equation in %
\eqref{Sys}, we obtain 
\begin{align}
d\partial _{x}^{i}z_{\varepsilon }(t,x)& =\partial _{x}^{i}A^{-1}\left(
-u_{\varepsilon }(t,x)\partial _{x}u_{\varepsilon }(t,x)-\partial
_{x}p_{\varepsilon }(t,x)+\varepsilon \partial _{xx}^{2}u_{\varepsilon
}(t,x)\right) \ dt  \notag \\
& +\sum_{k}\partial _{x}^{i}A^{-1}\sigma _{k}(x,u_{\varepsilon }(t,x))\,\ d{%
\mathcal{W}_{t}^{k}}\qquad \text{for each \ }i=0,1,2.  \notag
\end{align}%
For $S(z_{\varepsilon })$ defined by \eqref{17:01:17:22} the It\^{o}\
formula gives 
\begin{eqnarray}
\frac{1}{2}dS(z_{\varepsilon }) &=&4\left( dz_{\varepsilon },z_{\varepsilon
}\right) +5\left( d\partial _{x}z_{\varepsilon },\partial _{x}z_{\varepsilon
}\right) +\left( d\partial _{xx}^{2}z_{\varepsilon },\partial
_{xx}^{2}z_{\varepsilon }\right) \vspace{2mm}  \notag \\
&+&\frac{1}{2}\left[ 4\left( dz_{\varepsilon },dz_{\varepsilon })\right)
+5\left( d\partial _{x}z_{\varepsilon },d\partial _{x}z_{\varepsilon
}\right) +\left( d\partial _{xx}^{2}z_{\varepsilon },d\partial
_{xx}^{2}z_{\varepsilon }\right) \right] ,  \label{Hy21}
\end{eqnarray}%
where $\left( \varphi ,\psi \right) =$ $\int_{\mathbb{R}}\varphi (x)\psi
(x)\ dx$ is the the scalar product of two functions $\varphi ,\psi .$ \ 

On the other hand, integrating by parts we can verify that 
\begin{equation*}
\left( du_{\varepsilon },z_{\varepsilon }-\partial _{xx}^{2}z_{\varepsilon
}\right) =4\left( dz_{\varepsilon },z_{\varepsilon }\right) +5\left(
d\partial _{x}z_{\varepsilon },\partial _{x}z_{\varepsilon }\right) +\left(
d\partial _{xx}^{2}z_{\varepsilon },\partial _{xx}^{2}z_{\varepsilon
}\right) .
\end{equation*}%
Introducing this relation in \eqref{Hy21}, we obtain 
\begin{align*}
\frac{1}{2}dS(z_{\varepsilon })& =\left( du_{\varepsilon },z_{\varepsilon
}-\partial _{xx}^{2}z_{\varepsilon }\right) +\frac{1}{2}\left[ 4\left(
dz_{\varepsilon },dz_{\varepsilon })\right) +5\left( d\partial
_{x}z_{\varepsilon },d\partial _{x}z_{\varepsilon }\right) +\left( d\partial
_{xx}^{2}z_{\varepsilon },d\partial _{xx}^{2}z_{\varepsilon }\right) \right]
\\
& =\left\{ \left( -u_{\varepsilon }\partial _{x}u_{\varepsilon }-\partial
_{x}p_{\varepsilon },z_{\varepsilon }-\partial _{xx}^{2}z_{\varepsilon
}\right) +\left( \varepsilon \partial _{xx}^{2}u_{\varepsilon
},z_{\varepsilon }-\partial _{xx}^{2}z_{\varepsilon }\right) +\frac{1}{2}%
A_{0}(t)\right\} \ dt \\
& +\sum_{k}\int_{\mathbb{R}}\sigma _{k}(x,u_{\varepsilon }(x,t))\left(
z_{\varepsilon }-\partial _{xx}^{2}z_{\varepsilon }\right) \ dxd{\mathcal{W}%
_{t}^{k}}
\end{align*}%
with 
\begin{eqnarray*}
A_{0}(t) &=&\sum_{k}4\Vert A^{-1}\sigma _{k}(x,u_{\varepsilon }(x,t))\Vert
_{L^{2}}^{2}+5\Vert \partial _{x}A^{-1}\sigma _{k}(x,u_{\varepsilon
}(x,t))\Vert _{L^{2}}^{2} \\
&+&\Vert \partial _{xx}^{2}A^{-1}\sigma _{k}(x,u_{\varepsilon }(x,t))\Vert
_{L^{2}}^{2}\qquad \text{for }t\in (0,T).
\end{eqnarray*}

Again, integrating by parts and using the second equation of (\ref{Sys}), we
derive 
\begin{eqnarray*}
\left( -u_{\varepsilon }\partial _{x}u_{\varepsilon }-\partial
_{x}p_{\varepsilon },z_{\varepsilon }-\partial _{xx}^{2}z_{\varepsilon
}\right) &=&\left( -u_{\varepsilon }\partial _{x}u_{\varepsilon
},z_{\varepsilon }-\partial _{xx}^{2}z_{\varepsilon }\right) +\left(
-\partial _{x}(p_{\varepsilon }-\partial _{xx}^{2}p_{\varepsilon
}),z_{\varepsilon }\right) \\
&=&(-u_{\varepsilon }\partial _{x}u_{\varepsilon },4z_{\varepsilon
}-\partial _{xx}^{2}z_{\varepsilon })=(-u_{\varepsilon }\partial
_{x}u_{\varepsilon },u_{\varepsilon })=0
\end{eqnarray*}%
and 
\begin{equation}
\left( \partial _{xx}^{2}u_{\varepsilon },z_{\varepsilon }-\partial
_{xx}^{2}z_{\varepsilon }\right) (t)=-B_{0}(t)\qquad \text{for }t\in (0,T), 
\notag  \label{Hy23}
\end{equation}%
where $B_{0}(t)=S(\partial _{x}z_{\varepsilon })(t).$ Hence we have 
\begin{equation}
dS(z_{\varepsilon })+2\varepsilon B_{0}(t)\ dt=A_{0}(t)\ dt+2\sum_{k}\int_{%
\mathbb{R}}\sigma _{k}(x,u_{\varepsilon }(x,t))\left( z_{\varepsilon
}-\partial _{xx}^{2}z_{\varepsilon }\right) \ dxd{\mathcal{W}_{t}^{k}},
\label{23:37:20:01}
\end{equation}%
{which reads as 
\begin{equation*}
{S(z_{\varepsilon }(t))}+2\varepsilon \int_{0}^{t}B_{0}(\tau )\ d\tau
=\int_{0}^{t}A_{0}(\tau )\ d\tau +2\int_{0}^{t}\sum_{k}\int_{\mathbb{R}%
}\sigma _{k}(x,u_{\varepsilon }(x,\tau ))\left( z_{\varepsilon }-\partial
_{xx}^{2}z_{\varepsilon }\right) \ dxd{\mathcal{W}_{\tau }^{k}}.
\end{equation*}%
Let us introduce the following sequence of stopping times 
\begin{equation*}
\tau _{N}=\inf \{t\in \lbrack 0,T]:\,S(z_{\varepsilon }(t))\geqslant
N\}\wedge T,\qquad N\in \mathbb{N}.
\end{equation*}%
}

Taking the supremum with $s\in \lbrack 0,t\wedge \tau _{N}]$ and the
expectation, we obtain%
\begin{align*}
\mathbb{E}\sup_{s\in \lbrack 0,t\wedge \tau _{N}]}S(z_{\varepsilon }(s))&
+2\varepsilon \mathbb{E}\int_{0}^{t\wedge \tau _{N}}B_{0}(s)\ ds\leqslant 
\mathbb{E}S(z_{\varepsilon }(0))+\int_{0}^{t\wedge \tau _{N}}\mathbb{E}%
\sup_{\tau \in \lbrack 0,s]}A_{0}(\tau )\ ds \\
& +2\mathbb{E}\sup_{s\in \lbrack 0,t\wedge \tau
_{N}]}\int_{0}^{s}\sum_{k}\int_{\mathbb{R}}\sigma _{k}(x,u_{\varepsilon
}(x,\tau ))\left( z_{\varepsilon }-\partial _{xx}^{2}z_{\varepsilon }\right)
\ dxd{\mathcal{W}_{\tau }^{k}}.
\end{align*}%
{Recalling \eqref{ass2}, we have 
\begin{align*}
& \sum_{k}\left( \int_{\mathbb{R}}\sigma _{k}(x,u_{\varepsilon }(x,\tau
))\left( z_{\varepsilon }-\partial _{xx}^{2}z_{\varepsilon }\right) \
dx\right) ^{2} \\
& \leqslant \sum_{k}\Vert \sigma _{k}(x,u_{\varepsilon }(x,\tau ))\Vert
_{L^{2}}^{2}\Vert z_{\varepsilon }-\partial _{xx}^{2}z_{\varepsilon }\Vert
_{L^{2}}^{2}=D_{2}u_{\varepsilon }(\tau )\Vert z_{\varepsilon }-\partial
_{xx}^{2}z_{\varepsilon }\Vert _{L^{2}}^{2},
\end{align*}%
then the Burkholder-Davis-Gundy inequality gives%
\begin{align*}
\mathbb{E}\sup_{s\in \lbrack 0,t\wedge \tau _{N}]}&
\int_{0}^{s}\sum_{k}\int_{\mathbb{R}}\sigma _{k}(x,u_{\varepsilon }(x,\tau
))\left( z_{\varepsilon }-\partial _{xx}^{2}z_{\varepsilon }\right) \ dxd{%
\mathcal{W}_{\tau }^{k}} \\
& \leqslant \mathbb{E}\left( \int_{0}^{t\wedge \tau _{N}}\sum_{k}\left(
\int_{\mathbb{R}}\sigma _{k}(x,u_{\varepsilon }(x,\tau ))\left(
z_{\varepsilon }-\partial _{xx}^{2}z_{\varepsilon }\right) \ dx\right)
^{2}d\tau \right) ^{\frac{1}{2}} \\
& \leqslant \mathbb{E}\left( \int_{0}^{t\wedge \tau _{N}}\sum_{k}\Vert
\sigma _{k}(x,u_{\varepsilon }(x,\tau ))\Vert _{L^{2}}^{2}\Vert
z_{\varepsilon }-\partial _{xx}^{2}z_{\varepsilon }\Vert _{L^{2}}^{2}d\tau
\right) ^{\frac{1}{2}} \\
& \leqslant C\mathbb{E}\left( \int_{0}^{t\wedge \tau _{N}}\sup_{\tau \in
\lbrack 0,s]}D_{2}u_{\varepsilon }(\tau )\Vert z_{\varepsilon }-\partial
_{xx}^{2}z_{\varepsilon }\Vert _{L^{2}}^{2}\ d\tau \right) ^{\frac{1}{2}} \\
& \leqslant C\mathbb{E}\left( \int_{0}^{t\wedge \tau _{N}}\Vert
u_{\varepsilon }(\tau )\Vert _{2}^{2}S(z_{\varepsilon }(\tau ))\ d\tau
\right) ^{\frac{1}{2}}\leqslant \frac{1}{2}\mathbb{E}\sup_{s\in \lbrack
0,t\wedge \tau _{N}]}S(z_{\varepsilon }(s))+C\int_{0}^{t\wedge \tau _{N}}%
\mathbb{E}S(z_{\varepsilon }(\tau ))\ d\tau .
\end{align*}%
} Therefore, taking into account \eqref{17:01:17:33}, we deduce 
\begin{equation*}
\frac{1}{2}\mathbb{E}\sup_{s\in \lbrack 0,t\wedge \tau
_{N}]}S(z_{\varepsilon }(s))+2\varepsilon \mathbb{E}\int_{0}^{t\wedge \tau
_{N}}B_{0}(s)\ ds\leqslant \mathbb{E}S(z_{\varepsilon }(0))+C\int_{0}^{t}%
\mathbb{E}\sup_{\tau \in \lbrack 0,s\wedge \tau _{N}]}S(z_{\varepsilon
}(\tau ))\ ds.
\end{equation*}%
The Gronwall inequality and \eqref{02}, \eqref{03} imply 
\begin{equation}
\frac{1}{2}\mathbb{E}\sup_{s\in \lbrack 0,t\wedge \tau
_{N}]}S(z_{\varepsilon }(s))+2\varepsilon \mathbb{E}\int_{0}^{t\wedge \tau
_{N}}S(\partial _{x}z_{\varepsilon })\ ds\leqslant C\mathbb{E}%
S(z_{\varepsilon }(0)),\quad \forall N\in \mathbb{N},t\in \lbrack 0,T].
\label{13:07:25:004}
\end{equation}%
From this uniform estimate on $N$, we infer that $\tau _{N}\nearrow T$ in
probability, which gives the existence of a subsequence verifying $\tau
_{N_{k}}\nearrow T$ a.e. in $\Omega $. Writing \eqref{13:07:25:004} for the
subsequence $\{\tau _{N_{k}}\}$, passing to the limit as $k\rightarrow
+\infty $ and using the monotone convergence theorem, we deduce 
\begin{equation}
\frac{1}{2}\mathbb{E}\sup_{s\in \lbrack 0,t]}S(z_{\varepsilon
}(s))+2\varepsilon \mathbb{E}\int_{0}^{t}S(\partial _{x}z_{\varepsilon })\
ds\leqslant C\mathbb{E}S(z_{\varepsilon }(0)),\quad t\in \lbrack 0,T].
\label{13:07:25:04}
\end{equation}%
The inequality \eqref{20:01:14:59} and \eqref{02}, \eqref{03}\ imply the
first estimate of \eqref{Z} for the case $s=1$.

Now {\ applying} the It\^{o} formula to \eqref{23:37:20:01} with the
function $f(r)=r^{q}$, {$q\geqslant 2,$} we obtain 
\begin{align*}
dS(z_{\varepsilon })^{q}& =qS(z_{\varepsilon })^{q-1}dS(z_{\varepsilon })+%
\frac{q(q-1)}{2}S(z_{\varepsilon })^{q-2}\ dS(z_{\varepsilon
})dS(z_{\varepsilon }) \\
& =-2\varepsilon qS(z_{\varepsilon })^{q-1}B_{0}(t)\ dt+qS(z_{\varepsilon
})^{q-1}A_{0}(t)\ dt \\
& +2q(q-1)S(z_{\varepsilon })^{q-2}\sum_{k}\left( \int_{\mathbb{R}}\sigma
_{k}(x,u_{\varepsilon }(x,t))\left( z_{\varepsilon }-\partial
_{xx}^{2}z_{\varepsilon }\right) \ dx\right) ^{2}\ dt \\
& +2qS(z_{\varepsilon })^{q-1}\sum_{k}\int_{\mathbb{R}}\sigma
_{k}(x,u_{\varepsilon }(x,t))\left( z_{\varepsilon }-\partial
_{xx}^{2}z_{\varepsilon }\right) \ dxd{\mathcal{W}_{t}^{k}}.
\end{align*}%
Then%
\begin{align*}
dS(z_{\varepsilon })^{q}& +2\varepsilon qS(z_{\varepsilon })^{q-1}B_{0}(t)\
dt=qS(z_{\varepsilon })^{q-1}A_{0}(t)\ dt \\
& +2q(q-1)S(z_{\varepsilon })^{q-2}\sum_{k}\left( \int_{\mathbb{R}}\sigma
_{k}(x,u_{\varepsilon }(x,t))\left( z_{\varepsilon }-\partial
_{xx}^{2}z_{\varepsilon }\right) \ dx\right) ^{2}\ dt \\
& +2qS(z_{\varepsilon })^{q-1}\sum_{k}\int_{\mathbb{R}}\sigma
_{k}(x,u_{\varepsilon }(x,t))\left( z_{\varepsilon }-\partial
_{xx}^{2}z_{\varepsilon }\right) \ dxd{\mathcal{W}_{t}^{k}}.
\end{align*}%
Integrating on the interval $[0,s]$, for $0\leqslant s\leqslant t\wedge \tau
_{N}$, taking the supremum with $s\in \lbrack 0,t\wedge \tau _{N}]$ and the
expectation, we obtain 
\begin{align*}
& \mathbb{E}\sup_{s\in \lbrack 0,t\wedge \tau _{N}]}S(z_{\varepsilon
}(s))^{q}+2\varepsilon q\mathbb{E}\int_{0}^{t\wedge \tau
_{N}}S(z_{\varepsilon }(s))^{q-1}B_{0}(s)\ ds \\
& \leqslant \mathbb{E}S(z_{\varepsilon }(0))^{q}+q\int_{0}^{t\wedge \tau
_{N}}\mathbb{E}\sup_{\tau \in \lbrack 0,s]}\left[ S(z_{\varepsilon }(\tau
))^{q-1}A_{0}(\tau )\right] \ ds \\
& +2q(q-1)\int_{0}^{t\wedge \tau _{N}}\mathbb{E}\sup_{\tau \in \lbrack 0,s]}%
\biggl[S(z_{\varepsilon }(\tau ))^{q-2}\sum_{k}\left( \int_{\mathbb{R}%
}\sigma _{k}(x,u_{\varepsilon }(x,\tau ))\left( z_{\varepsilon }-\partial
_{xx}^{2}z_{\varepsilon }\right) \ dx\right) ^{2}\biggr]\ ds \\
& +2q\mathbb{E}\sup_{s\in \lbrack 0,t\wedge \tau
_{N}]}\int_{0}^{s}S(z_{\varepsilon }(\tau ))^{q-1}\sum_{k}\int_{\mathbb{R}%
}\sigma _{k}(x,u_{\varepsilon }(x,\tau ))\left( z_{\varepsilon }-\partial
_{xx}^{2}z_{\varepsilon }\right) \ \ dx\ d{\mathcal{W}_{\tau }^{k}} \\
& =I_{0}+I_{1}+I_{2}+I_{3}.
\end{align*}%
Due to \eqref{17:01:17:33}, we have%
\begin{equation*}
I_{1}=\int_{0}^{t\wedge \tau _{N}}\mathbb{E}\sup_{\tau \in \lbrack 0,s]}%
\left[ S(z_{\varepsilon }(\tau ))^{q-1}A_{0}(\tau )\right] \ ds\leqslant
C(q)\int_{0}^{t\wedge \tau _{N}}\mathbb{E}\sup_{\tau \in \lbrack
0,s]}S(z_{\varepsilon }(\tau ))^{q}\ ds
\end{equation*}%
for some constant $C(q)$ being independent of $\varepsilon $. Using %
\eqref{ass2}, we infer 
\begin{align*}
I_{2}& \leqslant C(q)\int_{0}^{t\wedge \tau _{N}}\mathbb{E}\sup_{\tau \in
\lbrack 0,s]}\biggl[S(z_{\varepsilon }(\tau ))^{q-2}\Vert z_{\varepsilon
}-\partial _{xx}^{2}z_{\varepsilon }\Vert _{L^{2}}^{2}\sum_{k}\Vert \sigma
_{k}(x,u_{\varepsilon }(x,\tau ))\Vert _{L^{2}}^{2}\ \biggr]\ ds \\
& \leqslant C\int_{0}^{t\wedge \tau _{N}}\mathbb{E}\sup_{\tau \in \lbrack
0,s]}\biggl[D_{2}u_{\varepsilon }(\tau )S(z_{\varepsilon }(\tau ))^{q-1}%
\biggr]\ ds \\
& \leqslant C\int_{0}^{t\wedge \tau _{N}}\mathbb{E}\sup_{\tau \in \lbrack
0,s]}\biggl[\Vert u_{\varepsilon }(\tau )\Vert _{L^{2}}^{2}S(z_{\varepsilon
}(\tau ))^{q-1}\biggr]\ ds \\
& \leqslant C\int_{0}^{t\wedge \tau _{N}}\mathbb{E}\sup_{\tau \in \lbrack
0,s]}S(z_{\varepsilon }(\tau ))^{q}\ ds.
\end{align*}%
The Burkholder-Davis-Gundy inequality gives%
\begin{align*}
I_{3}& =2q\mathbb{E}\sup_{s\in \lbrack 0,t\wedge \tau
_{N}]}\int_{0}^{s}S(z_{\varepsilon }(\tau ))^{q-1}\sum_{k}\int_{\mathbb{R}%
}\sigma _{k}(x,u_{\varepsilon }(x,\tau ))\left( z_{\varepsilon }-\partial
_{xx}^{2}z_{\varepsilon }\right) \ dxd{\mathcal{W}_{\tau }^{k}} \\
& \leqslant C(q)\mathbb{E}\left( \int_{0}^{t\wedge \tau _{N}}\sum_{k}\left(
S(z_{\varepsilon }(s))^{q-1}\int_{\mathbb{R}}\sigma _{k}(x,u_{\varepsilon
}(x,s))\left( z_{\varepsilon }-\partial _{xx}^{2}z_{\varepsilon }\right) \
dx\right) ^{2}\ ds\right) ^{\frac{1}{2}} \\
& \leqslant C(q)\mathbb{E}\left( \int_{0}^{t\wedge \tau
_{N}}S(z_{\varepsilon }(s))^{2q-2}D_{2}u_{\varepsilon }(s)\Vert
z_{\varepsilon }-\partial _{xx}^{2}z_{\varepsilon }\Vert _{L^{2}}^{2}\
dxds\right) ^{\frac{1}{2}} \\
& \leqslant C(q)\mathbb{E}\left( \int_{0}^{t\wedge \tau _{N}}\Vert
u_{\varepsilon }(\tau )\Vert _{L^{2}}^{2}S(z_{\varepsilon }(s))^{2q-1}\
ds\right) ^{\frac{1}{2}} \\
& \leqslant C(q)\mathbb{E}\left[ \sup_{s\in \lbrack 0,t\wedge \tau
_{N}]}\left( S(z_{\varepsilon }(s))\right) ^{\frac{q}{2}}\left(
\int_{0}^{t\wedge \tau _{N}}S(z_{\varepsilon }(s))^{q}\ ds\right) ^{\frac{1}{%
2}}\right] \\
& \leqslant \frac{1}{2}\mathbb{E}\sup_{s\in \lbrack 0,t\wedge \tau
_{N}]}S(z_{\varepsilon }(s))^{q}+C(q)\mathbb{E}\int_{0}^{t\wedge \tau
_{N}}S(z_{\varepsilon }(s))^{q}\ ds.
\end{align*}%
Therefore 
\begin{eqnarray*}
\frac{1}{2}\mathbb{E}\sup_{s\in \lbrack 0,t\wedge \tau
_{N}]}S(z_{\varepsilon }(t))^{q} &+&2\varepsilon q\mathbb{E}%
\int_{0}^{t\wedge \tau _{N}}S(z_{\varepsilon })^{q-1}B_{0}(s)\ ds \\
&\leqslant &\mathbb{E}S(z_{\varepsilon }(0))^{q}+C(q)\int_{0}^{t}\mathbb{E}%
\sup_{\tau \in \lbrack 0,s\wedge \tau _{N}]}S(z_{\varepsilon }(\tau ))^{%
\frac{q}{2}}\ ds.
\end{eqnarray*}%
The Gronwall inequality yields 
\begin{equation*}
\frac{1}{2}\mathbb{E}\sup_{s\in \lbrack 0,t\wedge \tau
_{N}]}S(z_{\varepsilon }(s))^{q}+2\varepsilon q\mathbb{E}\int_{0}^{t\wedge
\tau _{N}}S(z_{\varepsilon })^{q-1}S(\partial _{x}z_{\varepsilon })\
ds\leqslant C\mathbb{E}S(z_{\varepsilon }(0))^{q}.
\end{equation*}%
Passing to the limit (for a suitable subsequence) as $N\rightarrow +\infty $%
, we infer that 
\begin{equation*}
\frac{1}{2}\mathbb{E}\sup_{s\in \lbrack 0,t]}S(z_{\varepsilon
}(t))^{q}+2\varepsilon q\mathbb{E}\int_{0}^{t}S(z_{\varepsilon
})^{q-1}S(\partial _{x}z_{\varepsilon })\ ds\leqslant C\mathbb{E}%
S(z_{\varepsilon }(0))^{q}.
\end{equation*}%
Applying \eqref{02}, \eqref{20:01:14:59}, \eqref{03}, we derive the first
estimate of \eqref{Z} for the case {$s=q\geqslant 2$. The case $1<q<2$ can
be treated analogously, by first applying the It\^{o} formula to %
\eqref{23:37:20:01} with the function $f(r)=(r+\delta )^{q}$ for $\delta >0$
and then passing to the limit as $\delta \rightarrow 0.$}

\textit{Step 2.} Let us show the first estimate in \eqref{Z1} for the
function $\mathrm{p}_{\varepsilon }.$ Taking into account that $\mathrm{p}%
_{\varepsilon }$ has the following explicit representation%
\begin{equation*}
\mathrm{p}_{\varepsilon }(t,x)=\frac{3}{4}\int_{\mathbb{R}}u_{\varepsilon
}^{2}(t,y)\exp (-|x-y|)\ dy,\qquad (t,x)\in \mathbb{R}_{T},
\end{equation*}%
as the solution of the second equation of \eqref{Sys}, we derive that%
\begin{equation*}
\partial _{x}\mathrm{p}_{\varepsilon }(t,x)=-\frac{3}{4}\int_{\mathbb{R}%
}u_{\varepsilon }^{2}(t,y)\exp (-|x-y|)\mathrm{sign}(x-y)\ dy
\end{equation*}%
and%
\begin{equation}
\Vert \mathrm{p}_{\varepsilon }\Vert _{L^{\infty }},~\Vert \partial _{x}%
\mathrm{p}_{\varepsilon }\Vert _{L^{\infty }}\mathbf{\leqslant }C\int_{%
\mathbb{R}}u_{\varepsilon }^{2}(t,y)\ dy=C\Vert u_{\varepsilon }\Vert
_{L^{2}}^{2},\qquad \forall t\in \lbrack 0,T].  \label{10:24:21:01}
\end{equation}%
Moreover Fubini's theorem gives 
\begin{equation*}
\Vert \mathrm{p}_{\varepsilon }\Vert _{L^{1}},\Vert \partial _{x}\mathrm{p}%
_{\varepsilon }\Vert _{L^{1}}\leqslant C\int_{\mathbb{R}^{2}}(u_{\varepsilon
}(t,y)\exp (-|x-y|)\ dxdy\leqslant C\Vert u_{\varepsilon }\Vert _{L^{2}}^{2}.
\end{equation*}%
For any $r\in \lbrack 1,\infty )$ we have 
\begin{equation}
\Vert \mathrm{p}_{\varepsilon }\Vert _{L^{r}}\leqslant \left( \Vert \mathrm{p%
}_{\varepsilon }\Vert _{L^{\infty }}\right) ^{\frac{r-1}{r}}\left( \Vert 
\mathrm{p}_{\varepsilon }\Vert _{L^{1}}\right) ^{\frac{1}{r}}\leqslant
C(r)\Vert u_{\varepsilon }\Vert _{L^{2}}^{2},  \label{p1}
\end{equation}%
{then a similar argument gives } 
\begin{equation}
\Vert \partial _{x}\mathrm{p}_{\varepsilon }\Vert _{L^{r}}\leqslant
C(r)\Vert u_{\varepsilon }\Vert _{L^{2}}^{2}  \label{p2}
\end{equation}%
for constants $C(r)$, depending only on $r$. Combining \eqref{10:24:21:01}-%
\eqref{p2} and the first estimate of \eqref{Z} with $s=1$, we deduce %
\eqref{Z1}$_{1}$.

\textit{Step 3.} {\ Here, we show the second estimate in \eqref{Z} and the
second one in \eqref{Z1}. }

Applying the It\^{o} formula to (\ref{Sys})$_{1}$ for the function $%
f(r)=r^{2q}$, we obtain 
\begin{align}
d\left\vert u_{\varepsilon }(t,x)\right\vert ^{2q}& =2q(u_{\varepsilon
}(t,x))^{2q-1}\left( -u_{\varepsilon }(t,x)\partial _{x}u_{\varepsilon
}(t,x)-\partial _{x}\mathrm{p}_{\varepsilon }(t,x)+\varepsilon \partial
_{xx}^{2}u_{\varepsilon }(t,x)\right) \ dt  \notag \\
& +q(2q-1)\sum_{k}\sigma _{k}^{2}(t,x,u_{\varepsilon }(t,x))(u_{\varepsilon
}(t,x))^{2q-2}\ dt  \notag \\
& +2q\sum_{k}\sigma _{k}(x,u_{\varepsilon }(\tau ,x))(u_{\varepsilon
}(t,x))^{2q-1}\ d\mathcal{W}_{t}^{k}.  \label{26_11_55}
\end{align}%
Now, we argue as in the proof of Lemma 5.1 and Theorem 2.1 in 
\cite{K10}. First, we write the equatity \eqref{26_11_55} in the integral
form and integrate on $\mathbb{R}$. Next, the continuity and integrability
properties of $u_{\varepsilon }$ (see Theorem \ref{Visc_Sys} and Remark \ref%
{to zero}) allow to justify the application of the deterministic and
stochastic Fubini's theorems in order to deduce the following It\^{o}
formula on the Banach space $L^{2q}(\mathbb{R})$ 
\begin{align}
d\Vert u_{\varepsilon }(t,x)\Vert _{L^{2q}}^{2q}& =2q\int_{\mathbb{R}%
}(u_{\varepsilon }(t,x))^{2q-1}\left( -u_{\varepsilon }(t,x)\partial
_{x}u_{\varepsilon }(t,x)\right.  \notag \\
& \left. -\partial _{x}\mathrm{p}_{\varepsilon }(t,x)+\varepsilon \partial
_{xx}^{2}u_{\varepsilon }(t,x)\right) \ dxdt  \notag \\
& +q(2q-1)\sum_{k}\int_{\mathbb{R}}\sigma _{k}^{2}(x,u_{\varepsilon
}(t,x))(u_{\varepsilon }(t,x))^{2q-2}\ dxdt  \notag \\
& +2q\sum_{k}\int_{\mathbb{R}}\sigma _{k}(x,u_{\varepsilon
}(t,x))(u_{\varepsilon }(t,x))^{2q-1}\ dxd\mathcal{W}_{t}^{k}.  \label{u}
\end{align}%
We notice that 
\begin{equation*}
\int_{\mathbb{R}}(u_{\varepsilon }(t,x))^{2q-1}\left( u_{\varepsilon
}(t,x)\partial _{x}u_{\varepsilon }(t,x)\right) \ dx=\frac{1}{2q+1}\int_{%
\mathbb{R}}\partial _{x}(u_{\varepsilon }(t,x))^{2q+1}\ dx=0
\end{equation*}%
and 
\begin{equation*}
\varepsilon \int_{\mathbb{R}}(u_{\varepsilon }(t,x))^{2q-1}\partial
_{xx}^{2}u_{\varepsilon }(t,x)\ dx=-\varepsilon \left( 2q-1\right) \int_{%
\mathbb{R}}\left\vert u_{\varepsilon }(t,x)\right\vert ^{2(q-1)}\left\vert
\partial _{x}u_{\varepsilon }(t,x)\right\vert ^{2}\ dx.
\end{equation*}%
Now, let us consider the sequence of stopping times defined by 
\begin{equation*}
\tau _{N}=\inf \{t\in \lbrack 0,T]:\,\Vert u_{\varepsilon }(t)\Vert
_{L^{2q}}^{2q}\geqslant N\}\wedge T,\qquad N\in \mathbb{N}.
\end{equation*}

Hence integrating \eqref{u} on the interval $[0,s]$, for $0\leqslant
s\leqslant t\wedge \tau _{N}$, taking the supremum with $s\in \lbrack
0,t\wedge \tau _{N}]$ and the expectation, we obtain 
\begin{align*}
\mathbb{E}\sup_{s\in \lbrack 0,t\wedge \tau _{N}]}\Vert u_{\varepsilon
}(s,x)\Vert _{L^{2q}}^{2q}& +\varepsilon \left( 2q-1\right) \mathbb{E}%
\int_{0}^{t\wedge \tau _{N}}\int_{\mathbb{R}}\left\vert u_{\varepsilon
}(s,x)\right\vert ^{2(q-1)}\left\vert \partial _{x}u_{\varepsilon
}(s,x)\right\vert ^{2}\ dxds \\
& \leqslant \mathbb{E}\Vert u_{\varepsilon }(0,x)\Vert _{L^{2q}}^{2q} \\
& +2q\int_{0}^{t\wedge \tau _{N}}\mathbb{E}\sup_{\tau \in \lbrack 0,s]}\int_{%
\mathbb{R}}\left\vert u_{\varepsilon }(\tau ,x)\right\vert ^{2q-1}\left\vert
\partial _{x}\mathrm{p}_{\varepsilon }(\tau ,x)\right\vert \ dxds \\
& +q(2q-1)\int_{0}^{t\wedge \tau _{N}}\sum_{k}\int_{\mathbb{R}}\sigma
_{k}^{2}(x,u_{\varepsilon }(s,x))(u_{\varepsilon }(s,x))^{2q-2}\ dxds \\
& +2q\sup_{s\in \lbrack 0,t\wedge \tau _{N}]}\int_{0}^{s}\sum_{k}\mathbb{E}%
\int_{\mathbb{R}}\sigma _{k}(x,u_{\varepsilon }(\tau ,x))(u_{\varepsilon
}(\tau ,x))^{2q-1}\ dxd\mathcal{W}_{\tau }^{k} \\
& =J_{0}+J_{1}+J_{2}+J_{3}.
\end{align*}%
With the help of the H\"{o}lder inequality, the Young inequality, the
estimates \eqref{10:24:21:01} and \eqref{p2}, we derive 
\begin{equation*}
J_{1}\leqslant C(q)\Vert \partial _{x}\mathrm{p}_{\varepsilon }\Vert
_{L^{2q}}\Vert u_{\varepsilon }\Vert _{L^{2q}}^{2q-1}\leqslant C\Vert
u_{\varepsilon }\Vert _{L^{2}}\Vert u_{\varepsilon }\Vert
_{L^{2q}}^{2q-1}\leqslant C\left( \Vert u_{\varepsilon }\Vert
_{L^{2}}^{2q}+\Vert u_{\varepsilon }\Vert _{L^{2q}}^{2q}\right) .
\end{equation*}%
Also we have%
\begin{eqnarray*}
J_{2} &\leqslant &C(q)\int_{0}^{t\wedge \tau _{N}}\mathbb{E}\sup_{\tau \in
\lbrack 0,s]}\left( D_{2q}u_{\varepsilon }(\tau )\right) ^{1/2}\Vert
u_{\varepsilon }(\tau )\Vert _{L^{2q}}^{2(q-1)}\ ds \\
&\leqslant &C(q)\int_{0}^{t\wedge \tau _{N}}\mathbb{E}\sup_{\tau \in \lbrack
0,s]}\Vert u_{\varepsilon }(\tau )\Vert _{L^{2}}\Vert u_{\varepsilon }(\tau
)\Vert _{L^{2q}}^{2(q-1)}\ ds
\end{eqnarray*}%
In addition, the Burkholder-Davis-Gundy inequality and the assumption %
\eqref{ass2} give 
\begin{align*}
J_{3}& \leqslant C\mathbb{E}\sup_{s\in \lbrack 0,t\wedge \tau
_{N}]}\int_{0}^{s}\sum_{k}\int_{\mathbb{R}}\sigma _{k}(x,u_{\varepsilon
}(\tau ,x))(u_{\varepsilon }(\tau ,x))^{2q-1}\ dxd\mathcal{W}_{\tau }^{k} \\
& \leqslant \mathbb{E}\left( \int_{0}^{t\wedge \tau _{N}}\sum_{k}\left(
\int_{\mathbb{R}}\sigma _{k}(x,u_{\varepsilon }(\tau ,x))(u_{\varepsilon
}(\tau ,x))^{2q-1}\ dx\right) ^{2}\ d\tau \right) ^{\frac{1}{2}} \\
& \leqslant \mathbb{E}\left( \int_{0}^{t\wedge \tau
_{N}}D_{2q}u_{\varepsilon }(\tau )\Vert u_{\varepsilon }\Vert
_{L^{2q}}^{4q-2}\ d\tau \right) ^{\frac{1}{2}}\leqslant C\mathbb{E}\left(
\int_{0}^{t\wedge \tau _{N}}\Vert u_{\varepsilon }(\tau )\Vert
_{L^{2q}}^{2}\Vert u_{\varepsilon }\Vert _{L^{2q}}^{4q-2}\ d\tau \right) ^{%
\frac{1}{2}} \\
& \leqslant \frac{1}{2}\mathbb{E}\sup_{s\in \lbrack 0,t\wedge
\tau _{N}]}\Vert u_{\varepsilon }(s)\Vert _{L^{2q}}^{2q}+\mathbb{E}%
\int_{0}^{t\wedge \tau _{N}}\Vert u_{\varepsilon }(\tau )\Vert
_{L^{2q}}^{2q}\ d\tau  \\
& \leqslant \frac{1}{2}\mathbb{E}\sup_{s\in \lbrack 0,t\wedge
\tau _{N}]}\Vert u_{\varepsilon }(s)\Vert _{L^{2q}}^{2q}+\int_{0}^{t}{%
\mathbb{E}\sup_{\tau \in \lbrack 0,s\wedge \tau _{N}]}\Vert u_{\varepsilon
}(\tau )\Vert _{L^{2q}}^{2q}\ ds.}
\end{align*}%
Therefore 
\begin{align*}
\frac{1}{2}\mathbb{E}\sup_{s\in \lbrack 0,t\wedge \tau _{N}]}\Vert
u_{\varepsilon }(s)\Vert _{L^{2q}}^{2q}& +\varepsilon \left( 2q-1\right) 
\mathbb{E}\int_{0}^{t\wedge \tau _{N}}\int_{\mathbb{R}}\left\vert
u_{\varepsilon }(s,x)\right\vert ^{2(q-1)}\left\vert \partial
_{x}u_{\varepsilon }(s,x)\right\vert ^{2}\ dxds \\
& \leqslant \mathbb{E}\Vert u_{\varepsilon }(0)\Vert _{L^{2q}}^{2q}
+C\int_{0}^{t}{\mathbb{E}\sup_{\tau \in \lbrack 0,s\wedge \tau
_{N}]}\Vert u_{\varepsilon }(\tau )\Vert _{L^{2q}}^{2q}\ ds}.
\end{align*}%
Then the Gronwall inequality yields 
\begin{equation}
\mathbb{E}\sup_{s\in \lbrack 0,t\wedge \tau _{N}]}\Vert u_{\varepsilon
}(t,x)\Vert _{L^{2q}}^{2q}+\varepsilon \mathbb{E}\int_{0}^{t\wedge \tau
_{N}}\int_{\mathbb{R}}\left\vert u_{\varepsilon }(s,x)\right\vert
^{2(q-1)}\left\vert \partial _{x}u_{\varepsilon }(s,x)\right\vert ^{2}\
dxds\leqslant C,\quad \forall N\in \mathbb{N}.  \label{07_1}
\end{equation}%
Therefore, there exists a subsequence $\{\tau _{N_{k}}\}$ of $\{\tau _{N}\}$
such that $\tau _{N_{k}}\nearrow T$ a.e. in $\Omega $. Writing \eqref{07_1}
for the subsequence $\{\tau _{N_{k}}\}$ and passing to the limit as $%
k\rightarrow +\infty $, we infer that 
\begin{equation*}
\mathbb{E}\sup_{s\in \lbrack 0,t]}\Vert u_{\varepsilon }(t,x)\Vert
_{L^{2q}}^{2q}+\varepsilon \mathbb{E}\int_{0}^{t}\int_{\mathbb{R}}\left\vert
u_{\varepsilon }(s,x)\right\vert ^{2(q-1)}\left\vert \partial
_{x}u_{\varepsilon }(s,x)\right\vert ^{2}\ dxds\leqslant C,\quad t\in
\lbrack 0,T],
\end{equation*}%
by the monotone convergence theorem. Hence we deduce the second estimate in %
\eqref{Z}. Finally, the second equation of \eqref{Sys} implies the second
estimate in \eqref{Z1}.$\hfill \;\blacksquare $

\section{ Viscous kinetic transport equation}

\label{Sec4}

\setcounter{equation}{0}

In this section we reformulate the system (\ref{Sys}), showing that the pair 
$(u_{\varepsilon },p_{\epsilon })$ is the solution of a transport type
equation.

Let us introduce a so-called \textit{kinetic} function $f_{\varepsilon
}(t,x,c)=1_{u_{\varepsilon }(t,x)\geqslant c}$ and the measures 
\begin{equation}
\nu _{\varepsilon }^{(t,x)}(c)=\delta _{u_{\varepsilon }(t,x)=c}\text{\quad
and}\quad m_{\varepsilon }(t,x,c)=\varepsilon |\partial
_{x}u_{\varepsilon }(t,x)|^{2}\delta _{u_{\varepsilon }(t,x)=c}\mathcal{L}%
_{(t,x)} \text{\quad on}\quad (t,x,c)\in \mathbb{R}_{T}^{2},  \label{numu}
\end{equation}%
where $\delta _{u=c}$ is the Dirac measure at a point $u,$  $%
\mathcal{L}_{(t,x)}$\ is the Lebesgue measure on $(t,x)$ variables and $%
1_{u\geqslant c} $ is the characteristic function of the set $\left( -\infty
,u\right) .$

\begin{lemma}
\label{Lemma 4.4}1) There exist positive constants $C$, which are
independent of $\varepsilon $, such that%
\begin{eqnarray}
\mathbb{E}\sup_{t\in (0,T)}\int_{\mathbb{R}}(\int_{0}^{+\infty
}f_{\varepsilon }(t,x,c)p{c}^{p-1}dc &+&\int_{-\infty }^{0}(1-f_{\varepsilon
}(t,x,c))p|c|^{p-1}dc)dx\leqslant C,  \notag \\
\mathbb{E}\sup_{t\in (0,T)}\int_{\mathbb{R}}\left( \int_{0}^{+\infty
}f_{\varepsilon }(t,x,c)dc\right) ^{p} &+&\left( \int_{-\infty
}^{0}(1-f_{\varepsilon }(t,x,c))dc\right) ^{p}dx\leqslant C,  \notag \\
\mathbb{E}\sup_{t\in (0,T)}\int_{\mathbb{R}^{2}}|c|^{2s}\,d\nu _{\varepsilon
}^{(t,x)}(c)dx &\leqslant &C,\quad \mathbb{E}\int_{\mathbb{R}%
_{T}^{2}}|c|^{2(s-1)}\,dm_{\varepsilon }(t,x,c)\leqslant C  \label{mm}
\end{eqnarray}%
with $p=2s$ for $s=1$ and $s=q$. Moreover we have the relations 
\begin{equation}
f_{\varepsilon }(t,x,c)=\nu _{\varepsilon }^{(t,x)}(c,+\infty ),\quad \nu
_{\varepsilon }^{(t,x)}(c)=-\partial _{c}f_{\varepsilon }(t,x,c)\geqslant
0\quad \text{on }\Omega \times \mathbb{R}_{T}.  \label{nu}
\end{equation}

2) the functions $f_{\varepsilon }$ and $\widehat{f_{\varepsilon }}%
=1-f_{\varepsilon }$ verify the following transport equations%
\begin{eqnarray}
d\left( f_{\varepsilon }(t),\varphi \right) &=&\left[ \left( f_{\varepsilon
},c\partial _{x}\varphi \right) -\left( \partial _{x}p_{\varepsilon
}f_{\varepsilon },\partial _{c}\varphi \right) +\varepsilon \left(
f_{\varepsilon },\partial _{xx}^{2}\varphi \right) \right] dt+J_{\varepsilon
}(\varphi ),  \notag \\
d(\widehat{f_{\varepsilon }}(t),\varphi ) &=&[(\widehat{f_{\varepsilon }}%
,c\partial _{x}\varphi )-(\partial _{x}p_{\varepsilon }\widehat{%
f_{\varepsilon }},\partial _{c}\varphi )+\varepsilon (\widehat{%
f_{\varepsilon }},\partial _{xx}^{2}\varphi )]\,dt-J_{\varepsilon }(\varphi
),  \label{15:01:20:55b}
\end{eqnarray}%
for any $\varphi =\varphi (x,c)\in C_{c}^{\infty }(\mathbb{R}^{2}),$ where%
\begin{eqnarray}
J_{\varepsilon }(\varphi ) &=&\sum_{k}\int_{\mathbb{R}}\varphi
(x,u_{\varepsilon }(t,x))\sigma _{k}(x,u_{\varepsilon }(t,x))\,dxd{\mathcal{W%
}_{t}^{k}}  \notag \\
&+&\frac{1}{2}\int_{\mathbb{R}}\partial _{c}\varphi (x,u_{\varepsilon
}(t,x))\sigma ^{2}(x,u_{\varepsilon }(t,x))\,dxdt-m_{\varepsilon }(\partial
_{c}\varphi )  \notag \\
&=&\sum_{k}\int_{\mathbb{R}^{2}}\varphi (x,c)\sigma _{k}(x,c)\,d\nu
_{\varepsilon }^{(t,x)}(c)dxd{\mathcal{W}_{t}^{k}}  \notag \\
&+&\frac{1}{2}\int_{\mathbb{R}^{2}}\partial _{c}\varphi (x,c)\sigma
^{2}(x,c)\,d\nu _{\varepsilon }^{(t,x)}(c)dxdt-m_{\varepsilon }(\partial
_{c}\varphi )  \label{J=}
\end{eqnarray}

3) moreover the pair $(u_{\varepsilon },p_{\epsilon })$ fulfills%
\begin{eqnarray}
d\left( |u_{\varepsilon }-c|^{+},\psi \right) &=&(f_{\varepsilon }(\frac{%
u_{\varepsilon }^{2}}{2}-\frac{c^{2}}{2}),\partial _{x}\psi )dt-\left(
\partial _{x}p_{\varepsilon }|u_{\varepsilon }-c|^{+},\partial _{c}\psi
\right) \,dt  \notag \\
&-&\varepsilon \left( f_{\varepsilon }\partial _{x}u_{\varepsilon },\partial
_{x}\psi \right) dt  \notag \\
&+&\sum_{k}\int_{\mathbb{R}^{2}}f_{\varepsilon }\psi (x,c)\sigma
_{k}(x,u_{\varepsilon }(t,x))\,dcdxd{\mathcal{W}_{r}^{k}}  \notag \\
&+&\frac{1}{2}\int_{\mathbb{R}^{2}}f_{\varepsilon }\partial _{c}\psi
(x,c)\sigma ^{2}(x,u_{\varepsilon }(t,x))\,dcdxdt  \notag \\
&-&m_{\varepsilon }(\psi ),\qquad \forall \psi =\psi (x,c)\in C_{c}^{\infty
}(\mathbb{R}^{2}).  \label{15:01:20:550}
\end{eqnarray}%
Here we use a short notation%
\begin{equation*}
m_{\varepsilon }(\psi )=\int_{\mathbb{R}^{2}}\psi (x,c)\,dm_{\varepsilon
}(t,x,c)
\end{equation*}%
for this integral with respect of the measure $m_{\varepsilon }(t,x,c)$
defined by \eqref{numu}.
\end{lemma}

\medskip

\begin{remark}
Let us point that the equations \eqref{15:01:20:55b}, \eqref{15:01:20:550}
are understood in the distributional sense on the time variable. For instance, the first equation of \eqref{15:01:20:55b} 
should be read as 
\begin{eqnarray}
\int_{0}^{T}\left( f_{\varepsilon },\partial _{t}\varphi \right) dt
&+&\left( f_{\varepsilon, 0 },\varphi (0)\right) +\int_{0}^{T}\left[ \left(
f_{\varepsilon },c\partial _{x}\varphi \right) -\left( \partial
_{x}p_{\varepsilon }\;f_{\varepsilon },\partial _{c}\varphi \right)
+\varepsilon \left( f_{\varepsilon },\partial _{xx}^{2}\varphi \right) %
\right] dt  \notag \\
&=&-\sum_{k}\int_{\mathbb{R}_{T}^{2}}\varphi (t,x,c)\sigma _{k}(x,c)\,d\nu
_{\varepsilon }^{(t,x)}(c)dxd{\mathcal{W}_{t}^{k}}  \label{15:01:20:55a} \\
&-&\frac{1}{2}\int_{\mathbb{R}_{T}^{2}}\partial _{c}\varphi (t,x,c)\sigma
^{2}(x,c)\,d\nu _{\varepsilon }^{(t,x)}(c)dxdt+m_{\varepsilon }(\partial
_{c}\varphi ) \quad \text{with }f_{\varepsilon,
0}=sign^{+}(u_{\varepsilon, 0}-c)   \notag
\end{eqnarray}%
a.s. in $\Omega $ and $\forall \varphi \in
C^{1}([0,T];C_{c}^{\infty }(\mathbb{R}^{2})),\ $such that $\varphi
|_{t=T}=0. $
\end{remark}

\bigskip

\textbf{\ Proof of Lemma \ref{Lemma 4.4}.} 1) Let $p=2s$ for $s=1$ or $s=q$. 
{\ Performing standard } calculations for $f_{\varepsilon
}(t,x,c)=1_{u_{\varepsilon }(t,x)\geqslant c}$ we can verify that 
\begin{eqnarray*}
\int_{0}^{+\infty }f_{\varepsilon }(t,x,c)p{c}^{p-1}dc+\int_{-\infty
}^{0}(1-f_{\varepsilon }(t,x,c))p|c|^{p-1}dc &=&|u_{\varepsilon }(t,x)|^{p},
\\
\left( \int_{0}^{+\infty }f_{\varepsilon }(t,x,c)dc\right) ^{p}+\left(
\int_{-\infty }^{0}(1-f_{\varepsilon }(t,x,c))dc\right) ^{p}
&=&|u_{\varepsilon }(t,x)|^{p},
\end{eqnarray*}%
then the estimates \eqref{Z} imply the first two ones in \eqref{mm}. \
Moreover the definitions of the measures $m_{\varepsilon },$ $\nu
_{\varepsilon }$\ and \eqref{Z} give the last two estimates in \eqref{mm}.
We easily check that $f_{\varepsilon }$ and $\nu _{\varepsilon }^{(t,x)}$
satisfy the relations \eqref{nu}.

2) Let us consider $\theta =\theta (c)\in C_{c}^{\infty }(\mathbb{R})$. The
It\^{o} formula gives 
\begin{align}
d\theta (u_{\varepsilon }(t))& =\theta ^{\prime }(u_{\varepsilon })\left(
-u_{\varepsilon }\partial _{x}u_{\varepsilon }-\partial _{x}p_{\varepsilon
}\right) dt+\varepsilon \theta ^{\prime }(u_{\varepsilon })\partial
_{xx}^{2}u_{\varepsilon }dt  \notag \\
& +\sum_{k}\theta ^{\prime }(u_{\varepsilon })\sigma _{k}(x,u_{\varepsilon
}(t,x))\,d{\mathcal{W}_{t}^{k}}+\frac{1}{2}\sum_{k}\theta ^{\prime \prime
}(u_{\varepsilon })\sigma _{k}^{2}(x,u_{\varepsilon }(t,x))\,dt.
\label{10:09:27:01}
\end{align}%
Each term of this equation can be written as a duality relation. More
precisely, the left hand side corresponds to 
\begin{equation}
d\theta (u_{\varepsilon })=d\int_{-\infty }^{u_{\varepsilon }}\theta
^{\prime }(c)\,dc=d\int_{\mathbb{R}}\theta ^{\prime }(c)f_{\varepsilon
}(c)\,dc.  \label{Hy120}
\end{equation}%
For the terms in the right hand side of \eqref{10:09:27:01}, we have 
\begin{eqnarray}
\theta ^{\prime }(u_{\varepsilon })u_{\varepsilon }\partial
_{x}u_{\varepsilon } &=&\partial _{x}\left( \int_{-\infty }^{u_{\varepsilon
}}\theta ^{\prime }(c)c\,dc\right) =\partial _{x}\left( \int_{\mathbb{R}%
}\theta ^{\prime }(c)cf_{\varepsilon }(c)\,dc\right) ,  \notag \\
\theta ^{\prime }(u_{\varepsilon })\partial _{x}p_{\varepsilon } &=&\int_{%
\mathbb{R}}\partial _{x}p_{\varepsilon }\theta ^{\prime \prime
}(c)f_{\varepsilon }\,dc,  \notag \\
\theta ^{\prime }(u_{\varepsilon })\partial _{xx}^{2}u_{\varepsilon }
&=&\partial _{xx}^{2}\left( \theta (u_{\varepsilon })\right) -|\partial
_{x}u_{\varepsilon }|^{2}\theta ^{\prime \prime }(u_{\varepsilon })  \notag
\\
&=&\partial _{xx}^{2}\int_{\mathbb{R}}\theta ^{\prime }(c)f_{\varepsilon
}(c)\,dc-|\partial _{x}u_{\varepsilon }|^{2}\theta ^{\prime \prime
}(u_{\varepsilon }).  \label{10:19:26:01}
\end{eqnarray}

Introducing the relations \eqref{Hy120}-\eqref{10:19:26:01} in the equation %
\eqref{10:09:27:01}, multiplying by $\phi (x)\in C_{c}^{\infty }(\mathbb{R})$%
, and integrating over $\mathbb{R}$, we infer that 
\begin{eqnarray}
d\left( f_{\varepsilon }(t),\varphi \right) &=&\left( f_{\varepsilon
},c\partial _{x}\varphi \right) \,dt-\left( \partial _{x}p_{\varepsilon
}f_{\varepsilon },\partial _{c}\varphi \right) \,dt+\varepsilon \left(
f_{\varepsilon },\partial _{xx}^{2}\varphi \right) \,dt  \notag \\
&-&m_{\varepsilon }(\partial _{c}\varphi )+\sum_{k}\int_{\mathbb{R}}\varphi
(x,u_{\varepsilon }(t,x))\sigma _{k}(x,u_{\varepsilon }(t,x))\,dxd{\mathcal{W%
}_{t}^{k}}  \notag \\
&+&\frac{1}{2}\int_{\mathbb{R}}\partial _{c}\varphi (x,u_{\varepsilon
}(t,x))\sigma ^{2}(x,u_{\varepsilon }(t,x))\,dxdt  \label{10:15:27:01}
\end{eqnarray}%
for $\varphi (x,c)=\phi (x)\theta ^{\prime }(c).$ Since the test functions $%
\varphi (x,c)=\phi (x)\theta ^{\prime }(c)$ form a dense subset of $%
C_{c}^{\infty }(\mathbb{R}^{2}),$ we obtain the validality of the first
transport equation of \eqref{15:01:20:55b}.

Analogously, let us point%
\begin{eqnarray}
\theta ^{\prime \prime }(u_{\varepsilon })\sigma ^{2}(x,u_{\varepsilon
}(t,x)) &=&\int_{\mathbb{R}}\partial _{c}\left( \theta ^{\prime \prime
}(c)\sigma ^{2}(x,c)\right) f_{\varepsilon }(c)dc  \notag \\
&=&\int_{\mathbb{R}}\theta ^{\prime \prime }(c)\sigma ^{2}(x,c)d\nu
_{\varepsilon }^{(t,x)}(c),  \notag \\
\sum_{k}\theta ^{\prime }(u_{\varepsilon })\sigma _{k}(x,u_{\varepsilon
}(t,x))\,d{\mathcal{W}_{t}^{k}} &=& \sum_{k}\int_{\mathbb{R}%
}\theta ^{\prime }(c)\sigma _{k}(x,c)d\nu _{\varepsilon }^{(t,x)}(c)d{%
\mathcal{W}_{t}^{k}},  \label{10:16:26:01}
\end{eqnarray}%
that gives the equality \eqref{J=}.

Following the same reasoning and accounting that 
\begin{equation*}
\theta (u_{\varepsilon })=-\int_{\mathbb{R}}\theta ^{\prime
}(c)(1-f_{\varepsilon }(c))\,dc,
\end{equation*}%
we can also derive the second transport equation of \eqref{15:01:20:55b}.

3) Now let us use the relation 
\begin{equation}
\partial _{c}|u_{\varepsilon }-c|^{+}=-f_{\varepsilon }\qquad \text{in\quad }%
\mathcal{D}^{\prime }(\mathbb{R}_{T}^{2})  \label{Ue}
\end{equation}%
and the identity%
\begin{equation*}
\sum_{k}\varphi (x,u_{\varepsilon }(t,x))\sigma _{k}(x,u_{\varepsilon
}(t,x))\,d{\mathcal{W}_{t}^{k}}=\sum_{k}\int_{\mathbb{R}}f_{\varepsilon
}(c)\partial _{c}\varphi (x,c)\sigma _{k}(x,u_{\varepsilon }(t,x))\,dcd{%
\mathcal{W}_{t}^{k}}
\end{equation*}%
in the first transport equation of \eqref{15:01:20:55b}, we derive%
\begin{eqnarray}
d\left( |u_{\varepsilon }-c|^{+},\partial _{c}\varphi \right) &=&\left(
|u_{\varepsilon }-c|^{+},\partial _{x}\varphi +c\partial _{c}\partial
_{x}\varphi \right) dt  \notag \\
&-&\left( \partial _{x}p_{\varepsilon }|u_{\varepsilon }-c|^{+},\partial
_{cc}^{2}\varphi \right) \,dt  \notag \\
&+&\varepsilon \left( |u_{\varepsilon }-c|^{+},\partial _{c}\partial
_{xx}^{2}\varphi \right) dt  \notag \\
&+&\sum_{k}\int_{\mathbb{R}^{2}}f_{\varepsilon }(c)\partial _{c}\varphi
(x,c)\sigma _{k}(x,u_{\varepsilon }(t,x))\,dcdxd{\mathcal{W}_{r}^{k}}  \notag
\\
&+&\frac{1}{2}\int_{\mathbb{R}}f_{\varepsilon }(c)\partial _{cc}^{2}\varphi
(x,c)\sigma ^{2}(x,u_{\varepsilon }(t,x))\,dcdxdt  \notag \\
&-&m_{\varepsilon }(\partial _{c}\varphi )  \label{Cc}
\end{eqnarray}%
by the integration by parts on the parameter $c\in \mathbb{R}$. Moreover we
have 
\begin{eqnarray*}
\left( |u_{\varepsilon }-c|^{+},\partial _{x}\varphi \right) &=&\left(
\int_{c}^{+\infty }|u_{\varepsilon }-s|^{+}ds,\partial _{c}\partial
_{x}\varphi \right) \\
&=&\left( f_{\varepsilon }\left( \frac{u_{\varepsilon }^{2}}{2}%
-u_{\varepsilon }c+\frac{c^{2}}{2}\right) ,\partial _{x}\partial _{c}\varphi
\right) .
\end{eqnarray*}%
Substituting this relation inside of \eqref{Cc} and denoting by $\psi
=\partial _{c}\varphi ,$ we get \eqref{15:01:20:550}.$\hfill \;\blacksquare $

\bigskip

In order to pass to the limit on $\varepsilon \rightarrow 0$ in the
nonlinear terms of the transport equations \eqref{15:01:20:55b}, we need to
prove an additional uniform estimate with respect to the time variable.

\begin{lemma}
\label{10:42:26:01} Let us consider the sequence $\left\{ U_{\varepsilon
}=|u_{\varepsilon }-c|^{+}\right\} _{\varepsilon >0}.$ For each natural $%
n\in \mathbb{N}$ there exists some constant $C(n),$ which is independent of $%
\varepsilon $ (but may depend on $n$), such that%
\begin{equation}
\mathbb{E}\sup_{0\leqslant \theta \leqslant \delta }\int_{0}^{T-\delta
}||U_{\varepsilon }(t+\theta )-U_{\varepsilon
}(t)||_{H^{-3}(B_{n})}^{2}dt\leqslant C(n)\delta ,  \label{timeestimate1}
\end{equation}%
where 
\begin{equation*}
||h||_{H^{-3}(B_{n})}=\sup_{\psi \in H_{0}^{3}(B_{n})}\left\{ |(h,\psi
)|:\quad ||\psi ||_{H_{0}^{3}(B_{n})}\leqslant 1\right\}
\end{equation*}%
is the norm of $H^{-3}(B_{n})$ with \quad $B_{n}=\left\{ (c,x)\in
(-n,n)\times (-n,n)\right\} .$
\end{lemma}

\textbf{Proof}. \ For any positive $\theta $, according to the equality %
\eqref{15:01:20:550}, we can write 
\begin{equation}
\left( |u_{\varepsilon }-c|^{+}(t+\theta )-|u_{\varepsilon }-c|^{+}(t),\psi
\right) =\left( F_{_{\varepsilon ,\theta }}(t),\psi \right) +\left(
G_{_{\varepsilon ,\theta }}(t),\psi \right) ,  \label{L}
\end{equation}%
where%
\begin{eqnarray*}
\left( F_{_{\varepsilon ,\theta }}(t),\psi \right) &=&\int_{t}^{t+\theta
}\{(f_{\varepsilon }\left( \frac{u_{\varepsilon }^{2}}{2}-\frac{c^{2}}{2}%
\right) ,\partial _{x}\psi )-\left( \partial _{x}p_{\varepsilon
}|u_{\varepsilon }-c|^{+},\partial _{c}\psi \right) \\
&-&\varepsilon \left( f_{\varepsilon }\partial _{x}u_{\varepsilon },\partial
_{x}\psi \right) +\frac{1}{2}\int_{\mathbb{R}}\psi (x,u_{\varepsilon
}(s,x))\sigma ^{2}(x,u_{\varepsilon }(s,x))\,dx \\
&-&m_{\varepsilon }(\psi )\}\,ds,
\end{eqnarray*}%
and 
\begin{equation}
\left( G_{_{\varepsilon ,\theta }}(t),\psi \right)
=\sum_{k}\int_{t}^{t+\theta }\int_{\mathbb{R}^{2}}\psi (x,c)f_{\varepsilon
}(s,x,c)\sigma _{k}(x,u_{\varepsilon }(s,x))\,dcdxd{\mathcal{W}_{s}^{k}.}
\label{Gg}
\end{equation}

By the continuous embedding $H_{0}^{2}(B_{n})\hookrightarrow C(B_{n}),$
applied to $\psi \in H_{0}^{3}(B_{n})$, we infer that 
\begin{eqnarray*}
|\left( F_{_{\varepsilon ,\theta }}(t),\psi \right) |
&\leqslant &C(n)\int_{t}^{t+\theta }\{\left( 1+||u_{\varepsilon }||_{L^{2}(%
\mathbb{R})}^{2}+||\partial _{x}\mathrm{p}_{\varepsilon }||_{L^{\infty }(%
\mathbb{R})}||u_{\varepsilon }||_{L^{2}(\mathbb{R})}\right. \\
&+&\left. \varepsilon \left\Vert \partial _{x}u_{\varepsilon }\right\Vert
_{L^{2}(\mathbb{R})}\right) \left\Vert \nabla _{x,c}\psi \right\Vert
_{L^{\infty }(B_{n})} \\
&+&\left( \int_{\mathbb{R}}\sigma ^{2}(x,u_{\varepsilon
}(s,x))\,dx+\varepsilon \left\Vert \partial _{x}u_{\varepsilon }\right\Vert
_{L^{2}(\mathbb{R})}^{2}\right) \Vert \psi \Vert _{L^{\infty }(B_{n})}\}\,ds.
\end{eqnarray*}%
The assumption \eqref{ass1} and the Holder inequality imply 
\begin{equation*}
||F_{\varepsilon ,\theta }(t)||_{H^{-3}(B_{n})}\leqslant
C(n)\int_{t}^{t+\theta }{h_{\varepsilon }(s)~ds}
\end{equation*}%
with ${h_{\varepsilon }(s)=}1+\Vert u_{\varepsilon }\Vert _{L^{2}(\mathbb{R}%
)}^{2}+||\partial _{x}\mathrm{p}_{\varepsilon }||_{L^{\infty }(\mathbb{R}%
)}^{2}+\varepsilon \left\Vert \partial _{x}u_{\varepsilon }\right\Vert
_{L^{2}(\mathbb{R})}^{2}.$ Hence for a small fixed $\delta >0,$ taking the
supremum over $\theta \in (0,\delta )$\ and the integration on $t\in
(0,T-\delta )$, we derive 
\begin{eqnarray}
\mathbb{E}\int_{0}^{T-\delta }\left( \sup_{0\leqslant \theta \leqslant
\delta }||F_{_{\varepsilon ,\theta }}(t)||_{H^{-3}(B_{n})}\right) {dt}
&\leqslant &C(n)\int_{0}^{T-\delta }\int_{t}^{t+\delta }{h_{\varepsilon
}(s)~dsdt}  \notag \\
&\leqslant &C(n)\delta ,  \label{F}
\end{eqnarray}%
by the estimates \eqref{Z}$_{1}$ for $s=1$ and \eqref{Z1}$_{1}$ for $%
r=+\infty $.

Now we consider the operator $G_{\varepsilon ,\theta }$ given by \eqref{Gg}.
Let us consider an orthonormal basis $\{\psi _{j}\}$ of $L^{2}(B_{n}),$ such
that $\{\psi _{j}\}\subset H_{0}^{2}(B_{n}).$ We notice that $\psi _{j}$ can
be identified with its extension by zero to the domain $\mathbb{R}^{2}.$
With the help of the Burkholder-Davis-Gundy inequality and the Moore-Osgood
theorem on exchanging limits, we have the equalities%
\begin{align*}
& \mathbb{E}\int_{0}^{T-\delta }\sup_{0\leqslant \theta \leqslant \delta
}||G_{\varepsilon ,\theta }(t)||_{HS(L^{2}(B_{n}),\mathbb{R})}^{2}dt \\
& =\mathbb{E}\int_{0}^{T-\delta }\sup_{0\leqslant \theta \leqslant \delta
}\sum_{j}\left( \sum_{k}\int_{t}^{t+\theta }\int_{B_{n}}\psi
_{j}(x,c)f_{\varepsilon }(s,x,c)\sigma _{k}(x,u_{\varepsilon }(s,x))\,dcdxd{%
\mathcal{W}_{s}^{k}}\right) ^{2}{dt} \\
& =\sum_{j}\mathbb{E}\int_{0}^{T-\delta }\sup_{0\leqslant \theta \leqslant
\delta }\left( \sum_{k}\int_{t}^{t+\theta }\int_{B_{n}}\psi
_{j}(x,c)f_{\varepsilon }(s,x,c)\sigma _{k}(x,u_{\varepsilon }(s,x))\,dcdxd{%
\mathcal{W}_{s}^{k}}\right) ^{2}{dt} \\
& =\sum_{j}\mathbb{E}\int_{0}^{T-\delta }\int_{t}^{t+\delta }\sum_{k}\left(
\int_{B_{n}}\psi _{j}(x,c)f_{\varepsilon }(s,x,c)\sigma
_{k}(x,u_{\varepsilon }(s,x))\,dcdx\right) ^{2}~ds \\
& =\sum_{k}\mathbb{E}\int_{0}^{T-\delta }\int_{t}^{t+\delta }\sum_{j}(\psi
_{j},f_{\varepsilon }\sigma _{k}(\cdot ,u_{\varepsilon }(s,\cdot
)))^{2}~dsdt=J.
\end{align*}%
Hence applying the Parseval formula, the H\"{o}lder inequality, the
assumption \eqref{20:01:21:59} and the estimate \eqref{Z} for $s=1$, we
obtain 
\begin{align*}
J& =\mathbb{E}\int_{0}^{T-\delta }\int_{t}^{t+\delta }\sum_{k}\Vert
f_{\varepsilon }\sigma _{k}(\cdot ,u_{\varepsilon }(s,\cdot ))\Vert
_{L^{2}(B_{n})}^{2}~dsdt \\
& \leqslant C(n)\mathbb{E}\int_{0}^{T-\delta }\int_{t}^{t+\delta
}\sum_{k}\Vert \sigma _{k}(\cdot ,u_{\varepsilon }(s,\cdot ))\Vert _{L^{2}(%
\mathbb{R})}^{2}~dsdt \\
& \leqslant C(n)\mathbb{E}\int_{0}^{T-\delta }\int_{t}^{t+\delta }\Vert
u_{\varepsilon }(s,\cdot )\Vert _{L^{2}(\mathbb{R})}^{2}~dsdt\leqslant
\delta C(n)\mathbb{E}\Vert u_{\varepsilon }\Vert _{L^{\infty }(0,T;L^{2}(%
\mathbb{R}))}^{2}\leqslant C(n)\delta .
\end{align*}%
Since {%
\begin{equation*}
\Vert G_{\varepsilon ,\theta }(t)\Vert _{H^{-3}(B_{n})}\leqslant
||G_{\varepsilon ,\theta }(t)||_{HS(L^{2}(B_{n}),\mathbb{R})},
\end{equation*}
} we obtain the following estimate 
\begin{equation}
\mathbb{E}\int_{0}^{T-\delta }\sup_{0\leqslant \theta \leqslant \delta
}||G_{\varepsilon ,\theta }(t)||_{H^{-3}(B_{n})}^{2}dt\leqslant C(n)\delta .
\label{E}
\end{equation}%
Combining the derived estimates \eqref{F}, \eqref{E} with the relation %
\eqref{L}, we obtain \eqref{timeestimate1}.$\hfill \;\blacksquare $

\bigskip

\section{Tightness}

\label{Sec5}

\setcounter{equation}{0}

The passage to the limit, as $\varepsilon \rightarrow 0,$ will be performed
with the help of the Jakubowski-Skorokhod representation Theorem \ref%
{Jakubowski}, which relays on a tightness property over a sub-Polish
topological space $\mathcal{E}$. First, we introduce some appropriate
functional spaces and identify their compact sets. Then we define $\mathcal{E%
}$ as a product space and show that the family of $\mathcal{E}$-valued
random variables 
\begin{equation*}
\Phi _{\varepsilon }=(U_{\varepsilon },u_{\varepsilon },f_{\varepsilon
},\delta _{u_{\epsilon }=c},m_{\varepsilon },p_{\varepsilon },\mathcal{W})
\end{equation*}%
is a tight family (see \cite{BEF18}, Sec. 2.6-2.8, for a discussion on the
application of Jakubowski-Skorokhod theorem).

Let us consider the Sobolev space $H_{0}^{1}(B_{n})$ and its topological
dual space $H^{-1}(B_{n})=(H_{0}^{1}(B_{n}))^{\ast }$. Since each $v\in
H_{0}^{1}(B_{n})$ can be extended by zero outside of $B_{n}$ onto the domain 
$\mathbb{R}^{2}$, \ the sequence $\{H_{0}^{1}(B_{n})\}_{n\in \mathbb{N}}$ of
sets is an increasing one on the parameter $n\in \mathbb{N}$. As a direct
consequence we get that $\{H^{-1}(B_{n})\}_{n\in \mathbb{N}}$ forms a
decreasing sequence on $n,$ that permits to introduce the space 
\begin{equation*}
L^{2}(0,T;H_{loc}^{-1})=\bigcap_{n\in \mathbb{N}}L^{2}(0,T;H^{-1}(B_{n})).
\end{equation*}%
We consider the space $L^{2}(0,T;H_{loc}^{-1})$ with the topology induced by
the metric 
\begin{equation*}
d(f,g):=\sum_{n\in \mathbb{N}}\frac{1}{2^{n}}\max \left( {1,\Vert f-g\Vert
_{L^{2}(0,T;H^{-1}(B_{n}))}}\right) .
\end{equation*}

\medskip

In the following considerations the duality relation of a normed space $B$
and its dual space $B^{\ast }$ is denoted by $\langle \cdot ,\cdot \rangle .$

\medskip

For any $R>0$ we introduce the set 
\begin{equation*}
K_{R}^{1}=\left\{ v\in L^{2}(0,T;H_{loc}^{-1}):\;|||v|||\leqslant R\right\} ,
\end{equation*}%
endowed with the norm%
\begin{align*}
|||v|||& =\sum_{n\in \mathbb{N}}\frac{1}{n^{2}C(n)}(\Vert v\Vert
_{L^{2}(0,T;L^{2}(B_{n}))}+ \\
& +\sup_{\delta }\mathbb{\delta }^{-1}\{\sup_{0\leqslant \theta \leqslant
\delta }\int_{0}^{T-\delta }||v(t+\theta )-v(t)||_{H^{-3}(B_{n})}^{2}dt\}),
\end{align*}%
where $C(n)$ are the constants defined in Lemma \ref{10:42:26:01}. We have
the following continuous embeddings 
\begin{equation}
L^{2}(B_{n})\hookrightarrow H^{-1}(B_{n})\hookrightarrow H^{-3}(B_{n}),\quad
\forall n\in \mathbb{N}.  \label{embedding}
\end{equation}%
The first embedding is compact. Hence applying Lemma \ref{simon} and using
the diagonalization process on $n\in \mathbb{N}$, we verify that $K_{R}^{1}$
is a compact subset of $L^{2}(0,T;H_{loc}^{-1})$.

Let us consider a strictly increasing smooth function $\theta:\mathbb{R}\to
[-1,1].$

For each $n\in \mathbb{N}$ we can choose a countable dense subset $%
\{g_{i}^{n}\}_{i\in \mathbb{N}}$ of a separable space $%
L^{2}(0,T;H_{0}^{1}(B_{n}))$.\ Now we consider the duality relation $\langle
\cdot ,\cdot \rangle $ for the spaces $B=L^{2}(0,T;H_{0}^{1}(B_{n}))$ and $%
B^{\ast }=L^{2}(0,T;H^{-1}(B_{n}))$. Then the countable family of $\mathbb{R}%
-$valued continuous functions 
\begin{equation*}
\mathbb{M}_{1}:=\{\Psi _{i}^{n}:L^{2}(0,T;H_{loc}^{-1})\rightarrow {[-1,1]}%
,\quad \Psi _{i}^{n}\left( f\right) ={\theta(\langle f,g_{i}^{n}\rangle)}
,\;\forall i,n\in \mathbb{N}\}
\end{equation*}%
separate points of $L^{2}(0,T;H_{loc}^{-1}).$

\medskip

Next, for each $R>0$, the set%
\begin{equation*}
K_{R}^{2}=\biggl\{v\in L^{\infty }(0,T;L^{2}(\mathbb{R})\cap L^{2q}(\mathbb{R%
})):\;\Vert v\Vert _{L^{\infty }(0,T;L^{2}(\mathbb{R}))}+\Vert v\Vert
_{L^{\infty }(0,T;L^{2q}(\mathbb{R}))}\leqslant R\biggr\}
\end{equation*}%
is a weakly$-\ast $ compact set of $L^{\infty }(0,T;L^{2}(\mathbb{R})\cap
L^{2q}(\mathbb{R}))$. We recall that the spaces $L^{2}(\mathbb{R%
})\cap L^{2q}(\mathbb{R})$ and $L^{2}(\mathbb{R})\oplus L^{\lambda }(\mathbb{%
R})$, with $\frac{1}{\lambda }+\frac{1}{2q}=1,$ are duals of each other (see
in \cite{SS}, Chapter 1, p. 41, Exercise 26). 
Thus the space $L^{\infty }(0,T;L^{2}(\mathbb{R})\cap L^{2q}(\mathbb{R}))$
can be identified with the dual of $L^{1}(0,T;L^{2}(\mathbb{R})\oplus
L^{\lambda }(\mathbb{R}))$, with $\frac{1}{\lambda }+\frac{1}{2q}=1.$ Since
the space $L^{1}(0,T;L^{2}(\mathbb{R})\oplus L^{\lambda }(\mathbb{R}))$ is
separable, we can take a contable dense subset $\{\alpha _{i},i\in \mathbb{N}%
\}$ such that 
\begin{equation*}
\mathbb{M}_{2}:=\{\Psi _{i}:L^{\infty }(0,T;L^{2}(\mathbb{R})\cap L^{2q}(%
\mathbb{R}))\rightarrow {[-1,1]},\quad \Psi _{i}\left( f\right) ={\theta
(\langle f,\alpha _{i}\rangle )},\;i\in \mathbb{N}\}
\end{equation*}%
is a countable family of continuous funtions, which separate points of $%
L^{\infty }(0,T;L^{2}(\mathbb{R})\cap L^{2q}(\mathbb{R})).$

Now, let us consider $L^{\infty }(\mathbb{R}_{T})$ endowed with the weak$%
-\ast $ topology. The set 
\begin{equation*}
K_{R}^{3}=\biggl\{g\in L^{\infty }(\mathbb{R}_{T}):\Vert g\Vert _{L^{\infty
}(\mathbb{R}_{T})}\leqslant R\biggr\},
\end{equation*}%
is a compact set of $L^{\infty }(\mathbb{R}_{T})$, since any sequence of $%
K_{R}^{3}$ has a subsequence, which converges weakly$-\ast $ in $L^{\infty }(%
\mathbb{R}_{T}).$ Taking into account that $L^{\infty }(\mathbb{R}_{T})$ is
the dual of the separable space $L^{1}(\mathbb{R}_{T})$, using analogous
arguments, we can define a countable family $\mathbb{M}_{3}$ of $\mathbb{R}$%
-valued continuous functions defined on $L^{\infty }(\mathbb{R}_{T})$, which
separate points of $L^{\infty }(\mathbb{R}_{T})$.

\medskip

Let us consider the space $L^{\infty }(0,T;W^{1,1}(\mathbb{R})\cap
W^{1,\infty }(\mathbb{R}))$ endowed with the weak-$\ast $ topology. The set 
\begin{align*}
K_{R}^{4}=& \bigl\{\pi \in L^{\infty }(0,T;W^{1,1}(\mathbb{R})\cap
W^{1,\infty }(\mathbb{R})): \\
& \qquad \qquad \qquad \qquad \;\Vert \pi \Vert _{L^{\infty }(0,T;W^{1,1}(%
\mathbb{R}))}+\Vert \pi \Vert _{L^{\infty }(0,T;W^{1,\infty }(\mathbb{R}%
))}\leqslant R\bigr\},
\end{align*}%
is a compact set $L^{\infty }(0,T;W^{1,1}(\mathbb{R})\cap W^{1,\infty }(%
\mathbb{R}))$. Obviously, any sequence of  $K_{R}^{4}$  has a
subsequence, which converges weakly-$\ast $ in $L^{\infty }(0,T;W^{1,r}(%
\mathbb{R}))$ for any$\ r\in \lbrack 1,+\infty ].$ We know that for $r\in
]1,+\infty \lbrack $, the space $L^{\infty }(0,T;W^{1,r}(\mathbb{R}))$
corresponds to the topological dual of the separable space $%
L^{1}(0,T;W^{-1,r^{\prime }}(\mathbb{R}))$ with $\frac{1}{r}+\frac{1}{%
r^{\prime }}=1$.

Reasoning as above, taking a contable dense subset $\{h_{i},i\in \mathbb{N}%
\}, $ we define a countable family of $\mathbb{R}-$valued continuous
functions 
\begin{equation*}
\mathbb{M}_{4,r}:=\{\Psi _{i}:L^{\infty }(0,T;W^{1,r}(\mathbb{R}%
))\rightarrow {[-1,1]},\quad \Psi _{i}\left( f\right) ={\theta (\langle
f,h_{i}\rangle )},\;i\in \mathbb{N}\},
\end{equation*}%
which separates points of $L^{\infty }(0,T;W^{1,r}(\mathbb{R})).$ In case $%
r=+\infty $, the space $L^{\infty }(0,T;W^{1,\infty }(\mathbb{R}))$ can be
identified with a subset of $L^{\infty }(0,T;(L^{\infty }(\mathbb{R}))^{2})$
through the mapping $u\rightarrow (u,\partial _{x}u)$. In addition $%
L^{\infty }(0,T;(L^{\infty }(\mathbb{R}))^{2})$ can also be identified with
the dual of $L^{1}(0,T;(L^{1}(\mathbb{R}))^{2})$. Again, using a countable
dense subset $\{\tilde{h}_{i},i\in \mathbb{N}\}$ of $L^{1}(0,T;(L^{1}(%
\mathbb{R}))^{2})$, the set of functions 
\begin{equation*}
\mathbb{M}_{4,\infty }:=\{\Psi _{i}:L^{\infty }(0,T;W^{1,r}(\mathbb{R}%
))\rightarrow {[-1,1]},\quad f\rightarrow {\theta (\langle f,\tilde{h}%
_{i}\rangle )},\;i\in \mathbb{N}\},
\end{equation*}%
separates points of $L^{\infty }(0,T;W^{1,\infty }(\mathbb{R})).$ We set $%
\mathbb{M}_{4}=\mathbb{M}_{4,r},$ $1\leqslant r\leqslant +\infty $.

\medskip

The space of Young measures $\mathcal{M}_{Y}\left( \mathbb{R}_{T}\right) $
from $\mathbb{R}_{T}$ to $\mathcal{P}(\mathbb{R})$ (the space of probability
measures on $(\mathbb{R},\mathcal{B}(\mathbb{R})$) corresponds to $L^{\infty
}(\mathbb{R}_{T};\mathcal{P}(\mathbb{R}))$ which is non separable. The space 
$\mathcal{P}(\mathbb{R})$ is a subset of $\mathcal{M}_{R}\left( \mathbb{R}%
\right) $, the space of Radom measures on $(\mathbb{R},\mathcal{B}(\mathbb{R}%
)$). We recall that $\mathcal{M}_{R}\left( \mathbb{R}\right) $ is the
topological dual of $C_{c}\left( \mathbb{R}\right) $, the space of real
valued continuous functions defined on $\mathbb{R}$ vanishing at infinity.
Thus 
\begin{equation*}
\mathcal{M}_{Y}\left( \mathbb{R}_{T}\right) \subset (L^{1}(\mathbb{R}%
_{T};C_{c}\left( \mathbb{R}\right) )^{\ast }.
\end{equation*}%
The space $L^{1}(\mathbb{R}_{T};C_{c}\left( \mathbb{R}\right) )$ is
separable. Namely, the space $C_{c}\left( \mathbb{R}\right) $ is a separable
Banach space, then we consider a countable dense subset $\mathbb{C}%
_{0}:=\{\psi _{j}:\;j\in \mathbb{N}\}.$ The Borel $\sigma $-algebra $%
\mathcal{B}(\mathbb{R}_{T})$ is generated by the following countable set 
\begin{equation*}
\mathbb{A}_{0}=\{A=[a,b]\times \lbrack c,d]:\;\;\,a,b,c,d\in \mathbb{Q}%
,\;\;0\leqslant a<b\leqslant T\}.
\end{equation*}%
We consider $\mathcal{M}_{Y}\left( \mathbb{R}_{T}\right) $ endowed with the
so-called narrow topology, defined as the weakest topology for which the set
of following functions are continuous 
\begin{equation*}
\nu \rightarrow \int_{A}\int_{\mathbb{R}}\varphi (c)d\nu
^{(t,x)}(c)dxdt,\quad \varphi \in C_{b}(\mathbb{R}),\quad A\in \mathcal{B}(%
\mathbb{R}_{T}),
\end{equation*}%
where $C_{b}(\mathbb{R})$ is the space of continuous and bounded functions
on $\mathbb{R}$. Therefore 
\begin{align*}
\biggl\{& \Psi :\mathcal{M}_{Y}(\mathbb{R}_{T}) \rightarrow \lbrack
-1,1],\;\; \\
& \Psi (\nu )=\frac{1}{|A|\sup |\psi| }\int_{A}\left[ \int_{S}\psi (s)d\nu
^{z}(s)\right] d\mu (z),\;A\in \mathbb{A}_{0},\psi \in \mathbb{C}_{0}\biggr\}%
,
\end{align*}
where $|A|$ denotes the measure of $A$, is a set of countinuous functions,
which separates the points of $\mathcal{M}_{Y}(\mathbb{R}_{T}).$ This means
that {\ the space $\mathcal{M}_{Y}(\mathbb{R}_{T})$ is a sub-Polish then it
is suitable to apply the version of the Skorohod theorem established by
Jakubowski \cite{J98}. }

On the other hand, given $R>0$ and $q>1$, the function $h:\mathbb{R}%
_{T}^{2}\rightarrow \mathbb{R}$, $h(t,x,z)=|z|^{q}$, is inf-compact, i. e.
for each $\lambda >0$ and $(t,x)\in \mathbb{R}_{T}$, $\{z:\;h(t,x,z)%
\leqslant \lambda \}$ is a compact subset of $\mathbb{R}$, therefore the set
defined by 
\begin{equation*}
K_{R}^{5}=\left\{ \nu \in \mathcal{M}_{Y}\left( \mathbb{R}_{T}\right)
:\;|\nu |^{q}\leqslant R\right\}
\end{equation*}%
with 
\begin{equation*}
|\nu |^{q}=\int_{\mathbb{R}_{T}}\int_{\mathbb{R}}|c|^{q}d\nu ^{(t,x)}(c)dxdt,
\end{equation*}%
is a compact subset of $\mathcal{M}_{Y}\left( \mathbb{R}_{T}\right) $.

\medskip

Here, we consider $\mathcal{M}^{+}(\mathbb{R}_{T}^{2})$, the space of
positive finite measures on $\mathbb{R}_{T}^{2}.$ According to \cite%
{Prokhorov}, $\mathcal{M}^{+}(\mathbb{R}_{T}^{2})$ induced with the weak
topology is a Polish space, then $\mathcal{M}^{+}(\mathbb{R}_{T}^{2})$ is a
sub-Polish space. There exists a countable family $\mathbb{M}_{6}$ of $%
[-1,1] $-valued continuous functions defined on $\mathcal{M}^{+}(\mathbb{R}%
_{T}^{2}) $, which separate points of $\mathcal{M}^{+}(\mathbb{R}_{T}^{2}). $
On the other hand, for each $r>0$ and $R>0$, the set 
\begin{equation*}
K_{R}^{6}=\biggl\{m\in \mathcal{M}^{+}(\mathbb{R}_{T}^{2}):\int_{\mathbb{R}%
_{T}^{2}}|c|^{r}dm(t,x,c)\leqslant R\biggr\}
\end{equation*}%
is a compact set of $\mathcal{M}^{+}(\mathbb{R}_{T}^{2})$.

\medskip

Finally, the law $\mu $ of the Wiener process $\mathcal{W}$,
being a probability measure on $C([0,T],\mathbb{H})$, is a Radom measure and
the singular set $\{\mu \}$ is tight on $C([0,T],\mathbb{H})$. Then for
every $R>0$, there exists a compact set $K_{R}^{7}$ of $C([0,T],\mathbb{H})$
such that $\mu (K_{R}^{7})\geqslant 1-\frac{1}{R}$.

\medskip

Now, we consider the product space with the corresponding product topology 
\begin{align*}
\mathcal{E}& =L^{2}(0,T;H_{loc}^{-1})\times L^{\infty }(0,T;L^{2}(\mathbb{R}%
)\cap L^{2q}(\mathbb{R})) \\
& \times L^{\infty }(\mathbb{R}_{T})\times \mathcal{M}_{Y}\left( \mathbb{R}%
_{T}\right) \\
& \times \mathcal{M}^{+}(\mathbb{R}^{2}_T)\times L^{\infty }(0,T;W^{1,1}(%
\mathbb{R})\cap W^{1,\infty }(\mathbb{R}))\times C([0,T],\mathbb{H}).
\end{align*}
According to the the above considerations for the sets $K_{R}^{i},$ $%
i=1,...,7,$ the Cartesian product 
\begin{equation*}
K_{R}=K_{R}^{1}\times ...\times K_{R}^{7}
\end{equation*}%
is a compact set of the topological space $\mathcal{E}$. Moreover, the set 
\begin{equation*}
\mathbb{M}:=\{\widetilde{\Psi}_i: \mathcal{E}\to [-1,1], \quad (z_1,
z_2,z_3,z_4,z_5,z_6, z_7)\to \Psi_i(z_i), \quad \Psi_i\in \mathbb{M}_i,
\;i=1,\dots,7 \},
\end{equation*}
is a countable family of continuous functions, which separates points of $%
\mathcal{E}.$

Let us set by $\sigma(\mathbb{M})$, the $\sigma$-algebra generated by the
family of functions $\mathbb{M},$ and consider the measurable space $(%
\mathcal{E}, \sigma(\mathbb{M})). $ Let us introduce the family of random
variables 
\begin{equation*}
\Phi _{\varepsilon }=(U_{\varepsilon },u_{\varepsilon },f_{\varepsilon
},\delta _{u_{\epsilon }=c},m_{\varepsilon },p_{\varepsilon },\mathcal{W})
\end{equation*}%
with values in $\mathcal{E}$ and the family of the laws $\{\mu _{\varepsilon
}\}_{\varepsilon >0}$, defined by 
\begin{equation*}
\mu _{\varepsilon }(\Gamma )=P(\Phi _{\varepsilon }\in \Gamma ),\qquad
\Gamma \in \sigma (\mathbb{M}),
\end{equation*}%
which are probability measures on $(\mathcal{E},\sigma (\mathbb{M})).$

\bigskip

By the standard way applying Chebyshev%
\'{}%
s inequality and the apriori estimates\ \eqref{Z}-\eqref{Z1}, \eqref{mm}, %
\eqref{timeestimate1}, for arbitrary fixed $\delta >0$, we can choose $%
R_{\delta }$ big enough, such that 
\begin{equation*}
P\left( \Phi _{\varepsilon }\notin K_{R_{\delta }}\right) \leqslant \frac{C}{%
R_{\delta }^{2}}<\delta ,
\end{equation*}%
we refer for details to \textit{4.3.3. Tightness Property} of \cite{Ben} and
Lemma 4.8 of \cite{RS}. Thus $\mu _{\varepsilon }(K_{R_{\delta }})>1-\delta $
for all $\varepsilon >0$, which gives the tightness of the family of
measures $\{\mu _{\varepsilon }\}$ on the space $\mathcal{E}$.

\bigskip

\section{Existence of weak stochastic kinetic solution}

\label{Sec6} \setcounter{equation}{0}

This section is devoted to show the main result (Theorem \ref{the_1}) on the
existence of a weak stochastic kinetic solution for the equation, according
to Definition \ref{weak_kinetic}.

\subsection{Limit transition in viscous kinetic equations}

\label{sec32}

Taking into account the tightness of $\{\Phi _{\varepsilon }\}$, the
Skorohod-Jakubowski's theorem  (see also Theorem C.1 in \cite%
{brzezniak2018stochastic}), gives the existence of a probability space $%
\left( \widetilde{\Omega },\widetilde{\mathcal{F}},\widetilde{P}\right) $
and random variables 
\begin{equation*}
\widetilde{\Phi }_{\varepsilon }=(\widetilde{U}_{\varepsilon },\widetilde{u}%
_{\varepsilon },\widetilde{f}_{\varepsilon },\widetilde{p}_{\varepsilon },%
\widetilde{\nu }_{\varepsilon },\widetilde{m}_{\varepsilon },\widetilde{%
\mathcal{W}}_{t}),\qquad \widetilde{\Phi }=(\widetilde{U},\widetilde{u},%
\widetilde{f},\widetilde{p},\widetilde{\nu },\widetilde{m},\widetilde{%
\mathcal{W}}_{t})
\end{equation*}%
defined on this probability space with values in $\mathcal{E}$, such that

\begin{enumerate}
\item $\widetilde{\Phi }_{\varepsilon }$ and $\Phi _{\varepsilon }$ have the
same law $\mu _{\varepsilon }$.

\item There exists a subsequence of $\{\widetilde{\Phi }_{\varepsilon }\}$,
still denoted by $\{\widetilde{\Phi }_{\varepsilon }\}$, which converges $%
\widetilde{P}-$a.e. to $\widetilde{\Phi }$ in the topological space $%
\mathcal{E}$. Therefore, for $\widetilde{P}-$a.e. $\omega \in \widetilde{%
\Omega }$ 
\begin{align}
\widetilde{U}_{\varepsilon }& \rightarrow \widetilde{U}\quad \text{strongly
in $L^{2}(0,T;H_{loc}^{-1}(\mathbb{R}^{2})),$}  \notag \\
\widetilde{u}_{\varepsilon }& \rightharpoonup \widetilde{u}\quad \text{%
weakly -- $\ast $ in in }L^{\infty }(0,T;L^{2}(\mathbb{R})\cap L^{2q}(%
\mathbb{R})),  \notag \\
\widetilde{f}_{\varepsilon }& \rightharpoonup \widetilde{f}\quad \text{%
weakly -- $\ast $ in weakly in }L^{\infty }(\mathbb{R}_{T}^{2}),  \notag \\
\widetilde{\mathrm{p}}_{\varepsilon }& \rightharpoonup \widetilde{\mathrm{p}}%
\quad \text{weakly in }L^{\infty }(0,T;W^{1,1}(\mathbb{R})\cap W^{1,\infty }(%
\mathbb{R})),  \notag \\
\widetilde{\nu }_{\varepsilon }& \rightharpoonup \widetilde{\nu }\quad \text{%
weakly in }\mathcal{M}_{Y}\left( \mathbb{R}_{T}\right) ,  \notag \\
\widetilde{m}_{\varepsilon }& \rightharpoonup \widetilde{m}\quad \text{%
weakly in }\mathcal{M}^{+}(\mathbb{R}_{T}^{2}).  \label{19:17:08}
\end{align}
\end{enumerate}


Now, we define an appropriate filtration $\{\widetilde{\mathcal{F}}%
_{t}\}_{t\in \lbrack 0,T]}$ on the probability space $\left( \widetilde{%
\Omega },\widetilde{\mathcal{F}},\widetilde{P}\right) $ for which $%
\widetilde{\mathcal{W}}_{t}$ is a Wiener process and the random
distributions $\widetilde{U}$, $\widetilde{u}$, $\widetilde{f}$, $\widetilde{%
\mathrm{p}}$, $\widetilde{\nu }$, $\widetilde{m}$ are predictable (for
details see \cite{BEF18}, Section 2.2) 
\begin{eqnarray*}
\mathcal{G}_{\widetilde{U}}^{t} &=&\sigma \left\{ (\widetilde{U},\psi
),\;\psi \in C_{c}^{\infty }([0,t]\times \mathbb{R}^{2})\right\} , \\
\mathcal{G}_{\widetilde{u}}^{t} &=&\sigma \left\{ (\widetilde{u},\psi
),\;\psi \in C_{c}^{\infty }([0,t]\times \mathbb{R}^{2})\right\} , \\
\mathcal{G}_{\widetilde{f}}^{t} &=&\sigma \left\{ (\widetilde{f},\psi
),\;\psi \in C_{c}^{\infty }([0,t]\times \mathbb{R}^{2})\right\} , \\
\mathcal{G}_{\widetilde{p}}^{t} &=&\sigma \left\{ (\widetilde{p},\psi
),\;\psi \in C_{c}^{\infty }([0,t]\times \mathbb{R})\right\} , \\
\mathcal{G}_{\widetilde{\nu }}^{t} &=&\sigma \left\{ (\widetilde{\nu} ,\psi
),\;\psi \in C_{c}^{\infty }([0,t]\times \mathbb{R},C_{b}(\mathbb{R}%
))\right\} , \\
\mathcal{G}_{\widetilde{m}}^{t} &=&\sigma \left\{ (\widetilde{m},\psi
),\;\psi \in C_{c}^{\infty }([0,t]\times \mathbb{R}^{2})\right\} , \\
\mathcal{G}_{\widetilde{\mathcal{W}}}^{t} &=&\sigma \left\{ \widetilde{%
\mathcal{W}}_{\lambda },\;0\leqslant \lambda \leqslant t\right\} , \\
\widetilde{\mathcal{F}}_{t} &=&\bigcap_{s>t}\sigma \left\{ \mathcal{G}_{%
\widetilde{U}}^{s}\cup \mathcal{G}_{\widetilde{f}}^{s}\cup \mathcal{G}_{%
\widetilde{p}}^{s}\cup \mathcal{G}_{\widetilde{\nu }}^{s}\cup \mathcal{G}_{%
\widetilde{m}}^{s}\cup \mathcal{G}_{\widetilde{\mathcal{W}}}^{s}\cup
\{N:\;N\subset B\in \widetilde{\mathcal{F}},\quad \widetilde{P}%
(B)=0\}\right\} ,
\end{eqnarray*}%
where $\sigma \{Y_{i}:\,i\in \Gamma \}$ denotes the $\sigma $-algebra
generated by the family of random variables $Y_{i},$ $\,i\in \Gamma .$
\bigskip

Recalling that $U_{\varepsilon }=|u_{\varepsilon }-c|^{+}$ and relation %
\eqref{Ue}, let us set 
\begin{eqnarray*}
X_{\varepsilon } &=&\int_{\mathbb{R}_{T}^{2}}\left\vert U_{\varepsilon
}-|u_{\varepsilon }-c|^{+}\right\vert dtdxdc=0,\text{\quad a.s. in }\Omega 
\text{,} \\
Y_{\varepsilon } &=&\left\vert \int_{\mathbb{R}_{T}^{2}}U_{\varepsilon
}\partial _{c}\varphi dtdxdc-\int_{\mathbb{R}_{T}^{2}}f_{\varepsilon
}\varphi dtdxdc\right\vert =0,\text{\quad }\forall \varphi {\in }\mathcal{S}%
\subset C_{c}^{\infty }(\mathbb{R}_{T}^{2}),
\end{eqnarray*}%
where $\mathcal{S}$ is a countable dense subset of $C_{c}^{\infty }(\mathbb{R%
}_{T}^{2})$ and 
\begin{eqnarray*}
\widetilde{X}_{\varepsilon } &=&\int_{\mathbb{R}_{T}^{2}}\left\vert 
\widetilde{U}_{\varepsilon }-|\widetilde{u}_{\varepsilon }-c|^{+}\right\vert
dtdxdc,\text{\quad\ a.s. in }\widetilde{\Omega }, \\
\widetilde{Y}_{\varepsilon } &=&\left\vert \int_{\mathbb{R}_{T}^{2}}%
\widetilde{U}_{\varepsilon }\partial _{c}\varphi dtdxdc-\int_{\mathbb{R}%
_{T}^{2}}\widetilde{f}_{\varepsilon }\varphi dtdxdc\right\vert .
\end{eqnarray*}
Since ${X}_{\varepsilon }$ is a non-negative continuous
deterministic functional of ${\Phi }_{\varepsilon }$, then there exists a
non-negative, continuous, bounded deterministic functional 
\begin{equation*}
g:\mathcal{E}\to ]0, +\infty[ \quad\text{such that}\quad \frac{%
X_{\varepsilon }}{X_{\varepsilon }+1}=g({\Phi }_{\varepsilon }).
\end{equation*}
An analogous expression holds when we replace $X_{\varepsilon }$ and ${\Phi }%
_{\varepsilon }$ by $\widetilde{X}_{\varepsilon }$ and $\widetilde{\Phi }%
_{\varepsilon }$, respectively, namely 
\begin{equation*}
\frac{\widetilde{X}_{\varepsilon }}{\widetilde{X}_{\varepsilon }+1} =g(%
\widetilde{\Phi }_{\varepsilon }).
\end{equation*}
In addition, taking into account that $\Phi _{\varepsilon }$ and $\widetilde{%
\Phi }_{\varepsilon }$ are equal in law, we deduce 
\begin{eqnarray*}
0 &=&\mathbb{E}\left[ \frac{X_{\varepsilon }}{X_{\varepsilon }+1}\right] =%
\mathbb{E}\left[ g(\Phi _{\varepsilon })\right] =\widetilde{\mathbb{E}}\left[
g(\widetilde{\Phi }_{\varepsilon })\right] =\widetilde{\mathbb{E}}\left[ 
\frac{\widetilde{X}_{\varepsilon }}{\widetilde{X}_{\varepsilon }+1}\right],
\end{eqnarray*}
{where $\widetilde{\mathbb{E}}$ denotes the mathematical expectation with
respect to the measure $\widetilde{P}$.} The same argument shows that 
\begin{eqnarray*}
0 &=&\mathbb{E}\left[ \frac{Y_{\varepsilon }}{Y_{\varepsilon }+1}\right] =%
\widetilde{\mathbb{E}}\left[ \frac{\widetilde{Y}_{\varepsilon }}{\widetilde{Y%
}_{\varepsilon }+1}\right].
\end{eqnarray*}
Therefore $\widetilde{X}_{\varepsilon }=0$ and $\widetilde{Y}_{\varepsilon
}=0$ a.s.$-\widetilde{P}.$ Consequentlly, we obtain a.e. in $\widetilde{%
\Omega }\times \mathbb{R}_{T}$ 
\begin{equation}
\widetilde{U}_{\varepsilon }=|\widetilde{u}_{\varepsilon }-c|^{+}\quad \text{%
and}\quad \partial _{c}\widetilde{U}_{\varepsilon }=-\widetilde{f}%
_{\varepsilon }\text{\quad in\quad }\mathcal{D}^{\prime }(\mathbb{R}_{T}).
\label{fe0}
\end{equation}%
In addition, with the help of \eqref{19:17:08}$_{1,3}$, we pass to the limit
as $\varepsilon \rightarrow 0$ in these relations and deduce%
\begin{equation}
\widetilde{U}=|\widetilde{u}-c|^{+}\quad \text{and}\quad \partial _{c}%
\widetilde{U}=-\widetilde{f}\text{\quad\ ~\ in\quad } \mathcal{D%
}^{\prime }(\mathbb{R}_{T}).   \label{fe1}
\end{equation}

Applying a similar reasoning as above for $\widetilde{U}_{\varepsilon },%
\widetilde{U},$ as a direct consequence of \eqref{fe0}-\eqref{fe1}, we have
that 
\begin{eqnarray}
0 &\leqslant &\widetilde{f}_{\varepsilon }\leqslant 1,\quad 0\leqslant 
\widetilde{f}\leqslant 1\quad \text{a.e. in }\widetilde{\Omega }\times 
\mathbb{R}_{T}^{2}\quad \text{and}  \notag \\
\quad \widetilde{\nu }_{\varepsilon } &=&-\frac{\partial \widetilde{f}%
_{\varepsilon }}{\partial c}\geqslant 0\text{,}\quad \widetilde{\nu }=-\frac{%
\partial \widetilde{f}}{\partial c}\geqslant 0\quad \text{in \ }\mathcal{D}%
^{\prime }(\mathbb{R})\quad \text{a.e. in }\widetilde{\Omega }\times \mathbb{%
R}_{T},  \label{ff0}
\end{eqnarray}%
as well as the following system in the distributional sense a.s. in $%
\widetilde{\Omega }$ 
\begin{equation}
\left\{ 
\begin{array}{l}
d\widetilde{u}_{\varepsilon }=\left( -\widetilde{u}_{\varepsilon }\partial
_{x}\widetilde{u}_{\varepsilon }-\partial _{x}\widetilde{\mathrm{p}}%
_{\varepsilon }+\varepsilon \partial _{xx}^{2}\widetilde{u}_{\varepsilon
}\right) dt+\sigma (\widetilde{u}_{\varepsilon })d{\widetilde{\mathcal{W}}}%
_{t}\qquad \text{in\quad }\mathbb{R}_{T}, \\ 
\\ 
\left( 1-\partial _{xx}^{2}\right) \widetilde{\mathrm{p}}_{\varepsilon }=%
\frac{3}{2}\widetilde{u}_{\varepsilon }^{2}\qquad \text{in\quad }\mathbb{R}%
_{T}, \\ 
\\ 
\widetilde{u}_{\varepsilon }(0)=u_{\varepsilon ,0}\qquad \qquad \text{ in }%
\mathbb{R}.%
\end{array}%
\right.  \label{Sys11}
\end{equation}%
Arguing as in the proof of Lemmas \ref{lem1}, \ref{Lemma 4.4} and \ref%
{10:42:26:01}, the uniform estimates \eqref{Z}-\eqref{Z1}, \eqref{mm} and %
\eqref{timeestimate1} can be deduced for the pair $(\widetilde{u}%
_{\varepsilon },\widetilde{\mathrm{p}}_{\varepsilon })$, yielding%
\begin{align}
\varepsilon \nabla \widetilde{u}_{\varepsilon }& \rightarrow 0\quad \text{%
strongly in $L^{2}(\widetilde{\Omega }\times \mathbb{R}_{T})$},  \notag \\
\widetilde{u}_{\varepsilon }& \rightharpoonup \widetilde{u}\quad \text{%
weakly -- $\ast $ in $L^{2}(\widetilde{\Omega },L^{\infty }(0,T;L^{2}(%
\mathbb{R})))$},  \notag \\
\widetilde{u}_{\varepsilon }& \rightharpoonup \widetilde{u}\quad \text{%
weakly -- $\ast $ in $L^{2q}(\widetilde{\Omega },L^{\infty }(0,T;L^{2}(%
\mathbb{R})\cap L^{2q}(\mathbb{R})))$},  \notag \\
\widetilde{f}_{\varepsilon }& \rightharpoonup \widetilde{f}\quad \text{%
weakly -- $\ast $ in }L^{\infty }(\widetilde{\Omega }\times \mathbb{R}%
{}_{T}^{2}\text{$),$}  \notag \\
\widetilde{\mathrm{p}}_{\varepsilon }& \rightharpoonup \widetilde{\mathrm{p}}%
\quad \text{weakly -- $\ast $ in }L^{2}(\widetilde{\Omega },L^{\infty
}(0,T;W^{1,r}(\mathbb{R})),\qquad \forall r\in \lbrack 1,+\infty ],  \notag
\\
\partial _{xx}^{2}\widetilde{\mathrm{p}}_{\varepsilon }& \rightharpoonup
\partial _{xx}^{2}\widetilde{\mathrm{p}}_{\varepsilon }\quad \text{weakly -- 
$\ast $ in }L^{q}(\widetilde{\Omega },L^{\infty }(0,T;L^{q}(\mathbb{R})))%
\text{$,$}  \notag \\
\widetilde{\nu }_{\varepsilon }& \rightharpoonup \widetilde{\nu }\quad \text{%
weakly in }\mathcal{M}_{Y}\left( \mathbb{R}_{T}\right) ,  \notag \\
\widetilde{m}_{\varepsilon }& \rightharpoonup \widetilde{m}\quad \text{%
weakly in }\mathcal{M}^{+}(\mathbb{R}_{T}^{2}).  \label{limit0}
\end{align}%
Let us point that the last limit in (\ref{limit0}) is a direct consequence
of Theorem 1.41, p. 66 of \cite{EG}.

The uniform estimate \eqref{Z}$_{2}$ valid  for $\widetilde{u}%
_{\varepsilon }$ gives the following estimate 
\begin{equation*}
\widetilde{\mathbb{E}}\Vert \widetilde{U}_{\varepsilon }\Vert
_{L^{2}(0,T;H_{loc}^{-1}(\mathbb{R}^{2}))}^{2q}\leqslant C,
\end{equation*}%
where $C$ is a positive constant independent of $\varepsilon $. On the other
hand, the convergence \eqref{19:17:08}$_{1}$ and the Fatou's lemma yield 
\begin{equation*}
\widetilde{\mathbb{E}}\Vert \widetilde{U}\Vert _{L^{2}(0,T;H_{loc}^{-1}(%
\mathbb{R}^{2}))}^{2q}\leqslant C.
\end{equation*}%
Since $2q>2$, we infer that the set of random variables ${\Vert \widetilde{U%
}_{\varepsilon }-\widetilde{U}\Vert _{L^{2}(0,T;H_{loc}^{-1}(\mathbb{R}%
^{2}))}^{2}},${\ }${\varepsilon >0,}${\ is uniformly integrable. In
addition, due to \eqref{19:17:08}$_{1}$ we have 
\begin{equation*}
\Vert \widetilde{U}_{\varepsilon }-\widetilde{U}\Vert
_{L^{2}(0,T;H_{loc}^{-1}(\mathbb{R}^{2}))}^{2}\rightarrow 0\quad \text{a.s.
as }\quad \varepsilon \rightarrow 0.
\end{equation*}%
With the help of the Vitali theorem, we deduce 
\begin{equation*}
\widetilde{\mathbb{E}}\Vert \widetilde{U}_{\varepsilon }-\widetilde{U}\Vert
_{L^{2}(0,T;H_{loc}^{-1}(\mathbb{R}^{2}))}^{2}\rightarrow 0,\quad \text{ as }%
\quad \varepsilon \rightarrow 0,
\end{equation*}%
then 
\begin{equation}
\widetilde{U}_{\varepsilon }\rightarrow \widetilde{U}\quad \text{strongly in 
}L^{2}(\widetilde{\Omega };\text{$L^{2}(0,T;H_{loc}^{-1}(\mathbb{R}^{2})))$}.
\label{strong}
\end{equation}%
}

As above arguing that the random variables {$\Phi _{\varepsilon }$} and {$%
\widetilde{\Phi }_{\varepsilon }$} are equal in law, we can demonstrate that
the functions $\widetilde{f}_{\varepsilon }$ and $1-\widetilde{f}%
_{\varepsilon }$ satisfy the transport equations \eqref{15:01:20:55b} a.s.
in $\widetilde{\Omega }.$ These equations are understood in the
distributional sense on the time variable. For instance, $\widetilde{f}%
_{\varepsilon }$ fulfills 
\begin{eqnarray}
\int_{0}^{T}(\widetilde{f}_{\varepsilon },\partial _{t}\varphi )dt &+&\left(
f_{\varepsilon, 0},\varphi (0)\right) +\int_{0}^{T}\left[ (\widetilde{f}%
_{\varepsilon },c\partial _{x}\varphi )-(\partial _{x}\widetilde{p}%
_{\varepsilon }\;\widetilde{f}_{\varepsilon },\partial _{c}\varphi
)+\varepsilon (\widetilde{f}_{\varepsilon },\partial _{xx}^{2}\varphi )%
\right] dt  \notag \\
&=&-\sum_{k}\int_{\mathbb{R}_{T}^{2}}\varphi (t,x,c)\sigma _{k}(x,c)\,d%
\widetilde{\nu }_{\varepsilon }^{(t,x)}(c)dxd\widetilde{{\mathcal{W}}}{%
_{t}^{k}}  \notag \\
&-&\frac{1}{2}\int_{\mathbb{R}_{T}^{2}}\partial _{c}\varphi (t,x,c)\sigma
^{2}(x,c)\,d\widetilde{\nu }_{\varepsilon }^{(t,x)}(c)dxdt+\widetilde{m}%
_{\varepsilon }(\partial _{c}\varphi )  \label{1a}
\end{eqnarray}%
for any $\varphi \in C^{1}([0,T];C_{c}^{\infty }(\mathbb{R}%
^{2})),\ $such that$\ \varphi |_{t=T}=0.$

Using the convergences \eqref{strong}, \eqref{limit0}$_{4,5}$ and the
relations \eqref{fe0}, \eqref{fe1}, we deduce 
\begin{eqnarray*}
&&\int_{0}^{T}(\partial _{x}\widetilde{p}_{\varepsilon }\,\widetilde{f}%
_{\varepsilon },\partial _{c}\varphi )d\tau =-\int_{0}^{T}\int_{\mathbb{R}%
}(\partial _{x}\widetilde{p}_{\varepsilon }\,\widetilde{U}_{\varepsilon
},\partial _{cc}^{2}\varphi )dcdt\rightarrow \\
&&\qquad \qquad -\int_{0}^{T}\int_{\mathbb{R}}(\partial _{x}\widetilde{p}\,%
\widetilde{U},\partial _{cc}^{2}\varphi )dcdt=\int_{0}^{T}(\partial _{x}%
\widetilde{p}\,\widetilde{f},\partial _{c}\varphi )dt\quad \text{weakly in $%
L^{2}(\widetilde{\Omega }).$}
\end{eqnarray*}%
Moreover applying the weak convergences \eqref{limit0}, we easily pass on
the limit transition as $\varepsilon \rightarrow 0$ in other terms of %
\eqref{1a} weakly in $L^{2}(\widetilde{\Omega }).$\ Therefore the functions $%
\widetilde{u},\widetilde{f},\widetilde{\mathrm{p}},$ $\widetilde{\nu }$ and$%
\ \widetilde{m},$ obtained by the limit transition \eqref{limit0}, satisfy
for a.s. $\widetilde{\Omega }$ the integral equality%
\begin{eqnarray}
\int_{0}^{T}(\widetilde{f},\partial _{t}\varphi )dt &+&{\left( f_{0},\varphi
(0)\right) }+\int_{0}^{T}\left[ (\widetilde{f},c\partial _{x}\varphi
)-(\partial _{x}\widetilde{p}\;\widetilde{f},\partial _{c}\varphi )\right] dt
\notag \\
&=&-\sum_{k}\int_{\mathbb{R}_{T}^{2}}\varphi (t,x,c)\sigma _{k}(x,c)\,d%
\widetilde{\nu }^{(t,x)}(c)dxd\widetilde{{\mathcal{W}}}{_{t}^{k}}  \notag \\
&-&\frac{1}{2}\int_{\mathbb{R}_{T}^{2}}\partial _{c}\varphi (t,x,c)\sigma
^{2}(x,c)\,d\widetilde{\nu }^{(t,x)}(c)dxdt+\widetilde{m}(\partial
_{c}\varphi ).  \label{15:01:20:55e}
\end{eqnarray}%
Analogously, we deduce 
\begin{eqnarray}
\int_{0}^{T}(\widehat{\widetilde{f}},\partial _{t}\varphi )dt &+&{\left(
\left( 1-f_{0}\right) ,\varphi (0)\right) }+\int_{0}^{T}\left[ (\widehat{%
\widetilde{f}},c\partial _{x}\varphi )-(\partial _{x}\widetilde{p}\;\widehat{%
\widetilde{f}},\partial _{c}\varphi )\right] dt  \notag \\
&=&\sum_{k}\int_{\mathbb{R}_{T}^{2}}\varphi (t,x,c)\sigma _{k}(x,c)\,d%
\widetilde{\nu }^{(t,x)}(c)dxd\widetilde{{\mathcal{W}}}{_{t}^{k}}  \notag \\
&+&\frac{1}{2}\int_{\mathbb{R}_{T}^{2}}\partial _{c}\varphi (t,x,c)\sigma
^{2}(x,c)\,d\widetilde{\nu }^{(t,x)}(c)dxdt-\widetilde{m}(\partial
_{c}\varphi ),  \label{15:01:20:555e}
\end{eqnarray}%
hence we can formulate the following result.

\begin{lemma}
\label{limit}The limit functions $\widetilde{f}$ and $\widehat{\widetilde{f}}%
(t)=1-\widetilde{f}(t)$ verify the following transport equations%
\begin{eqnarray}
d(\widetilde{f}(t),\varphi ) &=&\left[ (\widetilde{f},c\partial _{x}\varphi
)-(\partial _{x}\widetilde{p}\widetilde{f},\partial _{c}\varphi )\right]
dt+J(\varphi ),  \notag \\
d(\widehat{\widetilde{f}}(t),\varphi ) &=&[(\widehat{\widetilde{f}}%
,c\partial _{x}\varphi )-(\partial _{x}\widetilde{p}\widehat{\widetilde{f}}%
,\partial _{c}\varphi )]\,dt-J(\varphi ),  \label{1b}
\end{eqnarray}%
for any $\varphi \in C_{c}^{\infty }(\mathbb{R}^{2}),$ where%
\begin{eqnarray*}
J(\varphi ) &=&\sum_{k}\int_{\mathbb{R}^{2}}\varphi (t,x,c)\sigma
_{k}(x,c)\,d\widetilde{\nu }^{(t,x)}(c)dxd\widetilde{{\mathcal{W}}}{_{t}^{k}}
\\
&+&\frac{1}{2}\int_{\mathbb{R}^{2}}\partial _{c}\varphi (t,x,c)\sigma
^{2}(x,c)\,d\widetilde{\nu }^{(t,x)}(c)dxdt-\widetilde{m}(\partial
_{c}\varphi )
\end{eqnarray*}%
and%
\begin{equation*}
\widetilde{m}(\partial _{c}\varphi )=\int_{\mathbb{R}^{2}}\partial
_{c}\varphi (x,c)\,d\widetilde{m}(t,x,c).
\end{equation*}
\end{lemma}

\bigskip

\subsection{$L^p$-integrability and the existence of the traces of $%
\widetilde{f}$}

\bigskip

In this section we establish technical results concerning the properties of
the function $\widetilde{f}$, and its relation with the function $\widetilde{%
u}$. \medskip

First we demonstrate the following result, which is a consequence of Lemma %
\ref{Lemma 4.4} and the weak convergences (\ref{limit0}).

\begin{lemma}
\label{Lemma 4.440} There exist positive constants $C$, such that%
\begin{eqnarray}
&&{\widetilde{\mathbb{E}}}\int_{\mathbb{R}_{T}}\left( \int_{0}^{+\infty }%
\widetilde{f}(t,x,c)p{c}^{p-1}dc+\int_{-\infty }^{0}(1-\widetilde{f}%
(t,x,c))p|c|^{p-1}dc\right) dxdt\leqslant C,  \notag \\
&&{\widetilde{\mathbb{E}}}\int_{\mathbb{R}_{T}}\left( \int_{0}^{+\infty }%
\widetilde{f}(t,x,c)dc\right) ^{p}+\left( \int_{-\infty }^{0}(1-\widetilde{f}%
(t,x,c))dc\right) ^{p}dxdt\leqslant C,  \notag \\
&&{\widetilde{\mathbb{E}}}\int_{\mathbb{R}_{T}^{2}}|c|^{2s}\,d\widetilde{\nu 
}^{(t,x)}(c)dxdt\leqslant C,\quad {\widetilde{\mathbb{E}}}\int_{\mathbb{R}%
_{T}^{2}}|c|^{2(s-1)}\,d\widetilde{m}(t,x,c)\leqslant C  \label{eq2.88}
\end{eqnarray}%
with $p=2s$ for $s=1$ and $s=q$. Moreover we have the relations 
\begin{equation}
\widetilde{f}(t,x,c)=\widetilde{\nu }^{(t,x)}(c,+\infty ),\quad \widetilde{%
\nu }^{(t,x)}(c)=-\partial _{c}\widetilde{f}(t,x,c)\geqslant 0\quad \text{on 
}\widetilde{\Omega }\times \mathbb{R}_{T}.  \label{eq2.888}
\end{equation}
\end{lemma}

\textbf{Proof}. Next we consider that $p=2s$ for $s=1$ or $s=q$. By the weak
convergence of $\widetilde{f}_{\varepsilon }$ to $\widetilde{f}$ in $%
L^{\infty }(\widetilde{\Omega }\times \mathbb{R}_{T}^{2}),$ \ for any $M>0$
and any $\delta >0$, there exists $\varepsilon _{1}=\varepsilon
_{1}(M,\delta )>0,$ such that for each positive $\varepsilon <\varepsilon
_{1}$, we have 
\begin{equation*}
\left\vert {\widetilde{\mathbb{E}}}\int_{0}^{T}\left( \int_{\mathbb{R}%
^{2}}1_{[-M,M]\times \lbrack 0,M]}\left( \widetilde{f}_{\varepsilon }-%
\widetilde{f}\right) pc^{p-1}dxdc\right) dt\right\vert \leqslant \frac{%
\delta }{3}
\end{equation*}%
where $1_{[-M,M]\times \lbrack 0,M]}$ is the characteristic function of the
set $[-M,M]\times \lbrack 0,M].$ In particular, if we take $\delta =1$ and
use the first estimate in (\ref{mm}), we obtain  
\begin{equation*}
{\widetilde{\mathbb{E}}}\int_{0}^{T}\int_{-M}^{M}\int_{0}^{M}\widetilde{f}%
pc^{p-1}dtdxdc\leqslant \frac{1}{3}+{\widetilde{\mathbb{E}}}\int_{\mathbb{R}%
_{T}}\int_{0}^{+\infty }\widetilde{f}_{\varepsilon }pc^{p-1}dtdxdc\leqslant 
\frac{1}{3}+C\quad \text{for all \ }M.
\end{equation*}
Taking $M\rightarrow \infty ,$ we deduce the estimate  
\begin{equation*}
{\widetilde{\mathbb{E}}}\int_{\mathbb{R}_{T}}\int_{0}^{+\infty }\widetilde{f}%
pc^{p-1}dtdxdc<C
\end{equation*}%
by the monotone convergence theorem. Dealing with the sequence $1-%
\widetilde{f}_{\varepsilon }$\ instead of $\widetilde{f}_{\varepsilon }$,
and applying the same reasoning, we can derive that the second term in the
first estimate of (\ref{eq2.88}) is bounded too.

By another hand if we define $\mathbb{M}_{T}=[0,T]\times \lbrack -M,M]$ for
arbitrary fixed $M>0,$ then the weak convergence of $\widetilde{f}%
_{\varepsilon }$ to $\widetilde{f}$ in $L^{\infty }(\widetilde{\Omega }%
\times \mathbb{R}_{T}^{2})$ gives 
\begin{equation*}
\lim_{\varepsilon \rightarrow 0}{\widetilde{\mathbb{E}}}\int_{\mathbb{M}%
_{T}}\left( \int_{0}^{M}\widetilde{f}_{\varepsilon }(t,x,c)\ dc\right) \psi
\ dxdt={\widetilde{\mathbb{E}}}\int_{\mathbb{M}_{T}}\left( \int_{0}^{M}%
\widetilde{f}(t,x,c)\ dc\right) \psi \ dxdt
\end{equation*}%
for any\ $\psi =\psi (t,x)\in L^{r}(\widetilde{\Omega }\times \mathbb{M}%
_{T}) $ with $\frac{1}{r}+\frac{1}{p}=1,$ since $1_{\mathbb{M}_{T}\times 
\mathbb{[}0,M\mathbb{]}}(t,x,c)\psi \in L^{1}(\widetilde{\Omega }\times 
\mathbb{R}_{T}^{2}).$ Therefore 
\begin{equation*}
\int_{0}^{M}\widetilde{f}_{\varepsilon }(\cdot ,\cdot ,c)\ dc\rightharpoonup
\int_{0}^{M}\widetilde{f}(\cdot ,\cdot ,c)\ dc\quad \text{weakly in \ }L^{p}(%
\widetilde{\Omega }\times \mathbb{M}_{T}).
\end{equation*}%
The second estimate in (\ref{mm}) and the lower semicontinuity of $L^{p}$%
-weakly convergent sequences imply%
\begin{equation*}
{\widetilde{\mathbb{E}}}\int_{\mathbb{M}_{T}}\left( \int_{0}^{M}\widetilde{f}%
(t,x,c)\ dc\right) ^{p}dxdt\leqslant \lim \inf_{\varepsilon \rightarrow 0}{%
\widetilde{\mathbb{E}}}\int_{\mathbb{M}_{T}}\left( \int_{0}^{M}\widetilde{f}%
_{\varepsilon }(t,x,c)dc\right) ^{p}dxdt\leqslant C.
\end{equation*}%
Passing on $M\rightarrow \infty $ and using the monotone convergence
theorem, we deduce the estimate%
\begin{equation*}
{\widetilde{\mathbb{E}}}\int_{\mathbb{R}_{T}}\left( \int_{0}^{+\infty }%
\widetilde{f}(t,x,c)\ dc\right) ^{p}dxdt\leqslant C.
\end{equation*}%
Applying the same reasoning to the sequence $1-\widetilde{f}_{\varepsilon }$%
\ instead of $\widetilde{f}_{\varepsilon }$, we also derive that the second
term in (\ref{eq2.88})$_{2}$ is bounded too.

If we use similar above arguments for the weak convergences (\ref{limit0})
for $\widetilde{\nu }_{\varepsilon }^{(t,x)}$ and $\widetilde{m}%
_{\varepsilon }$ with the help of (\ref{mm})$_{3, 4}$, we derive the last
two estimates for the limit measures $\widetilde{\nu }^{(t,x)}$ and $%
\widetilde{m}.$ The relations (\ref{eq2.888}) follow from (\ref{nu}) and the
weak convergences (\ref{limit0}) for $\widetilde{\nu }_{\varepsilon
}^{(t,x)} $ and $\widetilde{f}_{\varepsilon }$.\hfill $\blacksquare $

\bigskip

The equalities \eqref{15:01:20:55e} and \eqref{15:01:20:555e} yield the
existence of traces for the functions $\widetilde{f}(t)$ and $\widehat{%
\widetilde{f}}(t)$, respectively, for all $t\in \lbrack 0,T]$ (we refer to
cf. Theorem 2.1 in \cite{chen} and \cite{chem3}-\cite{chem1} for more
details).

\begin{lemma}
\label{lemma2 copy(1)} For all $t\in \lbrack 0,T]$, there exists $\widetilde{%
f}^{+}\in L^{\infty }(\widetilde{\Omega }\times {\mathbb{R}}_{T}^{2})$
satisfying 
\begin{equation*}
\widetilde{f}^{+}(t)=\lim_{\delta \rightarrow 0^{+}}\frac{1}{\delta }%
\int_{t}^{t+\delta }  \widetilde{f}(s)\,ds\text{\ \ a.e. on }%
\widetilde{\Omega }\times {\mathbb{R}}^{2},
\end{equation*}%
such that 
\begin{equation}
\widetilde{f}^{+}=\widetilde{f}\text{\ \ \ \ a.e. on }\widetilde{\Omega }%
\times {\mathbb{R}}_{T}^{2}\text{ \ \ \ \ and\ \ \ \ }\widetilde{f}%
^{+}(0)=f_{0}\text{\ \ a.e. on }{\mathbb{R}}^{2}.  \label{2.13}
\end{equation}
\end{lemma}

\noindent \textbf{Proof.} Let us consider $\varphi (t,x,c)=\phi (x,c)\alpha
(t)$ with arbitrary  positive $\phi \in C_{c}^{\infty }(\mathbb{%
R}^{2}),$ and $\alpha \in C^{1}([0,T])$ such that $\alpha (0)=\alpha (T)=0$.
If we define 
\begin{eqnarray*}
I_{\phi }(t) &=&\int_{0}^{t}\left[ (\widetilde{f},c\partial _{x}\phi
)-(\partial _{x}\widetilde{p}\widetilde{f},\partial _{c}\phi )\right] ds+%
\frac{1}{2}\int_{0}^{t}\int_{\mathbb{R}^{2}}\partial _{c}\phi (x,c)\sigma
^{2}(x,c)\,d\widetilde{\nu }^{(s,x)}(c)dxds, \\
J_{\phi }(t) &=&\int_{0}^{t}\sum_{k}\int_{\mathbb{R}^{2}}\phi (x,c)\sigma
_{k}(x,c)\,d\widetilde{\nu }^{(s,x)}(c)dxd\widetilde{{\mathcal{W}}}{_{s}^{k}}%
,
\end{eqnarray*}%
then the integration by parts in the equation \eqref{15:01:20:55e} gives%
\begin{equation*}
\int_{0}^{T}\left[ (\widetilde{f},\phi )(t)-I_{\phi }(t)-J_{\phi }(t)\right]
\partial _{t}\alpha \,dt=\int_{0}^{T}\alpha (t)d\widetilde{m}_{\phi }(t),
\end{equation*}%
where $\widetilde{m}_{\phi }(t)=\int_{\mathbb{R}^{2}}\partial _{c}\phi
(x,c)~d\widetilde{m}(t,x,c)$ is a Radon measure. Thus, we have 
\begin{equation*}
\partial _{t}\left( (\widetilde{f},\phi )(t)-I_{\phi }(t)-J_{\phi
}(t)\right) =-\widetilde{m}_{\phi }(t),
\end{equation*}%
 which allows to conclude that $(\widetilde{f},\phi )(\cdot )-I_{\phi
}(\cdot )-J_{\phi }(\cdot )\in BV([0,T])$, then the left and right limits
exist for any $t\in \lbrack 0,T]$, and the set of discontinuity points is
countable. Since the functions $I_{\phi }(\cdot ),\;J_{\phi }(\cdot )\in
C([0,T])$, we conclude that the function $(\widetilde{f},\phi )(\cdot )$ has
left and right limits at any $t\in \lbrack 0,T]$, and the set of
discontinuity points is countable. Hence there exists the limit%
\begin{equation*}
(\widetilde{f},\phi )(t^{+})=\lim_{\delta \rightarrow 0^{+}}\frac{1}{\delta }%
\int_{t}^{t+\delta }(\widetilde{f},\phi )(s)\,ds,\text{\qquad }\forall \phi
\in C_{c}^{\infty }(\mathbb{R}^{2}),\text{\quad a.s. in }\widetilde{\Omega }.
\end{equation*}%
Let us introduce 
\begin{equation*}
\widetilde{f}_{\delta }:=\frac{1}{\delta }\int_{t}^{t+\delta }\widetilde{f}%
(s)\,ds,\quad {\widetilde{\nu }}_{\delta }^{+(s,x)}:=\frac{1}{%
\delta }\int_{t}^{t+\delta }\widetilde{\nu }^{(s,x)}ds,
\end{equation*}%
satisfying the relation{\ 
\begin{equation}
 -\partial _{c}\widetilde{f}_{\delta }={\widetilde{\nu }}%
_{\delta }^{+(s,x)}.   \label{15:11}
\end{equation}%
} Taking into account that 
\begin{equation*}
0\leqslant \widetilde{f}_{\delta }=\frac{1}{\delta }\int_{t}^{t+\delta }%
\widetilde{f}(s)ds\leqslant 1\text{\quad a.e. on }\widetilde{\Omega }\times {%
\mathbb{R}}_{T}^{2},
\end{equation*}%
then up to a subsequence, $\{\widetilde{f}_{\delta }\}$ converges to some
function $\widetilde{f}^{+}$ in {$L^{\infty }(\widetilde{\Omega }\times {%
\mathbb{R}}_{T}^{2})$} weakly -- $\ast .${\ It follows that 
\begin{equation*}
(\widetilde{f}^{+}(t),\phi )=\lim_{\delta \rightarrow 0^{+}}\frac{1}{\delta }%
\int_{t}^{t+\delta }(\widetilde{f},\phi )(s)\,ds=(\widetilde{f},\phi
)(t^{+}),\text{\qquad }\forall \phi \in C_{c}^{\infty }(\mathbb{R}%
^{2}),\;\forall t\in \lbrack 0,T],
\end{equation*}%
therefore the limit is independent of the subsequence, and consequently }%
a.s. in $\widetilde{\Omega }${\ 
\begin{equation*}
\lim_{\delta \rightarrow 0^{+}}{\widetilde{f}}_{\delta }=\widetilde{f}%
^{+}\quad \text{weakly -- $\ast $ in }L^{\infty }(\mathbb{R}_{T}^{2})\text{.}
\end{equation*}%
The relation \eqref{15:11} implies that there exists a $\mathcal{M}%
_{Y}\left( \mathbb{R}_{T}\right) -$valued random variable ${\widetilde{\nu }}%
^{+(s,x)},$ such that }a.s. in $\widetilde{\Omega }$ 
\begin{equation*}
\lim_{\delta \rightarrow 0^{+}}{\widetilde{\nu }}_{\delta }^{+(s,x)}={%
\widetilde{\nu }}^{+(s,x)}\quad \text{in }\mathcal{M}_{Y}\left(\mathbb{R}%
_{T}\right) \quad \text{and}\quad -\partial _{c}\widetilde{f}^{+}={%
\widetilde{\nu }}^{+(s,x)}.
\end{equation*}%
The limit function $\widetilde{f}^{+}$ is right continuous, namely 
\begin{equation*}
(\widetilde{f}^{+}(t),\phi )=\lim_{s\rightarrow t+0}(\widetilde{f}%
^{+}(s),\phi ),\text{\qquad }\forall t\in \lbrack 0,T].
\end{equation*}%
Since the function $\widetilde{f}^{+}$ belongs to the class of equivalent
measurable functions of $\widetilde{f}$, we have the first relation in %
\eqref{2.13}. From now on, we consider the identification 
\begin{equation*}
\widetilde{f}(t)=\widetilde{f}^{+}(t),\qquad \forall t\in \lbrack 0,T].
\end{equation*}

For each $s\in (0,t)$, considering in \eqref{15:01:20:55e} the test function 
$\varphi (s,x,c)=\beta _{\epsilon }(s)\phi (x,c)$ with $\phi \in
C_{c}^{\infty }(\mathbb{R}^{2})$ and 
\begin{equation*}
\beta _{\epsilon }(s)=%
\begin{cases}
1,\quad 0\leqslant s\leqslant t, \\ 
1-\frac{s-t}{\epsilon },\quad t\leqslant s\leqslant t+\epsilon , \\ 
0,\quad s\geqslant t+\epsilon%
\end{cases}%
\end{equation*}%
we deduce 
\begin{eqnarray}
(\widetilde{f}(t),\varphi ) &=&\left( f_{0},\varphi (0)\right) +\int_{0}^{t}%
\left[ (\widetilde{f},c\partial _{x}\varphi )-(\partial _{x}\widetilde{p}\;%
\widetilde{f},\partial _{c}\varphi )\right] ds  \notag \\
&+&\sum_{k}\int_{0}^{t}\int_{\mathbb{R}^{2}}\varphi (s,x,c)\sigma
_{k}(x,c)\,d\widetilde{\nu }^{(t,x)}(c)dxd\widetilde{{\mathcal{W}}}{_{s}^{k}}
\notag \\
&+&\frac{1}{2}\int_{0}^{t}\int_{\mathbb{R}^{2}}\partial _{c}\varphi
(s,x,c)\sigma ^{2}(x,c)\,d\widetilde{\nu }^{(t,x)}(c)dxdt-\int_{0}^{t}\int_{%
\mathbb{R}^{2}}\partial _{c}\varphi \ d\widetilde{m}(s,x,c).  \label{f2}
\end{eqnarray}%
In particular, we obtain the second relation in \eqref{2.13}.

\hfill $\blacksquare $

\begin{remark}
\label{U1U} As a consequence of Lemma \ref{lemma2 copy(1)} we can understand
the transport equations \eqref{1b} for $\widetilde{f}$ in two equivalent
variational forms \eqref{15:01:20:55e} and \eqref{f2}, which we will use in
the next section.
\end{remark}

\subsection{The strong convergence of $\widetilde{u}_\protect\epsilon$}

\bigskip

The main task of this section consists in the verification that $\widetilde{f%
}=1_{\widetilde{u}>c}$. This identification will imply the strong
convergence of $\widetilde{u}_\epsilon$ to $\widetilde{u}$, as $\varepsilon
\rightarrow 0.$

First, we state the following lemma.

\begin{lemma}
\label{lemma4} There exists a random field $v=v(t,x),$ such that for $%
\widetilde{P}-$a.e. $\omega \in \widetilde{\Omega }$ 
\begin{equation}
\widetilde{f}(t,x,c)=1_{v(t,x)>c}\text{\ \ \ \ \ \ for a.e. \ }(t,x,c)\in 
\mathbb{R}_{T}^{2}.  \label{RR0}
\end{equation}
\end{lemma}

\textbf{Proof. } Let us consider 
\begin{equation*}
\varphi (x,y,c,\xi )=\phi (x,c)\psi (y,\xi )\quad \text{for }\phi ,\psi \in
C_{c}^{\infty }(\mathbb{R}^{2}).
\end{equation*}%
The It\^{o} formula, applied to the stochastic processes $\widetilde{f}$ and 
$\widehat{\widetilde{f}}(t)=1-\widetilde{f}(t)$ given by \eqref{1b}, yields%
\begin{align*}
d\left[ (\widetilde{f}(t),\phi )(\widehat{\widetilde{f}}(t),\psi )\right] &
=d(\widetilde{f}(t),\phi )(\widehat{\widetilde{f}}(t),\psi )+(\widetilde{f}%
(t),\phi )d(\widehat{\widetilde{f}}(t),\psi )+d(\widetilde{f}(t),\phi )d(%
\widehat{\widetilde{f}}(t),\psi ) \\
& =I_{0}+I_{p}+I_{W}+I_{\sigma }+I_{m}
\end{align*}%
where%
\begin{eqnarray*}
I_{0} &=&\int_{\mathbb{R}^{4}}\widetilde{f}(t,x,c)\widehat{\widetilde{f}}%
(t,y,\xi )\left[ c\partial _{x}\varphi +\xi \partial _{y}\varphi \right]
\,dxdcdyd\xi \,dt, \\
I_{p} &=&-\int_{\mathbb{R}^{4}}\widetilde{f}(t,x,c)\widehat{\widetilde{f}}%
(t,y,\xi )\left[ \partial _{x}\widetilde{p}(t,x)\partial _{c}\varphi
+\partial _{y}\widetilde{p}(t,y)\partial _{\xi }\varphi \right] \,dxdcdyd\xi
\,dt,
\end{eqnarray*}%
and%
\begin{align*}
I_{W}& =\sum_{k}\int_{\mathbb{R}^{4}}\sigma _{k}(x,c)\widehat{\widetilde{f}}%
(t,y,\xi )\varphi \,d\widetilde{\nu }^{(t,x)}(c)dxdyd\xi \, d{{\widetilde{%
\mathcal{W}}}_{t}} \\
& -\sum_{k}\int_{\mathbb{R}^{4}}\widetilde{f}(t,x,c)\sigma _{k}(y,\xi
)\varphi \,dxdcd\widetilde{\nu }^{(t,y)}(\xi )dyd{{\widetilde{\mathcal{W}}}%
_{t},} \\
I_{\sigma }& =\frac{1}{2}\int_{\mathbb{R}^{4}}\sigma ^{2}(x,c)\widehat{%
\widetilde{f}}(t,y,\xi )\partial _{c}\varphi \,d\widetilde{\nu }%
^{(t,x)}(c)dxdyd\xi dt \\
& -\frac{1}{2}\int_{\mathbb{R}^{4}}\widetilde{f}(t,x,c)\sigma ^{2}(y,\xi
)\partial _{\xi }\varphi \,dxdcd\widetilde{\nu }^{(t,y)}(\xi )dydt \\
& -\sum_{k}\int_{\mathbb{R}^{4}}\sigma _{k}(x,c)\sigma _{k}(y,\xi )\varphi
\,d\widetilde{\nu }^{(t,x)}(c)dxd\widetilde{\nu }^{(t,y)}(\xi )dydt, \\
I_{m}& =-\int_{\mathbb{R}^{4}}\,\widehat{\widetilde{f}}(t,y,\xi )\partial
_{c}\varphi \,d\widetilde{m}(t,x,c)dyd\xi +\int_{\mathbb{R}^{4}}\widetilde{f}%
(t,x,c)\partial _{\xi }\varphi \,dxdcd\widetilde{m}(t,y,\xi ).
\end{align*}%
Taking the expectation and integrating over the time interval $(0,t),$ we
have%
\begin{align}
& {\widetilde{\mathbb{E}}}\int_{\mathbb{R}^{4}}\widetilde{f}(t,x,c)\widehat{%
\widetilde{f}}(t,y,\xi )\varphi \,dxdydcd\xi ={\widetilde{\mathbb{E}}}\int_{%
\mathbb{R}^{4}}\widetilde{f}(0,x,c)\widehat{\widetilde{f}}(0,y,\xi )\varphi
\,dxdydcd\xi  \notag \\
& +{\widetilde{\mathbb{E}}}\int_{0}^{t}I_{0}(s)+{\widetilde{\mathbb{E}}}%
\int_{0}^{t}I_{p}(s)+{\widetilde{\mathbb{E}}}\int_{0}^{t}I_{\sigma }(s)+{%
\widetilde{\mathbb{E}}}\int_{0}^{t}I_{m}(s)  \notag \\
& =\mathbb{J}_{0}+\mathbb{J}_{p}+\mathbb{J}_{\sigma }+\mathbb{J}_{m}\quad 
\text{for a.e. }t\in \lbrack 0,T].  \label{J}
\end{align}%
By a density argument, this equality remains true for arbitrary test
function $\varphi \in C_{c}^{\infty }(\mathbb{R}^{4})$. In particular, we
can consider 
\begin{equation*}
\varphi (x,y,c,\xi )=\rho _{\epsilon }(x-y)\rho _{\delta }(c-\xi )\eta
_{n}(x)\eta _{k}(c)
\end{equation*}%
for small positive $\epsilon ,\delta $ and large natural $n,k$. Here $\rho
_{\epsilon }(s)=\frac{1}{\epsilon }\rho \left( \frac{s}{\epsilon }\right) $
and $\rho $ is a smooth nonnegative function with a support on $[-1,1]$,
such that $\int_{\mathbb{R}}\rho (s)ds=1,$ and $\eta _{n}$ is a smooth
nonnegative one, satisfying 
\begin{equation}
\eta _{n}(z)=%
\begin{cases}
1\quad \text{for}\quad |z|\leqslant n, \\ 
0\quad \text{for}\quad |z|\geqslant n+1, \\ 
|\eta _{n}^{\prime }(z)|\leqslant 2,\quad \forall z\in \mathbb{R}.%
\end{cases}
\label{eta}
\end{equation}%
The test function $\varphi $ fulfills the relations%
\begin{equation}
\partial _{x}\varphi +\partial _{y}\varphi =\rho _{\epsilon }\rho _{\delta
}\eta _{n}^{\prime }\eta _{k},\qquad \partial _{c}\varphi +\partial _{\xi
}\varphi =\rho _{\epsilon }\rho _{\delta }\eta _{n}\eta _{k}^{\prime }.
\label{08:11:31:01}
\end{equation}

Firstly, applying the first identity of \eqref{08:11:31:01}, we obtain%
\begin{align*}
\mathbb{J}_{0}& \leqslant {\widetilde{\mathbb{E}}}\int_{0}^{t}\int_{\mathbb{R%
}^{4}}\widetilde{f}(s,x,c)\widehat{\widetilde{f}}(s,y,\xi )\left( c-\xi
\right) \partial _{x}\varphi \,dxdydcd\xi ds \\
& +{\widetilde{\mathbb{E}}}\int_{0}^{t}\int_{\mathbb{R}^{4}}\widetilde{f}%
(s,x,c)\widehat{\widetilde{f}}(s,y,\xi )\ \xi \ \rho _{\epsilon }\rho
_{\delta }\eta _{n}^{\prime }\eta _{k}\,dxdydcd\xi ds \\
& =\mathbb{J}_{0,1}^{\epsilon ,\delta }+\mathbb{J}%
_{0,2}^{\epsilon ,\delta },
\end{align*}%
We have that there exists a positive constant $C_{n,k},$ being independent
of $\epsilon ,\delta $ (and depending of $n,k$), such that 
\begin{equation}
\mathbb{J}_{0,1}^{\epsilon ,\delta }\leqslant C_{n,k}\frac{\delta }{\epsilon 
}.  \label{j01}
\end{equation}%
Applying standard techniques of Kruzkov \cite{Kruzkov} (see also this
approach in details in Theorem 3 on pages 650-653 of \cite{Evans}), we can
derive 
\begin{equation}
lim_{\epsilon ,\delta \rightarrow 0}\mathbb{J}_{0,2}^{\epsilon ,\delta }={%
\widetilde{\mathbb{E}}}\int_{0}^{t}\int_{\mathbb{R}^{2}}\widetilde{f}(s,x,c)%
\widehat{\widetilde{f}}(s,x,c)\ c\ \eta _{n}^{\prime }\eta _{k}\,dxdcds=%
\mathbb{K}_{0,2}^{n,k}.  \label{j02}
\end{equation}%
The second identity of \eqref{08:11:31:01} gives  
\begin{align*}
\mathbb{J}_{p}& =-{\widetilde{\mathbb{E}}}\int_{0}^{t}\int_{\mathbb{R}^{4}}%
\widetilde{f}(s,x,c)\widehat{\widetilde{f}}(s,y,\xi )\left( \partial _{x}%
\widetilde{p}(s,x)-\partial _{y}\widetilde{p}(s,y)\right) \,\partial
_{c}\varphi \,dxdcdyd\xi ds \\
& -{\widetilde{\mathbb{E}}}\int_{0}^{t}\int_{\mathbb{R}^{4}}\widetilde{f}%
(s,x,c)\widehat{\widetilde{f}}(s,y,\xi )\partial _{y}\widetilde{p}(s,y)\rho
_{\epsilon }\rho _{\delta }\eta _{n}\eta _{k}^{\prime }\,dxdydcd\xi ds \\
& =\mathbb{J}_{p,1}+\mathbb{J}_{p,2}.
\end{align*}%
The embedding $W^{1,q}(\mathbb{R})\hookrightarrow C^{\alpha }(\mathbb{R})$
with $\alpha =1-\frac{1}{q}$ and the estimate \eqref{Z1}$_{2}$ imply 
\begin{equation*}
{\widetilde{\mathbb{E}}}\sup_{t\in \lbrack 0,T]}\Vert \partial _{x}%
\widetilde{p}\Vert _{C^{\alpha }}={\widetilde{\mathbb{E}}}\sup_{t\in \lbrack
0,T]}\left( \sup_{x,y\in \mathbb{R}}\frac{|\partial _{x}\widetilde{p}%
(t,x)-\partial _{x}\widetilde{p}(t,y)|}{|x-y|^{\alpha }}\right) \leqslant C{%
\widetilde{\mathbb{E}}}\sup_{t\in \lbrack 0,T]}\Vert \widetilde{p}\Vert
_{W^{2,q}}.
\end{equation*}%
Therefore integrating by parts on the parameter $c$ and using%
\begin{equation*}
\int_{0}^{t}\int_{\mathbb{R}^{4}}\widehat{\widetilde{f}}(s,y,\xi )\,\varphi
\,d\widetilde{\nu }^{(s,x)}(c)dxdyd\xi ds\leqslant C_{n,k},
\end{equation*}%
we obtain%
\begin{align}
\mathbb{J}_{p,1}& ={\widetilde{\mathbb{E}}}\int_{0}^{t}\int_{\mathbb{R}^{4}}%
\widehat{\widetilde{f}}(s,y,\xi )\left( \partial _{y}\widetilde{p}%
(s,x)-\partial _{x}\widetilde{p}(s,y)\right) \varphi \,d\widetilde{\nu }%
^{(s,x)}(c)dxdyd\xi ds  \notag \\
& \leqslant C_{n,k}{\widetilde{\mathbb{E}}}\sup_{t\in \lbrack 0,T]}\Vert 
\widetilde{p}\Vert _{W^{2,q}}\ \epsilon ^{\alpha }\leqslant C_{n,k}\epsilon
^{\alpha }.  \label{jp1}
\end{align}%
By the techniques of Kruzkov, we have 
\begin{equation}
lim_{\epsilon ,\delta \rightarrow 0}\mathbb{J}_{p,2}^{\epsilon ,\delta }=-{%
\widetilde{\mathbb{E}}}\int_{0}^{t}\int_{\mathbb{R}^{2}}\widetilde{f}(s,x,c)%
\widehat{\widetilde{f}}(s,x,c)\partial _{x}\widetilde{p}(s,x)\eta _{n}\eta
_{k}^{\prime }\,dxdcds=\mathbb{K}_{p,2}^{n,k}.  \label{jp2}
\end{equation}

By \eqref{08:11:31:01}$_{2}$, the integration by parts in $I_{\sigma }$ with
the help of the second line in \eqref{ff0}, we deduce%
\begin{equation*}
\mathbb{J}_{\sigma }=\mathbb{J}_{\sigma ,1}^{\epsilon ,\delta }+\mathbb{J}%
_{\sigma ,2}^{\epsilon ,\delta }
\end{equation*}%
with 
\begin{equation*}
\mathbb{J}_{\sigma ,1}^{\epsilon ,\delta }=\frac{1}{2}\sum_{k}{\widetilde{%
\mathbb{E}}}\int_{0}^{t}\int_{\mathbb{R}^{4}}|\sigma _{k}(x,c)-\sigma
_{k}(y,\xi )|^{2}\varphi \,d\widetilde{\nu }^{(s,x)}(c)dxd\widetilde{\nu }%
^{(s,y)}(\xi )dy\,ds
\end{equation*}%
and 
\begin{eqnarray*}
\mathbb{J}_{\sigma ,2}^{\epsilon ,\delta } &=&\frac{1}{2}{\widetilde{\mathbb{%
E}}}\int_{0}^{t}\int_{\mathbb{R}^{2}}\left( \int_{\mathbb{R}^{2}}\sigma
^{2}(x,c)\widehat{\widetilde{f}}(s,y,\xi )\rho _{\epsilon }\rho _{\delta
}\eta _{n}\eta _{k}^{\prime }dyd\xi \right) \,d\widetilde{\nu }%
^{(s,x)}(c)dxds \\
&-&\frac{1}{2}{\widetilde{\mathbb{E}}}\int_{0}^{t}\int_{\mathbb{R}%
^{2}}\left( \int_{\mathbb{R}^{2}}\widetilde{f}(s,x,c)\sigma ^{2}(y,\xi )\rho
_{\epsilon }\rho _{\delta }\eta _{n}\eta _{k}^{\prime }\,dxdc\right) d%
\widetilde{\nu }^{(s,y)}(\xi )dyds.
\end{eqnarray*}%
Applying the assumption \eqref{20:01:21:59} in the equality for $\mathbb{J}%
_{\sigma ,1}^{\epsilon ,\delta },$ we obtain  
\begin{eqnarray}
\mathbb{J}_{\sigma ,1}^{\epsilon ,\delta } &\leqslant &C{\widetilde{\mathbb{E%
}}}\int_{0}^{t}\int_{\mathbb{R}^{4}}\left( |x-y|^{2}+|c-\xi |h(|c-\xi
|)\right) \varphi \,d\widetilde{\nu }^{(s,x)}(c)dxd\widetilde{\nu }%
^{(s,y)}(\xi )dy\,ds  \notag \\
&=&C{\widetilde{\mathbb{E}}}\int_{0}^{t}\int_{\mathbb{R}^{2}}\left( \int_{%
\mathbb{R}^{2}}\left( |x-y|^{2}+|c-\xi |h(|c-\xi |)\right) \rho _{\delta
}\eta _{k}\,d\widetilde{\nu }^{(s,x)}(c)\widetilde{\nu }^{(s,y)}(\xi
)\right) \rho _{\epsilon }\eta _{n}dxdy\,ds  \notag \\
&\leqslant &C_{k}{\widetilde{\mathbb{E}}}\int_{0}^{t}\int_{\mathbb{R}%
^{2}}\left( |x-y|^{2}\delta ^{-1}+h(\delta )\right) \rho _{\epsilon
}(x-y)\eta _{n}(x)dxdy\,ds  \notag \\
&\leqslant &C_{n,k}\left( \epsilon ^{2}\delta ^{-1}+h(\delta )\right) .
\label{lll}
\end{eqnarray}%
By the techniques of Kruzkov, we also have 
\begin{eqnarray}
lim_{\epsilon ,\delta \rightarrow 0}\mathbb{J}_{\sigma ,2}^{\epsilon ,\delta
} &=&\frac{1}{2}{\widetilde{\mathbb{E}}}\int_{0}^{t}\int_{\mathbb{R}%
^{2}}\sigma ^{2}(x,c)(1-2\widetilde{f}(s,x,c))\eta _{n}\eta _{k}^{\prime }\,d%
\widetilde{\nu }^{(s,x)}(c)dxds  \notag \\
&=&\mathbb{K}_{\sigma ,2}^{n,k}.  \label{ks1}
\end{eqnarray}

By the second identity of \eqref{08:11:31:01} and the integration by parts
in $I_{m}$ with the help of \eqref{ff0}, we deduce%
\begin{equation*}
\mathbb{J}_{m}=\mathbb{J}_{m,1}^{\epsilon ,\delta }+\mathbb{J}%
_{m,2}^{\epsilon ,\delta }
\end{equation*}%
with 
\begin{eqnarray*}
\mathbb{J}_{m,1}^{\epsilon ,\delta } &=&{\widetilde{\mathbb{E}}}%
\int_{0}^{t}\int_{\mathbb{R}^{4}}\,\widehat{\widetilde{f}}(s,y,\xi )\partial
_{\xi }\varphi \,d\widetilde{m}(s,x,c)dyd\xi -{\widetilde{\mathbb{E}}}%
\int_{0}^{t}\int_{\mathbb{R}^{4}}\widetilde{f}(s,x,c)\partial _{c}\varphi
\,dxdcd\widetilde{m}(s,y,\xi ) \\
&=&-{\widetilde{\mathbb{E}}}\int_{0}^{t}\int_{\mathbb{R}^{4}}\,\varphi \,d%
\widetilde{\nu }^{(s,y)}(\xi )dyd\widetilde{m}(s,x,c)-{\widetilde{\mathbb{E}}%
}\int_{0}^{t}\int_{\mathbb{R}^{4}}\varphi \,d\widetilde{\nu }^{(s,x)}(c)dxd%
\widetilde{m}(s,y,\xi )\leqslant 0
\end{eqnarray*}%
and 
\begin{eqnarray*}
\mathbb{J}_{m,2}^{\epsilon ,\delta } &=&-{\widetilde{\mathbb{E}}}%
\int_{0}^{t}\int_{\mathbb{R}^{2}}\,\left( \int_{\mathbb{R}^{2}}\widehat{%
\widetilde{f}}(s,y,\xi )\rho _{\epsilon }\rho _{\delta }\eta _{n}\eta
_{k}^{\prime }dyd\xi \right) \,d\widetilde{m}(s,x,c) \\
&&+{\widetilde{\mathbb{E}}}\int_{0}^{t}\int_{\mathbb{R}^{2}}\left( \int_{%
\mathbb{R}^{2}}\widetilde{f}(s,x,c)\rho _{\epsilon }\rho _{\delta }\eta
_{n}\eta _{k}^{\prime }\,dxdc\right) d\widetilde{m}(s,y,\xi ),
\end{eqnarray*}%
such that 
\begin{equation}
lim_{\epsilon ,\delta \rightarrow 0}\mathbb{J}_{m,2}^{\epsilon ,\delta }=-{%
\widetilde{\mathbb{E}}}\int_{0}^{t}\int_{\mathbb{R}^{2}}\,(2\widetilde{f}%
(s,x,c)-1)\eta _{n}\eta _{k}^{\prime }d\widetilde{m}(t,x,c)=\mathbb{K}%
_{m,2}^{n,k}  \label{km1}
\end{equation}%
by standard techniques of Kruzkov.

Obviously%
\begin{equation*}
lim_{\epsilon ,\delta \rightarrow 0}{\widetilde{\mathbb{E}}}\int_{\mathbb{R}%
^{4}}\widetilde{f}(t,x,c)\widehat{\widetilde{f}}(t,y,\xi )\varphi
\,dxdydcd\xi ={\widetilde{\mathbb{E}}}\int_{\mathbb{R}^{2}}\widetilde{f}%
(t,x,c)\widehat{\widetilde{f}}(t,x,c)\eta _{n}\eta _{k}\,dxdc
\end{equation*}%
and 
\begin{equation*}
lim_{\epsilon ,\delta \rightarrow 0}{\widetilde{\mathbb{E}}}\int_{\mathbb{R}%
^{4}}\widetilde{f}(0,x,c)\widehat{\widetilde{f}}(0,y,\xi )\varphi
\,dxdydcd\xi ={\widetilde{\mathbb{E}}}\int_{\mathbb{R}^{4}}\widetilde{f}%
(0,x,c)\widehat{\widetilde{f}}(0,x,c)\eta _{n}\eta _{k}\,dxdc=0,
\end{equation*}%
since $\widetilde{f}(0)\widehat{\widetilde{f}}(0)=(f_{0}-f_{0}^{2})=0$ a.e.
in $\widetilde{\Omega }\times \mathbb{R}^{2}.$

Therefore choosing $\delta =\epsilon ^{3/2}$ in all above considerations and
using \eqref{j01}-\eqref{km1}, the limit transition on $\epsilon \rightarrow
0$ in the equality \eqref{J} gives%
\begin{equation}
{\widetilde{\mathbb{E}}}\int_{\mathbb{R}^{2}}\widetilde{f}(t,x,c)\widehat{%
\widetilde{f}}(t,x,c)\eta _{n}\eta _{k}\,dxdc\leqslant \mathbb{K}%
_{0,2}^{n,k}+\mathbb{K}_{p,2}^{n,k}+\mathbb{K}_{\sigma ,2}^{n,k}+\mathbb{K}%
_{m,2}^{n,k}.  \label{f}
\end{equation}

Now let us show that all terms in the right hand side of the last inequality
can be estimated by some integrable functions on $\widetilde{\Omega }\times 
\mathbb{R}_{T}^{2}$, being independent on $n,k$ parameters. This permits to
use the dominated convergence theorem and pass in these terms on the limit
transition on $n,k\rightarrow +\infty .$

Taking into account that $\widehat{\widetilde{f}}(t)=1-\widetilde{f}(t)$, $%
0\leqslant \widetilde{f}\leqslant 1$ and the definition \eqref{eta} for $%
\eta _{n},\eta _{k}$, we have that%
\begin{equation*}
|\mathbb{K}_{0,2}^{n,k}|\leqslant 2{\widetilde{\mathbb{E}}}\int_{\mathbb{R}%
_{T}}\left( \int_{-\infty }^{0}(1-\widetilde{f}(t,x,c))\ |c|\
dc+\int_{0}^{+\infty }\widetilde{f}(t,x,c)\ c\ \,dc\right) dtdx\leqslant C
\end{equation*}%
by (\ref{eq2.88})$_{1}$ for $p=2$. Since $\eta _{n}^{\prime }\eta
_{k}\rightarrow 0$ on $\mathbb{R}^{2}$ as $n,k\rightarrow +\infty ,$ the
dominated convergence theorem implies 
\begin{equation}
lim_{n,k\rightarrow \infty }\mathbb{K}_{0,2}^{n,k}=0.  \label{k02}
\end{equation}

There exist constants $C$, independent of $n,k$, such that 
\begin{eqnarray*}
|\mathbb{K}_{p,2}^{n,k}| &\leqslant &C{\widetilde{\mathbb{E}}}%
\int_{0}^{T}\int_{\mathbb{R}^{2}}\widetilde{f}(t,x,c)\widehat{\widetilde{f}}%
(t,x,c)|\partial _{x}\widetilde{p}(t,x)|\,dxdcdt \\
&\leqslant &C{\widetilde{\mathbb{E}}}\int_{\mathbb{R}_{T}}\left(
\int_{-\infty }^{0}(1-\widetilde{f}(t,x,c))\ dc+\int_{0}^{+\infty }%
\widetilde{f}(t,x,c)\ dc\right) |\partial _{y}\widetilde{p}(t,x)|\,dxds \\
&\leqslant &C{\widetilde{\mathbb{E}}}\int_{\mathbb{R}_{T}}\left[
(\int_{-\infty }^{0}(1-\widetilde{f}(t,x,c))\ dc)^{2}+(\int_{0}^{+\infty }%
\widetilde{f}(t,x,c)\ dc)^{2}+|\partial _{y}\widetilde{p}(t,x)|^{2}\right]
dxdt \\
&\leqslant &C
\end{eqnarray*}%
by Young's inequality for products and \eqref{Z1}$_{1},$ (\ref{eq2.88})$_{2}$
for $p=2$. Hence we can use the dominated convergence theorem to derive 
\begin{equation}
lim_{n,k\rightarrow \infty }\mathbb{K}_{p,2}^{n,k}=0.  \label{kp2}
\end{equation}

For the terms $\mathbb{K}_{\sigma ,2}^{n,k}$ and $\mathbb{K}_{m,2}^{n,k},$
we have 
\begin{equation*}
|\mathbb{K}_{\sigma ,2}^{n,k}|\leqslant C{\widetilde{\mathbb{E}}}%
\int_{0}^{T}\int_{\mathbb{R}^{2}}c^{2}\,d\widetilde{\nu }^{(t,x)}(c)dxdcdt%
\leqslant C
\end{equation*}%
by \eqref{ass1}, \eqref{eq2.88}$_{3}$ and 
\begin{equation*}
|\mathbb{K}_{m,2}^{n,k}|\leqslant C{\widetilde{\mathbb{E}}}\int_{0}^{T}\int_{%
\mathbb{R}^{2}}1d\widetilde{m}(t,x,c)\leqslant C.
\end{equation*}%
by \eqref{eq2.88}$_{4}$ for $s=1$ for some constants $C$, independent of $%
n,k $. Hence%
\begin{equation}
lim_{n,k\rightarrow \infty }\mathbb{K}_{\sigma ,2}^{n,k}=0\quad \text{and}%
\quad lim_{n,k\rightarrow \infty }\mathbb{K}_{m,2}^{n,k}=0.  \label{km2}
\end{equation}

Finally taking the limit transition on $n,k\rightarrow +\infty $ in the
inequality \eqref{f} and using \eqref{k02}-\eqref{km2}, we obtain 
\begin{equation*}
0\leqslant {\widetilde{\mathbb{E}}}\int_{\mathbb{R}^{2}}\widetilde{f}(1-%
\widetilde{f})dxdc\leqslant 0,\quad \text{that is}\quad \widetilde{f}(1-%
\widetilde{f})=0\quad \text{on \ }\widetilde{\Omega }\times \mathbb{R}%
_{T}^{2},
\end{equation*}%
which gives the existence a random field $v=v(t,x)$ satisfing \eqref{RR0}%
.\hfill $\blacksquare $

\bigskip

As a consequence of the previous considerations, we are able to solve the
main task of this section.

\begin{lemma}
\label{Lemma 4.44} For \ a.s. $\omega \in \widetilde{\Omega }$ \ the
stochastic process 
\begin{equation*}
\widetilde{u}\in L^{2}(\widetilde{\Omega };L^{\infty }(0,T;L^{2}(\mathbb{R}%
)\cap L^{2q}(\mathbb{R}))),
\end{equation*}%
defined in (\ref{limit0})$_{2,3}$ verifies 
\begin{equation}
v=\widetilde{u}\quad \text{a.e. in}\quad \mathbb{R}_{T},\quad \widetilde{f}%
=1_{\widetilde{u}>c},\quad \widetilde{\nu }^{(t,x)}(c)=\delta _{\widetilde{u}%
(t,x)=c}\quad \text{a.e. in}\quad \mathbb{R}_{T}^{2}  \label{24}
\end{equation}%
and 
\begin{equation}
\widetilde{u}_{\varepsilon }^{2}\rightharpoonup \widetilde{u}^{2}\quad \text{%
in \ }\mathcal{D}^{\prime }(\mathbb{R}_{T})\text{$.$}  \label{vv}
\end{equation}
\end{lemma}

\textbf{Proof}. Let us consider an arbitrary $p\in \lbrack 1,2q)$ and a
positive test function $\varphi =\varphi (t,x)\in C_{c}^{\infty }(\mathbb{R}%
_{T}).$ We have the estimate 
\begin{eqnarray*}
\int_{\mathbb{R}_{T}}\int_{0}^{+\infty }\varphi \widetilde{f}_{\varepsilon
}pc^{p-1}dtdxdc &=&\int_{\mathbb{R}_{T}}\varphi |\widetilde{u}_{\varepsilon
}|_{+}^{p}dtdx \\
&\leqslant &||\widetilde{u}_{\varepsilon }||_{L^{\infty }(0,T;L^{2q}(\mathbb{%
R}))}^{p}\Vert \varphi \Vert _{L^{1}(0,T;L^{r}(\mathbb{R}))}\leqslant
C(p,\varphi )
\end{eqnarray*}%
with $r,$ such that $1=\frac{p}{2q}+\frac{1}{r},$ and the constant $%
C(p,\varphi )$ depending only on $p$ and $\varphi $.

By the weak convergence of $\widetilde{f}_{\varepsilon }$ to $\widetilde{f}$
in $L^{\infty }(\mathbb{R}_{T}^{2}),$ \ for any $M>0$ and any $\delta >0$,
there exists $\varepsilon _{1}=\varepsilon _{1}(M,\delta )>0,$ such that for
each positive $\varepsilon <\varepsilon _{1}$, we have 
\begin{equation}
\left\vert \int_{\mathbb{R}_{T}}\int_{0}^{M}\varphi \widetilde{f}%
_{\varepsilon }pc^{p-1}dtdxdc-\int_{\mathbb{R}_{T}}\int_{0}^{M}\varphi 
\widetilde{f}pc^{p-1}dtdxdc\right\vert \leqslant \frac{\delta }{3}.
\label{aa66}
\end{equation}%
If we take $\delta =1$, we obtain that 
\begin{equation*}
\int_{\mathbb{R}_{T}}\int_{0}^{M}\varphi \widetilde{f}pc^{p-1}dtdxdc%
\leqslant \frac{1}{3}+\int_{\mathbb{R}_{T}}\int_{0}^{+\infty }\varphi 
\widetilde{f}_{\varepsilon }pc^{p-1}dtdxdc\leqslant \frac{1}{3}+C(p,\varphi
)\quad \text{for all \ }M,
\end{equation*}%
which implies 
\begin{equation}
\int_{\mathbb{R}_{T}}\int_{0}^{+\infty }\varphi \widetilde{f}%
pc^{p-1}dtdtdc<\infty  \label{aa44}
\end{equation}%
by the monotone convergence theorem.

For any $M>0$, using the H\"{o}lder inequality and the estimates (\ref{Z}),
\ we verify that 
\begin{equation}
\left\vert \int_{\mathbb{R}_{T}}\int_{M}^{+\infty }\varphi \widetilde{f}%
_{\varepsilon }pc^{p-1}dtdxdc\right\vert \leqslant \frac{C\Vert \varphi
\Vert _{L^{\infty }(\mathbb{R}_{T})}}{M^{q-p}}\int_{|\widetilde{u}%
_{\varepsilon }|_{+}\geqslant M}|\widetilde{u}_{\varepsilon
}|_{+}^{q}dtdx\leqslant \frac{C(\varphi ,p)}{M^{q-p}}.  \label{aa11}
\end{equation}%
With the help of the inequalities \eqref{aa66}-\eqref{aa11} and considering $%
M$ large enough, we deduce 
\begin{equation*}
\left\vert \int_{\mathbb{R}_{T}}\int_{0}^{+\infty }\varphi \widetilde{f}%
_{\varepsilon }pc^{p-1}dtdxdc-\int_{\mathbb{R}_{T}}\int_{0}^{+\infty
}\varphi \widetilde{f}pc^{p-1}dtdtdc\right\vert \leqslant \delta ,
\end{equation*}%
which corresponds to 
\begin{equation}
|\widetilde{u}_{\varepsilon }|_{+}^{p}=\int_{0}^{+\infty }\widetilde{f}%
_{\varepsilon }(\mathbf{\cdot },\mathbf{\cdot },c)pc^{p-1}\
dc\rightharpoonup |v|_{+}^{p}\quad \text{in \ }\mathcal{D}^{\prime }(\mathbb{%
R}_{T})  \label{a0}
\end{equation}%
by \eqref{RR0}. Dealing with the function $1-\widetilde{f}_{\varepsilon }$
instead of $\widetilde{f}_{\varepsilon }$, and applying the same reasoning,
we can derive 
\begin{equation}
|\widetilde{u}_{\varepsilon }|_{-}^{p}=\int_{-\infty }^{0}(1-\widetilde{f}%
_{\varepsilon }(\mathbf{\cdot },\mathbf{\cdot },c))p|c|^{p-1}\
dc\rightharpoonup |v|_{-}^{p}\quad \text{in \ }\mathcal{D}^{\prime }(\mathbb{%
R}_{T}).  \label{a00}
\end{equation}%
\hfill

Now let us consider 2 distinct cases:

1. Let us consider $p=1.$ Having the weak convergence of $\widetilde{u}%
_{\varepsilon }=|\widetilde{u}_{\varepsilon }|_{+}-|\widetilde{u}%
_{\varepsilon }|_{+}$ to $\widetilde{u}$ in $L^{\infty }(0,T;L^{2}(\mathbb{R}%
)\cap L^{2q}(\mathbb{R}))$ by \eqref{limit0}$_{2,3},$ we see that the
convergences \eqref{a0}-\eqref{a00} imply \eqref{24}.

2. Now let $p=2$. Due to the estimates (\ref{Z}), \eqref{a0}-\eqref{a00} and 
$|z|^{2}=|z|_{+}^{2}+|z|_{-}^{2}$, we conclude (\ref{vv}).\hfill $%
\blacksquare $

\bigskip \bigskip

\textbf{Proof of Theorem \ref{the_1}.} The relations \eqref{f2} and %
\eqref{24} imply that the pair $\left( \widetilde{u},\widetilde{\mathrm{p}}%
\right) $\ satisfies the transport equation \eqref{10:36:25:03:210}. Also
the convergences \eqref{limit0} and (\ref{vv}) imply that the limit pair $(%
\widetilde{u},\widetilde{\mathrm{p}})$ satisfies the equation %
\eqref{10:36:25:03:21}. Therefore, $\left( \widetilde{u},\widetilde{\mathrm{p%
}}\right) $ is a weak kinetic solution of \ the system \eqref{Hyp3}\ in the
sense of Definition \ref{weak_kinetic}.

Finally, taking the integration by parts on $c\in \mathbb{R}$ in %
\eqref{10:36:25:03:210}, a similar argumentation as under deduction of %
\eqref{15:01:20:550} implies the validality of \eqref{15:01:20:55}. \ \hfill 
$\blacksquare $

\appendix

\section{Appendix}

\label{appendix}

\setcounter{equation}{0}

\bigskip

\subsection{Young measures}

We consider the $\sigma -$finite measure space $(\Omega \times \mathbb{R}%
_{T},P\otimes \lambda )$, where $\lambda $ stands for the Lebesgue measure
on $\mathbb{R}_{T}$. Let us denote by $\mathcal{M}_{Y}\left( \Omega \times 
\mathbb{R}_{T}\right) $ the space of Young measures on $(\Omega \times 
\mathbb{R}_{T},P\otimes \lambda )$ which vanish at infinity. We recall that
a Young measure $\nu \in \mathcal{M}_{Y}\left( \Omega \times \mathbb{R}%
_{T}\right) $ is a function defined on $\Omega \times \mathbb{R}_{T}$ with
values in the space $\mathcal{P}_{1}(\mathbb{R})$ of the probability
measures on $\mathbb{R}$, which verifies 
\begin{equation*}
\mathbb{E}\int_{0}^{T}\int_{\mathbb{R}^{2}}|c|^{p}d\nu
^{(t,x)}(c)dxdt<\infty ,\qquad \forall 1\leqslant p<+\infty .
\end{equation*}

For $1\leqslant p<+\infty $, let us set 
\begin{equation}
|\nu |_{\infty }^{p}=\mathbb{E}\int_{0}^{T}\int_{\mathbb{R}^{2}}|c|^{p}d\nu
^{(t,x)}(c)dxdt.  \label{24:01:14:36}
\end{equation}%
It is known that for each $R>0$ and $1\leqslant p<+\infty $, the set 
\begin{equation}
K_{R}^{p}=\{\nu \in \mathcal{M}_{Y}\left( \Omega \times \mathbb{R}%
_{T}\right) :\;|\nu |_{\infty }^{p}<R\}  \label{24:01:14:37}
\end{equation}%
is a compact set. We refer to \cite{B94} for a more detailed description of
the Young measures.

\subsection{Compactness results}

\setcounter{equation}{0}

In this subsection we recall some compactness results, used in the article.

First we present some compactness results of probabilistic nature. Let $%
(S,\rho )$ be a topological space and $\mathcal{B}(S)$ be its Borel $\sigma
- $algebra.

\begin{definition}
\label{prob compact} (See in \cite{B94}, \cite{EG}) a) A family $\{\mu
_{n}\}_{n\in \mathbb{N}}$ of probability measures on $(S,\mathcal{B}(S))$ is
relatively compact if every sequence of elements of $\{\mu _{n}\}_{n\in 
\mathbb{N}}$ contains a subsequence $\{\mu _{n_{i}}\}_{i\in \mathbb{N}}$
which converges weakly to a probability measure $\mu $ that is, for any
bounded and continuous function $\varphi $ on $S$, we have%
\begin{equation*}
\int \varphi (s)d\mu _{_{n_{i}}}(s)\rightarrow \int \varphi (s)d\mu (s).
\end{equation*}%
b) A family $\left\{ \mu _{n}\right\} _{n\in \mathbb{N}}$ is said to be
uniformly tight if, for any $\varepsilon >0$, there exists a compact set $%
K_{\varepsilon }\subset S$ such that $\mu _{n}(K_{\varepsilon })\geqslant
1-\varepsilon $ for all $n\in \mathbb{N}$.
\end{definition}

\bigskip

We frequently use the following result of Prokhorov \cite{Prokhorov}.

\begin{theorem}
\label{prokhorov} If $S$ is a separable and complete metric space (Polish
space), a family $\{\mu _{n}\}_{n\in \mathbb{R}}$ of probability measures on 
$(S,\mathcal{B}(S))$ is relatively compact if and only if it is uniformly
tight.
\end{theorem}

Let us stress that the application of the Prohorov theorem as well as the
Skorohod theorem \cite{Skorokhod} require a separable and complete metric
space (Polish space). Nevertheless we deal with some spaces which are not
separable and others that are not metrizable.

Therefore, we consider a version of the Skorokhod established by Jakubowski 
\cite{J98}.

\begin{theorem}[Jakubowski-Skorokhod representation theorem]
\label{Jakubowski} Let $(\mathcal{X},\rho )$ be a topological space.

\begin{enumerate}
\item[(J)] Assume that there exists a countable set 
\begin{equation*}
\{\Psi _{i}:\mathcal{X}\rightarrow \lbrack -1,1],\,i\in \mathbb{N}\}
\end{equation*}%
of continuous functions, which separate points of $\mathcal{X}$, i e. if for 
$x,y\in \mathcal{X}$ we have the equalities $\Psi _{i}(x)=\Psi
_{i}(y),\;\forall i\in \mathbb{N}$, then $x=y.$
\end{enumerate}

Consider the measurable space $(\mathcal{X},\mathcal{G}_{\Psi })$, where $%
\mathcal{G}_{\Psi }$ denotes the $\sigma $-algebra on $\mathcal{X}$
generated by the sequence $(\Psi _{i})_{i\in \mathbb{N}}.$

Let $\{\mu _{n}\}_{n\in \mathbb{N}}$ be a uniformly tight sequence of laws
on $\mathcal{X}$ defined on $(\mathcal{X},\mathcal{G}_{\Psi })$. Then there
exists a subsequence $\mu _{n_{1}},\mu _{n_{2}},\dots $ which admits the
strong a.s. Skorokhod representation on $([0,1],\mathcal{B}([0,1]),d\omega )$%
. Namely, there exists a sequence $X,X_{n_{1}},X_{n_{2}},\dots $ of $%
\mathcal{X}$-valued random variables defined on the probability space $%
([0,1],\mathcal{B}([0,1]),d\omega )$ such that the law of $X_{n_{i}}$ is $%
\mu _{n_{i}}$, $i=1,2,\dots $, the law of $X$ is a radom measure, and%
\begin{equation*}
X_{n_{k}}(\omega )\rightarrow X(\omega ),\quad \forall \omega \in \lbrack
0,1].
\end{equation*}
\end{theorem}

A topological space satisfying the condition $(J)$ is called a sub-Polish
space.

\begin{remark}
\label{Remark-Jakubowski} Consider the case $\mathcal{X}=E^{\ast }$, where $%
E^{\ast }$ is the dual space of a separable Banach space $E$, and take the
ball of radius $R$, with respect to the strong topology of $E^{\ast }$, i.
e. 
\begin{equation*}
B_{R}=\{f\in E:\,\Vert f\Vert _{E^{\ast }}\leqslant R\}.
\end{equation*}
For each $R$, the ball $B_{R}$ endowed with the weak-* topology (induced by
the weak-* topology on all space) is a sub-Polish space.

In fact, considering a countable dense subset $\{x_{i},i\in \mathbb{N}\}$ of 
$E$, we can set a countable family 
\begin{equation*}
\left\{ \Psi _{i}:B_{R}\rightarrow \lbrack -1,1],\quad \Psi _{i}(f)=\frac{%
\langle f,x_{i}\rangle }{R\Vert x_{i}\Vert _{E}},\;i\in \mathbb{N}\right\}
\end{equation*}%
of continuous functions (for the weak-* topology), which separate points of $%
E^{\ast }.$
\end{remark}

\bigskip

Now let us present the result being a version of Theorem 5 of \cite{S87}.

\begin{lemma}
\label{simon}Let $Y$ be a bounded open subset of $\ \mathbb{R}^{2}.$ Let us
consider 
\begin{align*}
\mathcal{Z}=& \{z\in L^{\infty }(0,T;L^{2}(Y)): \\
& \qquad \delta ^{-1}\sup_{0\leqslant \theta \leqslant \delta
}\int_{0}^{T-\delta }||z(t+\theta )-z(t)||_{H^{-3}(B_{n})}^{2}dt\leqslant
C,\quad \forall \delta >0\},
\end{align*}%
where $C>0$ is a constant, independent of $\delta .$ The set $\mathcal{Z}$
is compact in $C(0,T;H^{-1}(Y))$.
\end{lemma}

\bigskip

\textbf{Acknowledgments.} {\ The authors would like to thank the anonymous
Referee for relevant recommendations which contributed to improve the
article. }

The work of N.V. Chemetov was supported by FAPESP (Funda\c{c}\~{a}o de
Amparo \`{a} Pesquisa do Estado de S\~{a}o Paulo), project 2021/03758-8,
"Mathematical problems in fluid dynamics".

The work of F. Cipriano is funded by national funds through the FCT - Funda%
\c{c}\~{a}o para a Ci\^{e}ncia e a Tecnologia, I.P., under the scope of the
projects UIDB/00297/2020 and UIDP/00297/2020 (Center for Mathematics and
Applications).

{\ A substantial part of this work was developed during N.V. Chemetov's
visit to the NOVAMath Research Center. He would like to thank the NOVAMath
for the financial support (through the projects UIDB/00297/2020 and
UIDP/00297/2020) and the very good working conditions. }

\end{document}